\documentclass[11pt]{amsart}
\usepackage[utf8]{inputenc}
\usepackage{amssymb,amsmath,amsthm,comment,multicol,mathdots}
\usepackage[hmargin=1 in, vmargin = 1 in]{geometry}
\usepackage[nobysame]{amsrefs}
\usepackage[svgnames]{xcolor}
\usepackage{tikz,hyperref}

\usepackage{ytableau}
\usepackage[all]{xy}
\usepackage{enumerate}
\usepackage[colorinlistoftodos]{todonotes}

\usepackage[capitalise,noabbrev]{cleveref}
\crefformat{equation}{(#2#1#3)}
\Crefformat{equation}{(#2#1#3)}

\newtheorem{theorem}{Theorem}[section]
\newtheorem*{theorem*}{Theorem}
\newtheorem{lemma}[theorem]{Lemma}
\newtheorem{proposition}[theorem]{Proposition}
\newtheorem{corollary}[theorem]{Corollary}

\theoremstyle{definition}
\newtheorem{definition}[theorem]{Definition}
\newtheorem{remark}[theorem]{Remark}
\newtheorem{problem}[theorem]{Problem}

\theoremstyle{remark}
\newtheorem{examplex}[theorem]{Example}
\newenvironment{example}
{\pushQED{\qed}\examplex}
{\popQED\endexamplex}

\definecolor{darkblue}{rgb}{0.0,0,0.7}

\newcommand{\newword}[1]{\textcolor{darkblue}{\textbf{\emph{#1}}}}

\DeclareMathOperator{\demprodfpf}{\delta^{\Sp}}

\newcommand{\asm}{\mathrm{ASM}}
\newcommand{\asmsp}{\mathrm{ASM}^{\Sp}}

\newcommand{\asmsk}{\mathrm{ASM}^{\Sk}}

\newcommand{\RPP}{{\mathrm{RPP}}}
\newcommand{\SetRPP}{{\mathrm{SetRPP}}}

\newcommand{\whitecirc}{\tikz\draw[fill=white] (0,0) circle (0.2em);}

\newlength{\cellsize} 
\newsavebox{\cell}
\sbox{\cell}{\begin{picture}(18,18)
\put(0,0){\line(1,0){18}}
\put(0,0){\line(0,1){18}}
\put(18,0){\line(0,1){18}}
\put(0,18){\line(1,0){18}}
\end{picture}}
\newcommand\cellify[1]{\def\thearg{#1}\def\nothing{}%
\ifx\thearg\nothing
\vrule width0pt height\cellsize depth0pt\else
\hbox to 0pt{\usebox{\cell} \hss}\fi%
\vbox to \cellsize{
\vss
\hbox to \cellsize{\hss$#1$\hss}
\vss}}
\newcommand\tableau[1]{\vtop{\let\\\cr
\baselineskip -16000pt \lineskiplimit 16000pt \lineskip 0pt
\ialign{&\cellify{##}\cr#1\crcr}}}

\newcommand{\Z}{\mathbb{Z}}

\newcommand{\Cyc}{\mathrm{Cyc}}

\newcommand{\wt}{\mathrm{wt}}

\newcommand{\mb}{\mathcal{B}}

\newcommand{\md}{\mathcal{D}}
\newcommand{\mf}{\mathcal{F}}

\newcommand{\mr}{\mathcal{R}}

\newcommand{\mm}{\mathcal{M}}

\newcommand{\D}{\mathcal{D}}
\newcommand{\bfa}{\textbf{a}}
\newcommand{\bfb}{\textbf{b}}

\newcommand{\fkS}{\mathfrak{S}}
\newcommand{\fkSsp}{\mathfrak{S}^{\mathrm{Sp}}}

\newcommand{\fkG}{\mathfrak{G}}

\newcommand{\fkGsp}{\mathfrak{G}^{\mathrm{Sp}}}

\newcommand{\rothe}{\md}
\newcommand{\Dsp}{\rothe^{\Sp}}
\newcommand{\negset}{\mathcal{N}}
\newcommand{\negsetsp}{\mathcal{N}^{\Sp}}

\newcommand{\key}{\kappa}
\newcommand{\keysp}{\key^{\Sp}}

\newcommand{\ellfpf}{\ell^\mf}

\newcommand{\iprod}{\ \hat{\circ}\ }

\newcommand{\Sp}{\mathrm{Sp}}
\newcommand{\Sk}{\mathrm{Sk}}
\newcommand{\bpd}{\mathrm{BPD}}
\newcommand{\bpdsp}{\mathrm{BPD}^\Sp}
\newcommand{\bpdspback}{\overleftarrow{\mathrm{BPD}}^\Sp}

\newcommand{\lrtri}{\begin{tikzpicture}
	\draw (0,0) -- (0,.25) -- (.25,0) -- (0,0);
\end{tikzpicture}}

\newcommand{\maybeaddspace}{%
	\futurelet\nexttoken\checknexttoken
}

\newcommand{\checknexttoken}{%
	\ifx\nexttoken.\else
	\ifx\nexttoken,\else
	\ifx\nexttoken;\else
	\ifx\nexttoken:\else
	\ifx\nexttoken\}\else
	\;\fi
	\fi
	\fi
	\fi
	\fi
}

\newcommand{\inlinescale}{.3}

\newcommand{\cross}{\begin{tikzpicture}[scale = \inlinescale]
	\draw (0,0) -- (1,0) -- (1,1) -- (0,1) -- (0,0);
	\draw[thick,ForestGreen] (.5,0) -- (.5,1);
	\draw[thick,ForestGreen] (0,.5) -- (1,.5);
\end{tikzpicture}}
\newcommand{\seelbow}{\begin{tikzpicture}[scale = \inlinescale]
	\draw (0,0) -- (1,0) -- (1,1) -- (0,1) -- (0,0);
	\draw[thick,ForestGreen] (.5,0) -- (.5,.5) -- (1,.5);
	\draw[fill=black] (.5,.5) circle (0.5em);
\end{tikzpicture}}
\newcommand{\nwelbow}{\begin{tikzpicture}[scale = \inlinescale]
	\draw (0,0) -- (1,0) -- (1,1) -- (0,1) -- (0,0);
	\draw[thick,ForestGreen] (0,.5) -- (.5,.5) -- (.5,1);
    \draw[fill=white] (.5,.5) circle (0.5em);
\end{tikzpicture}}
\newcommand{\hwire}{\begin{tikzpicture}[scale = \inlinescale]
	\draw (0,0) -- (1,0) -- (1,1) -- (0,1) -- (0,0);
	\draw[thick,ForestGreen] (0,.5) -- (1,.5);
\end{tikzpicture}}
\newcommand{\bump}{\begin{tikzpicture}[scale = \inlinescale]
	\draw (0,0) -- (1,0) -- (1,1) -- (0,1) -- (0,0);
	\draw[thick,ForestGreen] (0,.5) to [bend right=45] (.5,1);
    \draw[thick,ForestGreen] (.5,0) to [bend left=45] (1,.5);
\end{tikzpicture}}
\newcommand{\vwire}{\begin{tikzpicture}[scale = \inlinescale]
	\draw (0,0) -- (1,0) -- (1,1) -- (0,1) -- (0,0);
	\draw[thick,ForestGreen] (.5,0) -- (.5,1);
\end{tikzpicture}}
\newcommand{\nowire}{\begin{tikzpicture}[scale = \inlinescale]
	\draw (0,0) -- (1,0) -- (1,1) -- (0,1) -- (0,0);
\end{tikzpicture}}
\newcommand{\notri}{\begin{tikzpicture}[scale = \inlinescale]
	\draw (0,0) -- (1,0) -- (0,1) -- (0,0);
\end{tikzpicture}}
\newcommand{\yestri}{\begin{tikzpicture}[scale = \inlinescale]
	\draw (0,0) -- (1,0) -- (0,1) -- (0,0);
	\draw[thick,ForestGreen] (.5,0) -- (.5,.5) -- (0,.5);
\end{tikzpicture}}

\begin{document}

\title{Symplectic bumpless pipe dreams}

\author{Zachary Hamaker}
\address[ZH]{Department of Mathematics, University of Florida, Gainesville FL 32611}
\email{zhamaker@ufl.edu}

\author{Anna Weigandt}
\address[AW]{School of Mathematics, University of Minnesota, Minneapolis MN 55455}
\email{weigandt@umn.edu}

\thanks{ZH was partially funded by NSF grant DMS-2054423. AW was partially funded by NSF grant DMS-2344764.}

\begin{abstract}
    Bumpless pipe dreams are combinatorial objects in bijection with alternating sign matrices.
    Double Grothendieck polynomials can be expressed as weighted sums over bumpless pipe dreams, a perspective that has led to many advances in the combinatorics of Schubert calculus and connections with integrable systems.
    We introduce symplectic bumpless pipe dreams and present an analogous formula for symplectic Grothendieck polynomials, which represent the equivariant $K$-theoretic classes of $\mathrm{Sp}_{2n}$-orbit closures in the flag variety.
    To prove this result, we construct an analogue of Lascoux's inflation on these objects, providing a combinatorial realization of the symplectic transition equations that appear in work of the first author, Marberg, and Pawlowski.
    We also show that certain marked symplectic bumpless pipe dreams correspond to set-valued reverse plane partitions.
    This allows us to give a tableau formula for a class of symplectic Grothendieck polynomials and to compute principal specializations of the corresponding symplectic Schubert polynomials, answering a question of Hamaker--Marberg--Pawlowski in the fixed-point-free case.
\end{abstract}

\maketitle

\vspace{-1em}

\section{Introduction}

The \emph{double Grothendieck polynomials} $\fkG_w$, introduced by Lascoux and Sch\"utzenberger in~\cite{lascoux1982structure}, represent equivariant $K$-theory classes of Schubert varieties in the flag variety.
Originally defined in terms of divided difference operators, they now have several combinatorial descriptions, such as the pipe dream formula~\cites{bergeron1993rc,fomin1996yang,FK-1,knutson2005grobner,km04}.
In~\cite{weigandt2021bumpless}, the second author showed how to interpret Lascoux's \emph{alternating sign matrix (ASM)} formula for $\fkG_w$~\cite{lascoux2002chern} as a $K$-theoretic extension of the \emph{bumpless pipe dream (BPD)} formula for Schubert polynomials in~\cite{lam2021back}.
Since then, bumpless pipe dreams have been studied extensively.
A few of the many advances include diagonal degenerations of matrix Schubert varieties~\cites{HPW22, klein2025bumpless, Klein23}, back stable Schubert calculus formulas~\cite{lam2023back}, formulas for certain Schubert structure coefficients~\cite{huang2023schubert}, integrable systems interpretations of $\fkG_w$~\cites{buciumas2022double,brubaker2020frozen}, and bijections from pipe dreams to bumpless pipe dreams~\cites{gao2021canonical, huang2024marked, Knutson.Udell2023}.

In this paper, we study the \emph{symplectic Grothendieck polynomials} $\fkGsp_\tau$, which were introduced in~\cite{wyser2017polynomials} to represent the equivariant $K$-theory classes of orbit closures of the symplectic group $\Sp_{2n}$ acting on the flag variety.
These polynomials, originally defined in terms of divided difference operators, are indexed by fixed-point-free involutions, or \emph{$\mf$-involutions} for short.
They also have descriptions as sums of Grothendieck polynomials~\cite{marberg2020k}*{Thm.~3.12}, recursively with symplectic Lascoux--Sch\"utzenberger (LS) transition equations~\cites{hamaker2018transition,hamaker2020fixed,marberg2021some}, and combinatorially via symplectic pipe dreams~\cite{marberg2021principal}*{Thm.~4.5}.  Symplectic Grothendieck polynomials are representatives for skew-symmetric matrix Schubert varieties~\cites{hamaker2022involution,marberg2020k}.
In this latter setting, symplectic pipe dreams index components in certain antidiagonal degenerations of skew-symmetric matrix Schubert varieties~\cite{marberg2022grobner}, giving the analogous result to~\cite{knutson2005grobner}.

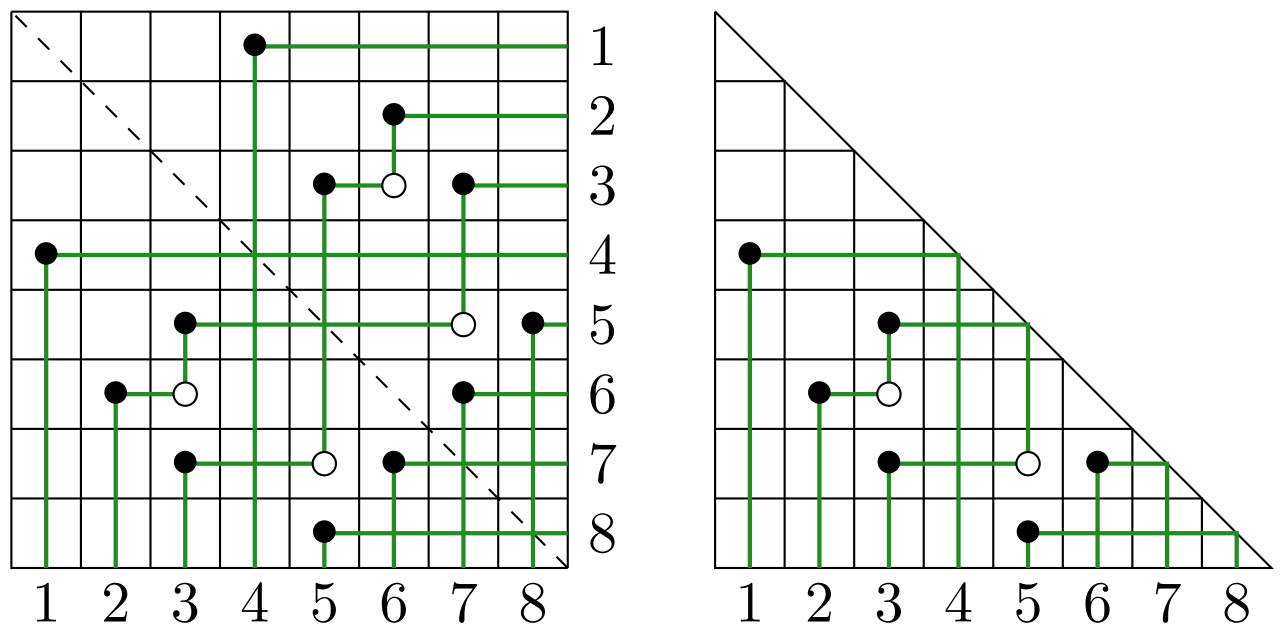
\begin{figure}
\raisebox{6em}{$\begin{bmatrix}
    0 & 0 & 0 & 1 & 0 & 0 & 0 &0\\
    0 & 0 & 0 & 0 & 0 & 1 & 0 &0\\
    0 & 0 & 0 & 0 & 1 & -1 & 1 &0\\
    1 & 0 & 0 & 0 & 0 & 0 & 0 &0\\
    0 & 0 & 1 & 0 & 0 & 0 & -1 &1\\
    0 & 1 & -1 & 0 & 0 & 0 & 1 &0\\
    0 & 0 & 1 & 0 & -1 & 1 & 0 &0\\
    0 & 0 & 0 & 0 & 1 & 0 & 0 &0
\end{bmatrix}\quad$}
\begin{tikzpicture}[scale=.46]
\draw (.5,7.5) -- (.5,-.5) -- (8.5,-.5);
\draw[dashed] (8.5,-.5) -- (.5,7.5);
\draw (.5,6.5) -- (1.5,6.5) -- (1.5,-.5);
\draw (.5,5.5) -- (2.5,5.5) -- (2.5,-.5);
\draw (.5,4.5) -- (3.5,4.5) -- (3.5,-.5);
\draw (.5,3.5) -- (4.5,3.5) -- (4.5,-.5);
\draw (.5,2.5) -- (5.5,2.5) -- (5.5,-.5);
\draw (.5,1.5) -- (6.5,1.5) -- (6.5,-.5);
\draw (.5,0.5) -- (7.5,0.5) -- (7.5,-.5);
\draw[thick,ForestGreen] (1,-.5) -- (1,4)  -- (4,4)--(4,-.5);
\draw[thick,ForestGreen] (2,-.5) --(2,2)--(3,2)--(3,3)--(5,3)--(5,1)--(3,1)--(3,-.5);
\draw[thick,ForestGreen] (5,-.5) --(5,0)--(8,0)--(8,-.5);
\draw[thick,ForestGreen] (6,-.5) --(6,1)--(7,1)--(7,-.5);
\node at (1,4) {$\bullet$};
\node at (2,2) {$\bullet$};
\node at (3,1) {$\bullet$};
\node at (3,2) {$\whitecirc$};
\node at (5,1) {$\whitecirc$};
\node at (3,3) {$\bullet$};
\node at (5,0) {$\bullet$};
\node at (6,1) {$\bullet$};
\node at (1,-1) {1};
\node at (2,-1) {2};
\node at (3,-1) {3};
\node at (4,-1) {4};
\node at (5,-1) {5};
\node at (6,-1) {6};
\node at (7,-1) {7};
\node at (8,-1) {8};
\begin{scope}[cm={0,-1,-1,0,(8,8)}]
\draw (.5,7.5) -- (.5,-.5) -- (8.5,-.5);
\draw (.5,6.5) -- (1.5,6.5) -- (1.5,-.5);
\draw (.5,5.5) -- (2.5,5.5) -- (2.5,-.5);
\draw (.5,4.5) -- (3.5,4.5) -- (3.5,-.5);
\draw (.5,3.5) -- (4.5,3.5) -- (4.5,-.5);
\draw (.5,2.5) -- (5.5,2.5) -- (5.5,-.5);
\draw (.5,1.5) -- (6.5,1.5) -- (6.5,-.5);
\draw (.5,0.5) -- (7.5,0.5) -- (7.5,-.5);
\draw[thick,ForestGreen] (1,-.5) -- (1,4)  -- (4,4)--(4,-.5);
\draw[thick,ForestGreen] (2,-.5) --(2,2)--(3,2)--(3,3)--(5,3)--(5,1)--(3,1)--(3,-.5);
\draw[thick,ForestGreen] (5,-.5) --(5,0)--(8,0)--(8,-.5);
\draw[thick,ForestGreen] (6,-.5) --(6,1)--(7,1)--(7,-.5);
\node at (1,4) {$\bullet$};
\node at (2,2) {$\bullet$};
\node at (3,1) {$\bullet$};
\node at (3,2) {$\whitecirc$};
\node at (5,1) {$\whitecirc$};
\node at (3,3) {$\bullet$};
\node at (5,0) {$\bullet$};
\node at (6,1) {$\bullet$};
\node at (1,-1) {1};
\node at (2,-1) {2};
\node at (3,-1) {3};
\node at (4,-1) {4};
\node at (5,-1) {5};
\node at (6,-1) {6};
\node at (7,-1) {7};
\node at (8,-1) {8};
\end{scope}
\end{tikzpicture}
\quad 
\begin{tikzpicture}[scale=.46]
\draw (.5,7.5) -- (.5,-.5) -- (8.5,-.5) -- (.5,7.5);
\draw (.5,6.5) -- (1.5,6.5) -- (1.5,-.5);
\draw (.5,5.5) -- (2.5,5.5) -- (2.5,-.5);
\draw (.5,4.5) -- (3.5,4.5) -- (3.5,-.5);
\draw (.5,3.5) -- (4.5,3.5) -- (4.5,-.5);
\draw (.5,2.5) -- (5.5,2.5) -- (5.5,-.5);
\draw (.5,1.5) -- (6.5,1.5) -- (6.5,-.5);
\draw (.5,0.5) -- (7.5,0.5) -- (7.5,-.5);
\draw[thick,ForestGreen] (1,-.5) -- (1,4)  -- (4,4)--(4,-.5);
\draw[thick,ForestGreen] (2,-.5) --(2,2)--(3,2)--(3,3)--(5,3)--(5,1)--(3,1)--(3,-.5);
\draw[thick,ForestGreen] (5,-.5) --(5,0)--(8,0)--(8,-.5);
\draw[thick,ForestGreen] (6,-.5) --(6,1)--(7,1)--(7,-.5);
\node at (1,4) {$\bullet$};
\node at (2,2) {$\bullet$};
\node at (3,1) {$\bullet$};
\node at (3,2) {$\whitecirc$};
\node at (5,1) {$\whitecirc$};
\node at (3,3) {$\bullet$};
\node at (5,0) {$\bullet$};
\node at (6,1) {$\bullet$};
\node at (1,-1) {1};
\node at (2,-1) {2};
\node at (3,-1) {3};
\node at (4,-1) {4};
\node at (5,-1) {5};
\node at (6,-1) {6};
\node at (7,-1) {7};
\node at (8,-1) {8};
\end{tikzpicture} 
\caption{A symplectic alternating sign matrix with corresponding bumpless pipe dream and symplectic bumpless pipe dream for $\tau = (1,4)(2,3)(5,8)(6,7)$.}
\label{fig:fold-bpd}
\end{figure}
We introduce \emph{symplectic bumpless pipe dreams}, a novel family of combinatorial objects.
These can be viewed as half of a symmetric bumpless pipe dream satisfying the additional conditions in \cref{subsection:sympBPD}, and thus they correspond to certain symmetric alternating sign matrices; see \cref{fig:fold-bpd}.
Each symplectic bumpless pipe dream has an associated $\mf$-involution $\tau$, obtained by reading off a perfect matching of its pipes.
Let $\bpdsp(\tau)$ denote the set of symplectic bumpless pipe dreams associated to $\tau$.
Our first main result is:
\begin{theorem*}[\cref{c:sbpd-formula}]
For $\tau$ an $\mf$-involution and $\wt$ an appropriate weight function,
\[
\fkGsp_\tau = \beta^{-\ellfpf(\tau)}\sum_{B \in \bpdsp(\tau)} \wt(B).
\]
\end{theorem*}
We give an example of the main theorem below.
\begin{example}
    Let $\tau=(1,4)(2,3)(5,7)(6,8)$.  Pictured below are the elements of $\bpdsp(\tau)$.
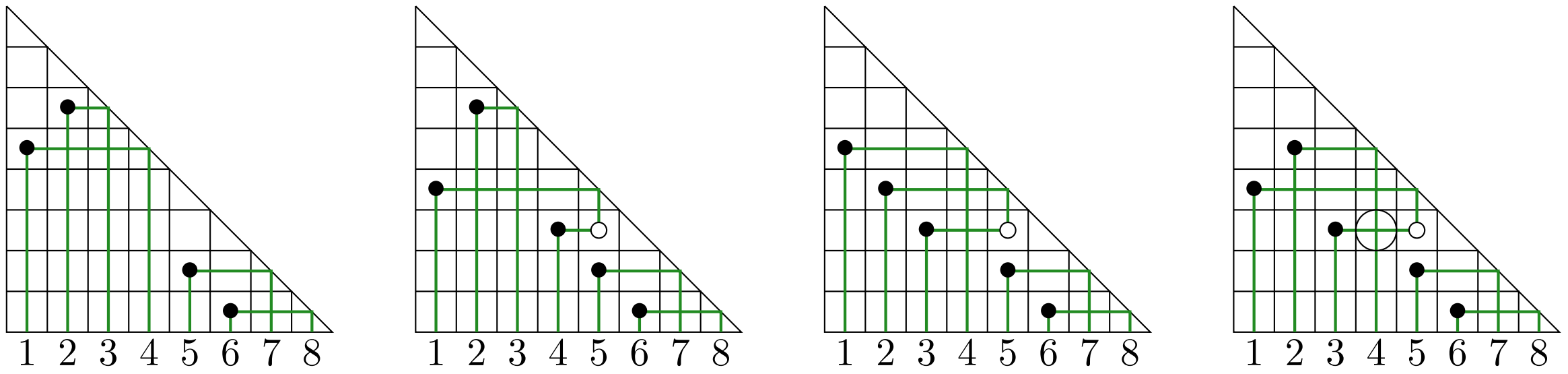
\[	\begin{tikzpicture}[scale=.4]
\draw (.5,7.5) -- (.5,-.5) -- (8.5,-.5) -- (.5,7.5);
\draw (.5,6.5) -- (1.5,6.5) -- (1.5,-.5);
\draw (.5,5.5) -- (2.5,5.5) -- (2.5,-.5);
\draw (.5,4.5) -- (3.5,4.5) -- (3.5,-.5);
\draw (.5,3.5) -- (4.5,3.5) -- (4.5,-.5);
\draw (.5,2.5) -- (5.5,2.5) -- (5.5,-.5);
\draw (.5,1.5) -- (6.5,1.5) -- (6.5,-.5);
\draw (.5,0.5) -- (7.5,0.5) -- (7.5,-.5);
\draw[thick,ForestGreen] (1,-.5) -- (1,4) -- (4,4) -- (4,-.5);
\draw[thick,ForestGreen] (2,-.5) -- (2,5) -- (3,5) -- (3,-.5);
\draw[thick,ForestGreen] (6,-.5) -- (6,0) -- (8,0) -- (8,-.5);
\draw[thick,ForestGreen] (5,-.5) -- (5,1) -- (7,1) -- (7,-.5);
\node at (1,4) {$\bullet$};
\node at (2,5) {$\bullet$};
\node at (5,1) {$\bullet$};
\node at (6,0) {$\bullet$};
\node at (1,-1) {1};
\node at (2,-1) {2};
\node at (3,-1) {3};
\node at (4,-1) {4};
\node at (5,-1) {5};
\node at (6,-1) {6};
\node at (7,-1) {7};
\node at (8,-1) {8};
\end{tikzpicture} 
\hspace{2em}
\begin{tikzpicture}[scale=.4]
\draw (.5,7.5) -- (.5,-.5) -- (8.5,-.5) -- (.5,7.5);
\draw (.5,6.5) -- (1.5,6.5) -- (1.5,-.5);
\draw (.5,5.5) -- (2.5,5.5) -- (2.5,-.5);
\draw (.5,4.5) -- (3.5,4.5) -- (3.5,-.5);
\draw (.5,3.5) -- (4.5,3.5) -- (4.5,-.5);
\draw (.5,2.5) -- (5.5,2.5) -- (5.5,-.5);
\draw (.5,1.5) -- (6.5,1.5) -- (6.5,-.5);
\draw (.5,0.5) -- (7.5,0.5) -- (7.5,-.5);
\draw[thick,ForestGreen] (1,-.5) -- (1,3) -- (5,3) -- (5,2)--(4,2)--(4,-.5);
\draw[thick,ForestGreen] (2,-.5) -- (2,5) -- (3,5) -- (3,-.5);
\draw[thick,ForestGreen] (6,-.5) -- (6,0) -- (8,0) -- (8,-.5);
\draw[thick,ForestGreen] (5,-.5) -- (5,1) -- (7,1) -- (7,-.5);
\node at (1,3) {$\bullet$};
\node at (2,5) {$\bullet$};
\node at (4,2) {$\bullet$};
\node at (5,2) {$\whitecirc$};
\node at (5,1) {$\bullet$};
\node at (6,0) {$\bullet$};
\node at (1,-1) {1};
\node at (2,-1) {2};
\node at (3,-1) {3};
\node at (4,-1) {4};
\node at (5,-1) {5};
\node at (6,-1) {6};
\node at (7,-1) {7};
\node at (8,-1) {8};
\end{tikzpicture}
\hspace{2em}
\begin{tikzpicture}[scale=.4]
\draw (.5,7.5) -- (.5,-.5) -- (8.5,-.5) -- (.5,7.5);
\draw (.5,6.5) -- (1.5,6.5) -- (1.5,-.5);
\draw (.5,5.5) -- (2.5,5.5) -- (2.5,-.5);
\draw (.5,4.5) -- (3.5,4.5) -- (3.5,-.5);
\draw (.5,3.5) -- (4.5,3.5) -- (4.5,-.5);
\draw (.5,2.5) -- (5.5,2.5) -- (5.5,-.5);
\draw (.5,1.5) -- (6.5,1.5) -- (6.5,-.5);
\draw (.5,0.5) -- (7.5,0.5) -- (7.5,-.5);
\draw[thick,ForestGreen] (1,-.5) -- (1,4) -- (4,4) -- (4,-.5);
\draw[thick,ForestGreen] (2,-.5) -- (2,3) -- (5,3) --(5,2)--(3,2)-- (3,-.5);
\draw[thick,ForestGreen] (6,-.5) -- (6,0) -- (8,0) -- (8,-.5);
\draw[thick,ForestGreen] (5,-.5) -- (5,1) -- (7,1) -- (7,-.5);
\node at (1,4) {$\bullet$};
\node at (2,3) {$\bullet$};
\node at (3,2) {$\bullet$};
\node at (5,2) {$\whitecirc$};
\node at (5,1) {$\bullet$};
\node at (6,0) {$\bullet$};
\node at (1,-1) {1};
\node at (2,-1) {2};
\node at (3,-1) {3};
\node at (4,-1) {4};
\node at (5,-1) {5};
\node at (6,-1) {6};
\node at (7,-1) {7};
\node at (8,-1) {8};
\end{tikzpicture} 
\hspace{2em}
\begin{tikzpicture}[scale=.4]
\draw (.5,7.5) -- (.5,-.5) -- (8.5,-.5) -- (.5,7.5);
\draw (.5,6.5) -- (1.5,6.5) -- (1.5,-.5);
\draw (.5,5.5) -- (2.5,5.5) -- (2.5,-.5);
\draw (.5,4.5) -- (3.5,4.5) -- (3.5,-.5);
\draw (.5,3.5) -- (4.5,3.5) -- (4.5,-.5);
\draw (.5,2.5) -- (5.5,2.5) -- (5.5,-.5);
\draw (.5,1.5) -- (6.5,1.5) -- (6.5,-.5);
\draw (.5,0.5) -- (7.5,0.5) -- (7.5,-.5);
\draw[thick,ForestGreen] (2,-.5) -- (2,4) -- (4,4) -- (4,-.5);
\draw[thick,ForestGreen] (1,-.5) -- (1,3) -- (5,3) --(5,2)--(3,2)-- (3,-.5);
\draw[thick,ForestGreen] (6,-.5) -- (6,0) -- (8,0) -- (8,-.5);
\draw[thick,ForestGreen] (5,-.5) -- (5,1) -- (7,1) -- (7,-.5);
\node at (1,3) {$\bullet$};
\node at (2,4) {$\bullet$};
\node at (3,2) {$\bullet$};
\node at (5,2) {$\whitecirc$};
\node at (5,1) {$\bullet$};
\node at (6,0) {$\bullet$};
\node at (1,-1) {1};
\node at (2,-1) {2};
\node at (3,-1) {3};
\node at (4,-1) {4};
\node at (5,-1) {5};
\node at (6,-1) {6};
\node at (7,-1) {7};
\node at (8,-1) {8};
\path [draw] (4,2) circle [radius=.5cm];
\end{tikzpicture} 
\]
Let $x_i \oplus x_j = x_i + x_j + \beta x_ix_j$.
Here $\wt(B)$ assigns a factor of $\beta (x_i \oplus x_j)$ if the $(j,i)$th cell is $\nowire$, and $1 + \beta (x_i \oplus x_j)$ if it is $\nwelbow$.  The formula includes the global factor $\beta^{-3}$, since $\ellfpf(\tau)=3$.
From the main theorem, we conclude 
    \begin{align*}
\fkGsp_\tau 
  &= (x_1 \oplus x_2)(x_1 \oplus x_3)(x_5 \oplus x_6) \\
  &\quad + (x_1 \oplus x_2)(x_1 \oplus x_3)(x_1 \oplus x_4)\bigl(1 + \beta(x_5 \oplus x_6)\bigr) \\
  &\quad + (x_1 \oplus x_2)(x_1 \oplus x_3)(x_2 \oplus x_3)\bigl(1 + \beta(x_5 \oplus x_6)\bigr) \\
  &\quad + \beta(x_1 \oplus x_2)(x_1 \oplus x_3)(x_1 \oplus x_4)(x_2 \oplus x_3)\bigl(1 + \beta(x_5 \oplus x_6)\bigr).
\end{align*}
The first three symplectic bumpless pipe dreams are \emph{reduced}, and one may readily deduce that they belong to $\bpdsp(\tau)$ by tracing along pipes starting at the bottom of the grid and recording which pairs of pipes meet at the diagonal.
The last is not reduced, and reading its associated matching is more subtle.
Here one treats the pipes in the circled $\cross$ as not crossing; we describe the general procedure to deduce this in \cref{subsection:sympBPD}. 
\end{example}

To prove this result, we characterize the alternating sign matrices corresponding to symplectic bumpless pipe dreams, which we call \emph{symplectic ASMs}, and define \emph{symplectic inflation} operations analogous to Lascoux's inflation on alternating sign matrices \cite{lascoux2002chern}.
These operations are combinatorial realizations of a subset of symplectic LS transition equations called tree equations, and decompose into inverse \emph{symplectic droop moves} for the corresponding symplectic bumpless pipe dreams.
These recurrences determine the polynomials $\fkGsp_\tau$, so the result follows.
In fact, one can identify symplectic bumpless pipe dreams with paths in the symplectic LS tree.
This proof strategy mirrors the second author's proof in~\cite{weigandt2021bumpless}, which in turn follows an outline of Lascoux~\cite{lascoux2002chern}.
As a byproduct, we show that $\bpdsp(\tau)$ is connected by droop moves, symplectic droop moves, and their $K$-theoretic analogues; see \cref{c:droop-graph} and \cref{f:bpd,f:bpd-all-moves}.

Extending work of Buch~\cite{buch2002} and Wachs~\cite{wachs1985flagged}, Knutson--Miller--Yong~\cite{knutson2009grobner} showed that double Grothendieck polynomials for vexillary permutations have a set-valued tableau formula.
The second author showed there is a transparent bijection between such tableaux and bumpless pipe dreams for vexillary permutations with certain $\nwelbow$ tiles marked~\cite{weigandt2021bumpless}.
Our other main result is an analogous correspondence between certain set-valued shifted reverse plane partitions $\SetRPP^E(\tau)$ and marked symplectic bumpless pipe dreams for a subset of the $\mf$-vexillary involutions of~\cite{hamaker2020fixed}, which we call \emph{flagged}; see \cref{p:flagged-vexillary}.

\begin{theorem*}[\cref{t:flagged-setrpp}]
    For a flagged $\mf$-involution $\tau$, with $\wt_\SetRPP$ defined below,
    \[
    \fkGsp_\tau = \beta^{-\ellfpf(\tau)}\sum_{P \in \SetRPP^E(\tau)} \wt_\SetRPP(P).
    \]
\end{theorem*}

All of our results specialize to \emph{symplectic Schubert polynomials}; see \cref{c:schubert-bpd} for the BPD formula.
In particular, \cref{t:flagged-rpp}, combined with~\cite{okada2021intermediate}*{Thm.~5.1}, essentially solves the $\mf$ case of Problem~6.4 from~\cite{hamaker2022involution}, computing the principal specialization of a symplectic Schubert polynomial for any flagged $\mf$-involution in terms of reverse plane partitions; see \cref{c:rpp_count}.

\begin{remark}
    \label{rem:orthogonal}
    There is a parallel story for $\mathrm{O}_n$-orbit closures of the flag variety, which are indexed by involutions and have as equivariant cohomology representatives $2^{c_2(\tau)}\hat{\fkS}_\tau$, where $\hat\fkS_\tau$ is the \emph{involution Schubert polynomial} for $\tau$ and $c_2(\tau)$ is the number of $2$-cycles in $\tau$.
    These satisfy orthogonal transition equations~\cite{hamaker2018transition}*{Thm.~3.30}, with reduced orthogonal bumpless pipe dreams corresponding to these equations introduced in~\cite{marberg2026bumpless}.
    One reason we originally chose to restrict our focus to the symplectic case is that, while \emph{involution Grothendieck polynomials} are defined in~\cite{marberg2021principal}*{\S4}, their transition equations, though expected to be tractable, are not written down explicitly.
    Moreover, these polynomials do \emph{not} represent $K$-theory classes of $\mathrm{O}_n$-orbit closures, nor is it known how to build such representatives using transition recurrences, except in the vexillary case (see~\cite{marberg2027on}*{Thm.~4.5}).
    That said, orthogonal bumpless pipe dreams should correspond directly to symmetric alternating sign matrices \cite{kuperberg2002symmetry}, and hence would be of great interest. 
\end{remark}

\subsection*{Organization}
In \cref{section:fpfasmatchings}, we recall necessary background for $\mf$-involutions, including standard properties of their Bruhat order.
\cref{section:transitionformatchings} gives a precise definition of $\fkGsp_\tau$ and introduces the symplectic Lascoux--Sch\"utzenberger transition equations.
Symplectic alternating sign matrices are introduced in \cref{section:symplecticASM}, where we develop their basic properties and define symplectic inflation operations to prove a first combinatorial formula for $\fkGsp_\tau$.
In \cref{section:symplecticbpd}, we introduce symplectic bumpless pipe dreams and show they are in bijection with symplectic ASMs.
We use this correspondence to prove our first main result, as well as some additional properties of symplectic ASMs.
\cref{section:flagged} develops the connection to reverse plane partitions, including our second main result.
We conclude in \cref{section:future} with some final remarks and future directions.

\subsection*{Acknowledgements}
The authors would like to thank Valentin Buciumas, Patricia Klein, Brendan Pawlowski, and Travis Scrimshaw
for helpful conversations.
We are especially grateful to Eric Marberg for pointing out a mistake in our original description of symplectic droop moves and for comments on a later draft of this manuscript.

\subsection*{AI Statement}
Most of the work in this paper was completed before LLMs became useful for mathematics.
However, we have used LLMs (Claude and ChatGPT) extensively for copy-editing, pre-refereeing, and generating figures.
In particular, several proofs were revised with input from LLMs, though the overall structure of the arguments was unchanged.

\section{Fixed-point-free involutions as matchings}
\label{section:fpfasmatchings}

Write $\mathbb N=\{0,1,2,\ldots\}$ and $[n] := \{1,2,\dots,n\}$. Write $\ell(w)$ for the usual Coxeter length of $w$.
	Let $\mf_{2n}$ denote the set of fixed-point-free involutions in $S_{2n}$, or equivalently, the set of perfect matchings of $[2n]$.
We describe $\tau \in \mf_{2n}$ by its set of two-cycles, which we denote $\Cyc(\tau)$.
For example, $\tau = (1,4)(2,6)(3,5) \in \mf_{6}$ and $\Cyc(\tau) = \{(4,1),(6,2),(5,3)\}$.
We will always write elements of $\Cyc(\tau)$ in descending order.
Two cycles $(j,i)$ and $(\ell,k)$ are \newword{crossing} if $i < k < j < \ell$.
The \newword{crossing number} $c(\tau)$ of $\tau$ is the number of crossing pairs in $\Cyc(\tau)$.

Let $\lrtri_{2n} = \{(i,j) \in [2n] \times [2n]: i > j\}$, with the partial order $\preceq$ where $(i,j) \preceq (k,\ell)$ if $i\leq k$ and $j \leq \ell$; we write $(i,j) \prec (k,\ell)$ if in addition $(i,j) \neq (k,\ell)$.
For $p \preceq q$, we write $[p,q]_\preceq = \{r : p \preceq r \preceq q\}$ for the corresponding interval.
We use the same order and the same interval notation for positions in an $n \times n$ grid.
The \newword{Rothe diagram} of a permutation $w \in S_n$ is
\[
\rothe_w := \{(i,j) \in [n] \times [n]: w(i) > j, w^{-1}(j) > i\}.
\]
Similarly, the \newword{diagram} of the matching $\tau \in \mf_{2n}$ is
\[
\Dsp_\tau := \rothe_\tau \cap \lrtri_{2n} = \{(i,j) \in \lrtri_{2n}: \tau(i) > j, \tau(j) > i\}.
\]
For example, with $\tau = (1,4)(2,6)(3,5)$ as above, we have $\Dsp_\tau = \{(2,1),(3,1),(3,2),(5,2)\}$.
We depict $\rothe_\tau$ and $\Dsp_\tau$ below using matrix coordinates:
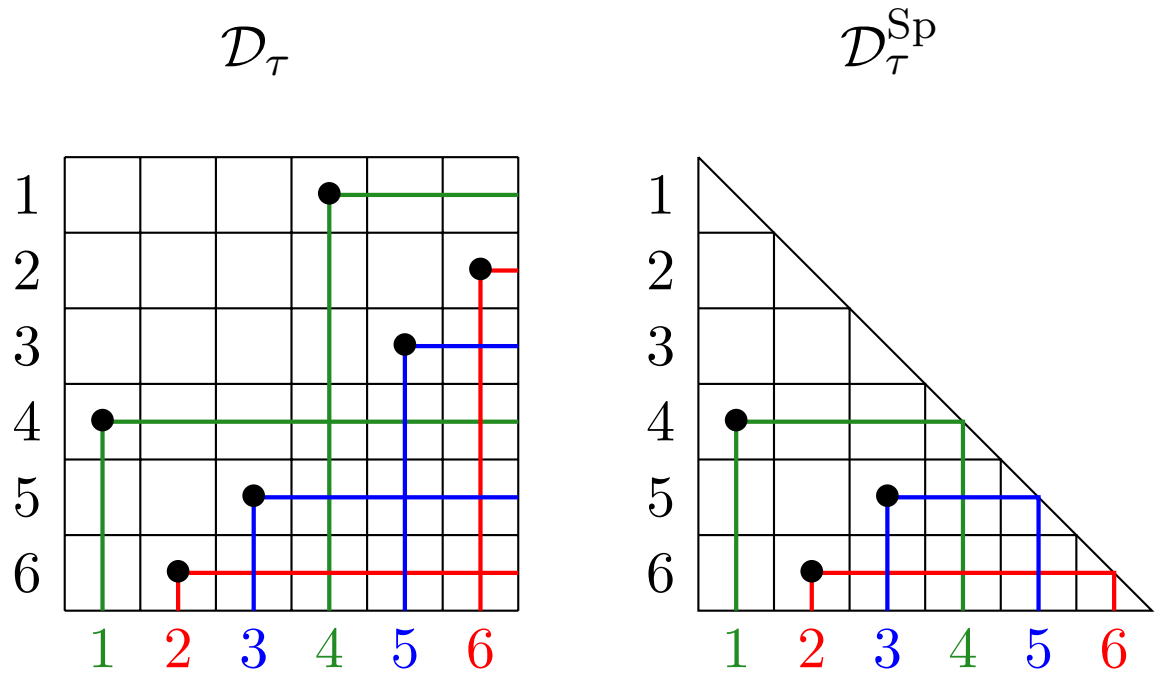
\begin{center}
\begin{tabular}{ccc}
	$\rothe_\tau$ & & $\Dsp_\tau$\\ \\
	\begin{tikzpicture}[scale=.5]
	\draw[step=1.0,black,thin] (0,0) grid (6,6);
	\node[ForestGreen] at (.5,-.5) (1) {1};
	\node[red] at (1.5,-.5) (2) {2};
	\node[blue] at (2.5,-.5) (3) {3};
	\node[ForestGreen] at (3.5,-.5) (4) {4};
	\node[blue] at (4.5,-.5) (5) {5};
	\node[red] at (5.5,-.5) (6) {6};
	\node at (-.5,5.5) {1};
	\node at (-.5,4.5) {2};
	\node at (-.5,3.5) {3};
	\node at (-.5,2.5) {4};
	\node at (-.5,1.5) {5};
	\node at (-.5,.5) {6};
	\draw[ForestGreen,thick] (.5,0) -- (.5,2.5) -- (6,2.5);
	\draw[ForestGreen,thick] (3.5,0) -- (3.5,5.5) -- (6,5.5);
	\draw[red,thick] (1.5,0) -- (1.5,.5) -- (6,.5);
	\draw[red,thick] (5.5,0) -- (5.5,4.5) -- (6,4.5);
	\draw[blue,thick] (2.5,0) -- (2.5,1.5) -- (6,1.5);
	\draw[blue,thick] (4.5,0) -- (4.5,3.5) -- (6,3.5);
	\node at (.5,2.5) (14) {$\bullet$};
	\node at (1.5,.5) (26) {$\bullet$};
	\node at (2.5,1.5) (35) {$\bullet$};
	\node at (3.5,5.5) (41) {$\bullet$};
	\node at (4.5,3.5) (62) {$\bullet$};
	\node at (5.5,4.5) (53) {$\bullet$};
\end{tikzpicture}
&&
\begin{tikzpicture}[scale=.5]
	\node[ForestGreen] at (.5,-.5) (1) {1};
	\node[red] at (1.5,-.5) (2) {2};
	\node[blue] at (2.5,-.5) (3) {3};
	\node[ForestGreen] at (3.5,-.5) (4) {4};
	\node[blue] at (4.5,-.5) (5) {5};
	\node[red] at (5.5,-.5) (6) {6};
	\node at (-.5,5.5) {1};
	\node at (-.5,4.5) {2};
	\node at (-.5,3.5) {3};
	\node at (-.5,2.5) {4};
	\node at (-.5,1.5) {5};
	\node at (-.5,.5) {6};
	\draw (0,6) -- (0,0) -- (6,0) -- (0,6);
	\draw (0,5) -- (1,5) -- (1,0);
	\draw (0,4) -- (2,4) -- (2,0);
	\draw (0,3) -- (3,3) -- (3,0);
	\draw (0,2) -- (4,2) -- (4,0);
	\draw (0,1) -- (5,1) -- (5,0);
	\draw[ForestGreen,thick] (.5,0) -- (.5,2.5) -- (3.5,2.5) -- (3.5,0);
	\draw[red,thick] (1.5,0) -- (1.5,.5) -- (5.5,.5) -- (5.5,0);
	\draw[blue,thick] (2.5,0) -- (2.5,1.5) -- (4.5,1.5) -- (4.5,0); 
    \node at (.5,2.5) (14) {$\bullet$};
	\node at (1.5,.5) (26) {$\bullet$};
	\node at (2.5,1.5) (35) {$\bullet$};
\end{tikzpicture}
\end{tabular}
\end{center}
To construct $\Dsp_\tau$ directly, place $\bullet$ in the cell $(j,i)$ for each $(j,i) \in \Cyc(\tau)$, draw rays extending below and to the right of each $\bullet$, and reflect the rightward ray down when it reaches the main diagonal.
The cells of $\lrtri_{2n}$ that do not intersect a ray form $\Dsp_\tau$.
Depending on the context, we use $\Dsp_\tau$ to refer both to the set as defined above and to this diagrammatic depiction.  Note that the visualization of $\Dsp_\tau$ is obtained from the visualization of $\rothe_\tau$ by restricting to the lower triangular part of the grid.

The \newword{shifted Young diagram} of a strict partition $\lambda = (\lambda_1 > \dots > \lambda_k)$ is the diagram $D_\lambda = \{(j,i) : i < j \leq i + \lambda_i\}$.
Note that this is the transpose of the usual convention, shifted down by one row.
A matching $\tau$ is \newword{$\mf$-dominant} if $\Dsp_\tau = D_\lambda$ for some strict partition $\lambda$.
In general, permutations associated to $\mf$-dominant matchings are \textbf{not} dominant in the usual sense.
	
The restriction of Bruhat order from permutations to fixed-point-free involutions induces a partial order $<$ on $\mf_{2n}$.
This partial order is graded, with minimal element $\Theta_{2n} = (1,2)(3,4)\cdots (2n{-1},2n)$ and rank function 
\[
\ellfpf(\tau) = |\Dsp_\tau| = \left(\sum_{(j,i) \in \Cyc(\tau)} (j - i - 1)\right) - c(\tau).
\]
For $\tau \in \mf_{2n}$, let $u_{ij}(\tau) = (i,j) \tau (i,j)$, i.e., the fixed-point-free involution obtained by conjugating $\tau$ by the transposition $(i,j)$.
Equivalently, $u_{ij}$ acts on $\Cyc(\tau)$ by swapping $i$ and $j$ if they are in different cycles, and fixes $\tau$ if $\tau(i)=j$.
The cover relations in $(\mf_{2n},<)$ are of the form 
\[
\tau \lessdot u_{ij} (\tau) \quad \mbox{where} \quad \ellfpf(\tau) + 1 = \ellfpf(u_{ij}(\tau)).
\]
These have the following explicit description, which is implicit in~\cite{incitti2004bruhat} (see also~\cite{hamaker2020fixed}*{Lem.~3.2 and Rem.~3.3}):
\begin{lemma}
	\label{l:bruhat_cover}
	Let $\tau \in \mf_{2n}$, and let $i < j$ be positive integers.
	Then $\tau \lessdot u_{ij}(\tau)$ if and only if both (1) and (2) hold:
	\begin{enumerate}
	\item (Type I) $\tau(i) < \tau(j) < i < j$ or $i < j < \tau(i) < \tau(j)$; or (Type II) $\tau(i) < i < j < \tau(j)$ or $i < \tau(i) < \tau(j) < j$;
	\item if $i < e < j$, then $\tau(e) < \tau(i)$ or $\tau(e) > \tau(j)$.
\end{enumerate}
	\end{lemma}

Type I moves correspond to the case where $(i,\tau(i))$ and $(j,\tau(j))$ are on the \emph{same} side of the diagonal in $\rothe_\tau$.  
Type II moves correspond to the case where they are on opposite sides.
They can be visualized as:
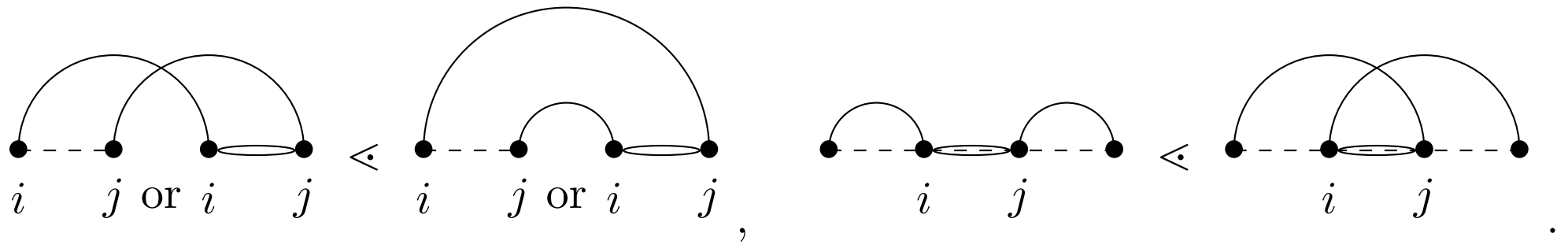
\begin{equation}
	\label{eq:cover}
	\begin{tikzpicture}[scale=.4]
	\node at (-2,-1) {$i$};
	\node at (0,-1) {$j$};
	\node at (1,-1) {or};
	\node at (2,-1) {$i$};
	\node at (4,-1) {$j$};
	\node at (0,0) {$\bullet$};
	\node at (2,0) {$\bullet$};
	\node at (-2,0) {$\bullet$};
	\node at (4,0) {$\bullet$};
	\draw (3,0) ellipse (.8cm and .1cm);
	\draw (-2,0) arc (180:0:2cm);
	\draw (0,0) arc (180:0:2cm);
	\draw[dashed] (-2,0) -- (0,0);
\end{tikzpicture}
\raisebox{1.4em}{$\ \lessdot\ $} 
	\begin{tikzpicture}[scale=.4]
	\node at (-2,-1) {$i$};
	\node at (0,-1) {$j$};
	\node at (1,-1) {or};
	\node at (2,-1) {$i$};
	\node at (4,-1) {$j$};
	\node at (0,0) {$\bullet$};
	\node at (2,0) {$\bullet$};
	\node at (-2,0) {$\bullet$};
	\node at (4,0) {$\bullet$};
	\draw (3,0) ellipse (.8cm and .1cm);
	\draw (-2,0) arc (180:0:3cm);
	\draw (0,0) arc (180:0:1cm);
	\draw[dashed] (-2,0) -- (0,0);
\end{tikzpicture}, \quad
\begin{tikzpicture}[scale=.4]
	\node at (0,-1) {$i$};
	\node at (2,-1) {$j$};
	\node at (0,0) {$\bullet$};
	\node at (2,0) {$\bullet$};
	\node at (-2,0) {$\bullet$};
	\node at (4,0) {$\bullet$};
	\draw (1,0) ellipse (.8cm and .1cm);
	\draw (-2,0) arc (180:0:1cm);
	\draw (2,0) arc (180:0:1cm);
	\draw[dashed] (-2,0) -- (4,0);
\end{tikzpicture}
\raisebox{1.4em}{$\ \lessdot\ $} 
\begin{tikzpicture}[scale=.4]
	\node at (0,-1) {$i$};
	\node at (2,-1) {$j$};
	\node at (0,0) {$\bullet$};
	\node at (2,0) {$\bullet$};
	\node at (-2,0) {$\bullet$};
	\node at (4,0) {$\bullet$};
	\draw (1,0) ellipse (.8cm and .1cm);
	\draw (-2,0) arc (180:0:2cm);
	\draw (0,0) arc (180:0:2cm);
	\draw[dashed] (-2,0) -- (4,0);
\end{tikzpicture}.
\end{equation}
The first cover relation is of Type I and the second is of Type II.
For Type II we have drawn only the case $\tau(i) < i < j < \tau(j)$; the case $i < \tau(i) < \tau(j) < j$, in which both arcs lie inside the interval $[i,j]$, is obtained by taking $i$ and $j$ to be the outer two points instead.
Here, Condition (2) in \cref{l:bruhat_cover} ensures that values in the circled regions cannot match with values in the dashed regions.
These cover relations can also be visualized on diagrams:
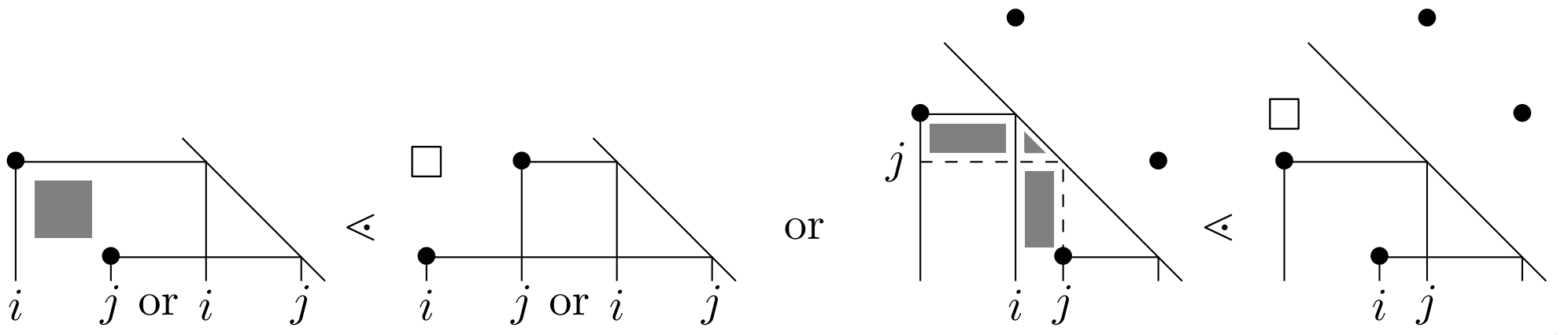
\begin{equation}
	\label{eq:dsp_cover}
\begin{tikzpicture}[scale=.4]
		\node at (0,2) {$\bullet$};
		\node at (2,0) {$\bullet$};
		\node at (0,-1) {$i$};
		\node at (2,-1) {$j$};
		\node at (3,-1) {or};
		\node at (4,-1) {$i$};
		\node at (6,-1) {$j$};
		\draw (3.5,2.5) -- (6.5,-.5);
		\draw (0,-.5) -- (0,2) -- (4,2) -- (4,-.5);
		\draw (2,-.5) -- (2,0) -- (6,0) -- (6,-.5);
		\fill [gray] (.4,.4) rectangle (1.6,1.6);
	\end{tikzpicture}
	\raisebox{2.2em}{$\ \lessdot\ $} 
	\begin{tikzpicture}[scale=.4]
		\node at (2,2) {$\bullet$};
		\node at (0,0) {$\bullet$};
		\node at (0,2) {$\square$};
		\node at (0,-1) {$i$};
		\node at (2,-1) {$j$};
		\node at (3,-1) {or};
		\node at (4,-1) {$i$};
		\node at (6,-1) {$j$};
		\draw (3.5,2.5) -- (6.5,-.5);
		\draw (0,-.5) -- (0,0) -- (6,0) -- (6,-.5);
		\draw (2,-.5) -- (2,2) -- (4,2) -- (4,-.5);
	\end{tikzpicture}
		\quad\raisebox{2.2em}{or}\quad
	\begin{tikzpicture}[scale =.4]
		\node at (0,3) {$\bullet$};
		\node at (3,0) {$\bullet$};
		\node at (2,-1) {$i$};
		\node at (3,-1) {$j$};
		\node at (-.5,2) {$j$};
		\draw (0,-.5) -- (0,3) -- (2,3) -- (2,-.5);
		\draw (3,-.5) -- (3,0) -- (5,0) -- (5,-.5);
		\draw (.5,4.5) -- (5.5,-.5);
		\draw[dashed] (3,0) -- (3,2) -- (0,2);
		\fill [gray] (.2,2.2) rectangle (1.8,2.8);
		\fill [gray] (2.2,.2) rectangle (2.8,1.8);
        \filldraw[gray] (2.2,2.2) -- (2.2,2.6) -- (2.6,2.2) -- (2.2,2.2);
        \node at (2,5) {$\bullet$};
        \node at (5,2) {$\bullet$};
	\end{tikzpicture}
	\raisebox{2.2em}{$\ \lessdot\ $} 
	\begin{tikzpicture}[scale=.4]
		\node at (0,2) {$\bullet$};
		\node at (2,0) {$\bullet$};
		\node at (2,-1) {$i$};
		\node at (3,-1) {$j$};
		\node at (0,3) {$\square$};
		\draw (0,-.5) -- (0,2) -- (3,2) -- (3,-.5);
		\draw (2,-.5) -- (2,0) -- (5,0) -- (5,-.5);
		\draw (.5,4.5) -- (5.5,-.5);
        \node at (3,5) {$\bullet$};
        \node at (5,3) {$\bullet$};
	\end{tikzpicture}.
\end{equation}
The constraints on cover relations imposed by avoiding matching dashed and circled regions are equivalent to the grey regions not containing any $\bullet$'s.
Viewed in the full diagram $\D_\tau$, a Type II move in $\Dsp_\tau$ corresponds to two Type I moves between $\bullet$'s on opposite sides of the diagonal. 
We have omitted information about other cells in the diagrams, which are not necessarily stationary.
See \cref{ex:transition} for further examples of cover relations as viewed on diagrams.

The simple transpositions $s_i = (i,i{+1})$ act on $\tau \in \mf_{2n} \cup \{\varnothing\}$ via
\[
\varnothing \iprod s_i = \varnothing \quad \mbox{and} \quad 
\tau \iprod s_i = \begin{cases}
	\varnothing & \tau(i) = i{+1},\\
	\tau & \tau > s_i \tau s_i,\\
	s_i \tau s_i & \tau < s_i \tau s_i.
\end{cases}
\]
For $\tau \in \mf_{2n} \cup\{\varnothing\}$,  $\mathbf{a}=(a_1,\ldots,a_p)$ is a \newword{symplectic Hecke word} for $\tau$ if $\tau = \Theta_{2n} \iprod s_{a_1} \iprod \cdots \iprod s_{a_p}$, where we have omitted parentheses.
In this case, we write $\demprodfpf(\mathbf{a})=\tau$ and call it the \newword{symplectic Demazure product} of $\mathbf{a}$.
The word $\mathbf{a}$ is \newword{reduced} if 
$\tau \neq \varnothing$
 and $p = \ellfpf(\tau)$.
By specializing a more general action from~\cite{rains2013deformations}, this extends to an action of the 0-Hecke monoid of $S_{2n}$, so words for $\tau\neq \varnothing$ are closed under the relations
\begin{equation}
\label{eq:hecke}
\begin{aligned}
\dots,i,i,\ldots &\leftrightarrow \dots, i,\dots, \\ 
\dots,i,j,\ldots &\leftrightarrow \dots, j,i,\ldots \qquad (|i-j|>1),\\  \dots, i, i{+1},i,\ldots &\leftrightarrow \dots,i{+1},i,i{+1},\ldots.
\end{aligned}
\end{equation}
Here $\dots$ indicates unchanged segments of the word.
As shown in~\cite{hamaker2017involution}, the words for $\tau \in \mf_{2n}$ are connected by these and the additional relation 
\begin{equation}
\label{eq:sp-hecke}
(2k,2k{+1},\dots) \leftrightarrow (2k,2k{-1},\dots).
\end{equation}

Let $\mathbf{a}$ and $\mathbf{b}$ be words.
We write $\mathbf{a}\sim \mathbf{b}$, and call $\mathbf{a}$ and $\mathbf{b}$ \newword{Hecke equivalent}, if they are related by a sequence of the relations in~\cref{eq:hecke}.
We write $\mathbf{a}\sim_\Sp \mathbf{b}$, and call $\mathbf{a}$ and $\mathbf{b}$ \newword{symplectic Hecke equivalent}, if they are related by a sequence of the relations in~\cref{eq:hecke} together with~\cref{eq:sp-hecke}.
Hecke equivalence is the stronger condition: $\mathbf{a}\sim \mathbf{b}$ implies $\mathbf{a}\sim_\Sp \mathbf{b}$, but not conversely.
By construction $\mathbf{a}\sim_\Sp \mathbf{b}$ implies $\demprodfpf(\mathbf{a}) = \demprodfpf(\mathbf{b})$, and conversely any two words for the same $\tau \in \mf_{2n}$ are symplectic Hecke equivalent~\cite{marberg2020symplectic}*{Thm.~2.4}.
Moreover, $\demprodfpf(\bfa) = \varnothing$ if and only if $\bfa \sim_\Sp \bfb$ for some $\bfb$ with $b_1$ odd~\cite{marberg2020symplectic}.

\section{Transition for matchings}
\label{section:transitionformatchings}

Associated to each $\tau \in \mf_{2n}$ is a \newword{symplectic $\beta$-Grothendieck polynomial} $\fkGsp_\tau$, also known as a fixed-point-free involution Grothendieck polynomial.
We first recall the original definition of $\fkGsp_\tau$ from~\cite{wyser2017polynomials}, which uses \emph{$K$-theoretic divided difference operators}.

Let $\beta$ be a variable.
Define the $i$th \newword{divided difference operator} and its \newword{$K$-theoretic} counterpart 
\[
\partial_i(f) = \frac{f - f(\dots,x_{i+1},x_i,\dots)}{x_i - x_{i+1}} \quad  \text{and} \quad  \partial_i^{(\beta)}(f)  = \partial_i((1+\beta x_{i+1})f).
\]
Let $x_i \oplus x_j = x_i + x_j + \beta x_ix_j$.
\begin{definition}
	\label{d:fkGsp}
For $\tau$ an $\mf$-dominant matching, define
	\begin{equation}
	\fkGsp_\tau = \prod_{(j,i) \in \Dsp_\tau} x_i \oplus x_j.
    \label{eq:dominant}
	\end{equation}
Additionally, for $u_{i,i{+}1}(\tau) \lessdot \tau$, let $\fkGsp_{u_{i,i{+}1}(\tau)} = \partial^{(\beta)}_i \fkGsp_\tau$.
These polynomials are well defined; $\fkGsp_\tau$ does not depend on the initial choice of $\mf$-dominant matching~\cite{marberg2020k}*{Thm.~3.8}, nor on the sequence of divided differences applied (see \cite{wyser2017polynomials}*{Thm.~2.2}).
\end{definition}

Symplectic Grothendieck polynomials satisfy a family of recurrence relations known as the \newword{symplectic transition equations}.
For the unrestricted transition equation we regard $\rho\in\mf_{2n}$ as a matching of the positive integers by adjoining the cycles $(2n{+}1,2n{+}2),(2n{+}3,2n{+}4),\ldots$.
Adjoining these cycles does not change $\fkGsp_\rho$.
For $k\geq1$, define
\[
I_{\Sp}(\rho,k) = \{i\geq1: \rho \lessdot u_{ik}(\rho),\ i<k\},
\qquad
R_{\Sp}(\rho,k) = \{r\geq1: \rho \lessdot u_{kr}(\rho),\ r>k\}.
\]
These sets are finite by \cref{l:bruhat_cover}, and all terms may therefore be computed in a sufficiently large finite matching set.

\begin{theorem}[{\cite{marberg2021some}*{Thm.~3.8}}]
\label{t:g_monk}
	Let $\tau \in \mf_{2n}$ and $(k,j) \in \Cyc(\tau)$.
	Then
\begin{equation}
\label{eq:g_monk}
		(1{+}\beta (x_j \oplus x_k)) \left[\sum_{\{i_1 < \dots < i_p\} \subseteq I_\Sp(\tau,j)} \beta^p\fkGsp_{u_{i_pj}\dots u_{i_1j}(\tau)} \right] = \sum_{\{r_1<\dots<r_q\} \subseteq R_\Sp(\tau,k)} \beta^q\fkGsp_{u_{kr_1}\cdots u_{kr_q}(\tau)}.
\end{equation}
\end{theorem}

\cref{t:g_monk} only holds when the $u$ operators are applied as written, from right to left.

We now explain how \cref{t:g_monk} uniquely determines $\fkGsp_\tau$ using the specific recurrence from~\cite{marberg2021some}*{Cor.~3.11}.
To do so, we must fix some notation.
Throughout, we order positions lexicographically: $(i,j)$ precedes $(k,\ell)$ if $i < k$, or if $i = k$ and $j < \ell$.
We call this the \newword{reading order}.
Note that a lexicographically maximal element of a subset of $\lrtri_{2n}$ is in particular $\preceq$-maximal.
For $\tau \in \mf_{2n}$ \emph{not} $\mf$-dominant, let $(k,j)$ be the lexicographically maximal entry in $\Dsp_\tau$.
By~\cite{hamaker2020fixed}*{Lem.~5.6}, we see $\rho = u_{k\tau(j)}(\tau)$ has the property that $R_\Sp(\rho,k) = \{\tau(j)\}$.
One can construct $\Dsp_\rho$ from $\Dsp_\tau$ by the operations on diagrams
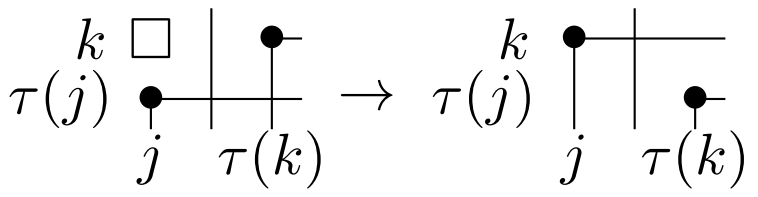
\begin{equation}
	\label{eq:ls_cover}
\begin{tikzpicture}[scale=.4]
	\node at (-1.5,0) {$\tau(j)$};
	\node at (-1,1) {$k$};
	\node at (0,-1) {$j$};
	\node at (2,-1) {$\tau(k)$};
	\node at (0,0) {$\bullet$};
	\node at (2,1) {$\bullet$};
	\node at (0,1) {$\square$};
	\draw (1,-.5) -- (1,1.5);
	\draw (0,-.5) -- (0,0) -- (2.5,0);
	\draw (2,-.5) -- (2,1) -- (2.5,1);
	\node at (3.57,0) {$\to$};
	\node at (5.5,0) {$\tau(j)$};
	\node at (6,1) {$k$};
	\node at (7,-1) {$j$};
	\node at (9,-1) {$\tau(k)$};
	\node at (9,0) {$\bullet$};
	\node at (7,1) {$\bullet$};
	\draw (8,1.5) -- (8,-.5);
	\draw (7,-.5) -- (7,1) -- (9.5,1);
	\draw (9,-.5) -- (9,0) -- (9.5,0);
\end{tikzpicture}
\quad\quad \raisebox{1.4em}{or} \quad\quad
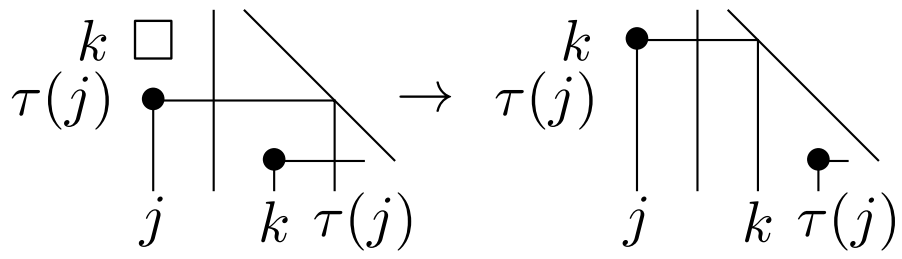
\begin{tikzpicture}[scale=.4]
	\node at (-1.5,0) {$\tau(j)$};
	\node at (-1,1) {$k$};
	\node at (0,-2) {$j$};
	\node at (2,-2) {$k$};
	\node at (3.5,-2) {$\tau(j)$};
	\node at (0,0) {$\bullet$};
	\node at (2,-1) {$\bullet$};
	\node at (0,1) {$\square$};
	\draw (1.5,1.5) -- (4,-1);
	\draw (1,1.5) -- (1,-1.5);
	\draw (0,-1.5) -- (0,0) -- (3,0) -- (3,-1.5);
	\draw (2,-1.5) -- (2,-1) -- (3.5,-1);
	\node at (4.5,0) {$\to$};
	\node at (6.5,0) {$\tau(j)$};
	\node at (7,1) {$k$};
	\node at (8,-2) {$j$};
	\node at (10,-2) {$k$};
	\node at (11.5,-2) {$\tau(j)$};
	\node at (11,-1) {$\bullet$};
	\node at (8,1) {$\bullet$};
	\draw (9.5,1.5) -- (12,-1);
	\draw (9, 1.5) -- (9,-1.5);
	\draw (8,-1.5) -- (8,1) -- (10,1) -- (10,-1.5);
	\draw (11,-1.5) -- (11,-1) -- (11.5,-1);
\end{tikzpicture},
\end{equation}
where the former case reverses a Type I cover relation and the latter reverses a Type II cover relation.
To see this, note that any $\bullet$ in the prescribed grey region from~\cref{eq:dsp_cover} would imply the existence of $(q,p) \in \Dsp_\tau$ with $(k,j) \prec (q,p)$.
Since $(k,j)$ is $\preceq$-maximal, the other cells in $\Dsp_\tau$ are unchanged by this operation, i.e., $\Dsp_\rho = \Dsp_\tau \setminus \{(k,j)\}$.
Summarizing, we have:
\begin{lemma}
\label{l:v_matching}
Let $\tau \in \mf_{2n}$ with $\Dsp_\tau \neq \varnothing$, and let $(k,j)$ be the lexicographically maximal entry in $\Dsp_\tau$.
	For $\rho = u_{k \tau(j)}(\tau)$,
	\[
	R_\Sp(\rho,k) = \{\tau(j)\} \quad \mbox{and} \quad \Dsp_\rho = \Dsp_\tau \setminus \{(k,j)\}.
	\]
\end{lemma}

The following is~\cite{marberg2021some}*{Cor.~3.11}, rearranged.

\begin{corollary}
	\label{c:g_tree}
Let $\tau \in \mf_{2n}$ with $\Dsp_\tau \neq \varnothing$, and let $(k,j)$ be lexicographically maximal in $\Dsp_\tau$.
	For $\rho = u_{k \tau(j)}(\tau)$, 
	\begin{equation}
		\fkGsp_\tau = (x_j \oplus x_k )\fkGsp_\rho + (1+\beta (x_j \oplus x_k)) \left[\sum_{\varnothing \neq \{i_1 < \dots < i_p\} \subseteq I_\Sp(\rho,j)} \beta^{p-1}\fkGsp_{u_{i_p j}\cdots u_{i_1 j}(\rho)} \right].
		\label{eq:symplectic-LS-tree}
	\end{equation}
\end{corollary}

\begin{proof}
Apply \cref{t:g_monk} to $\rho$ and the cycle $(k,j) \in \Cyc(\rho)$. Since $u_{k\tau(j)}(\rho) = \tau$, \cref{l:v_matching} shows that the right-hand side is $\fkGsp_\rho + \beta \fkGsp_\tau$.
	After subtracting $\fkGsp_\rho$ from both sides, we then divide by $\beta$ to obtain the stated equation.
\end{proof}

This equation is the \newword{symplectic LS tree transition equation}.
Here, LS stands for Lascoux--Sch\"utzenberger, who introduced the analogous recurrence for Schubert polynomials~\cite{lascoux1985schubert}.
For $\tau \in \mf_{2n}$, call each summand on the right hand side of~\cref{eq:symplectic-LS-tree} a \newword{child} of $\tau$.

\begin{lemma}
\label{l:children}
Let $\tau \in \mf_{2n}$ with lexicographically maximal entry $(k,j)$ in $\Dsp_\tau$, and let $\rho = u_{k\tau(j)}(\tau)$.
Then every diagram cell of every child of $\tau$ precedes $(k,j)$ in reading order.
\end{lemma}
\begin{proof}
Let $\sigma$ be such a matching.
By \cref{l:u_matchings}, each selected index $i$ satisfies $i<j$ and $\rho(i)<k$, so the conjugations only rearrange cycles contained in $[1,k]$.
These changes create no cells below row $k$: those rows are unchanged and every changed column $c$ has $\sigma(c)\le k$.
Moreover, $\sigma(k)$ is the least selected index, or $j$ if none is selected, so row $k$ has no cells in columns at least $j$.
Since $\Dsp_\rho=\Dsp_\tau \setminus \{(k,j)\}$ by \cref{l:v_matching}, the claim follows.
\end{proof}

The symplectic LS tree and~\cref{eq:dominant} uniquely determine each $\fkGsp_\tau$, an observation morally equivalent to \cite{hamaker2020fixed}*{Thm.~5.21}.
To see this, observe by \cref{l:children} that for $\sigma$ a child of $\tau$, every cell of $\Dsp_\sigma$ precedes $(k,j)$.
This implies $\sigma$ is closer to $\mf$-dominant than $\tau$.
Therefore, repeated applications of symplectic transition eventually terminate with all leaves corresponding to $\mf$-dominant matchings, which we can compute using~\cref{eq:dominant}.

We now translate the restrictions on cover relations to identify the elements of $I_\Sp(\rho,j)$ diagrammatically.
\begin{lemma}
\label{l:u_matchings}
	Let $\tau \in \mf_{2n}$ with $\Dsp_\tau \neq \varnothing$, and let $(k,j)$ be the lexicographically maximal entry in $\Dsp_\tau$.
	For $\rho = u_{k\tau(j)}(\tau)$, we have
	 $i \in I_{\Sp}(\rho,j)$ if and only if  either
	 \begin{itemize}
	 	\item[(I)] $j < \tau(i)$ and $(\tau(i),i)$ is the only element of $\Cyc(\tau)$ in $[(\tau(i),i),(k,j)]_\preceq$.
	 	\item[(II)] $i < j$ and $(i,\tau(i))$ is the only element of $\Cyc(\tau)$ in $[(i,\tau(i)),(j,j)]_\preceq \cup [(j,i),(k,j)]_\preceq$.
	 \end{itemize}
\end{lemma}

\begin{proof}
	Since $\rho(j) = k > j$, the Type II pattern $i < \rho(i) < \rho(j) < j$ of \cref{l:bruhat_cover} cannot occur here, and a Type I move must have $i < j < \tau(i) < k$.
	It is now easy to see that (I) corresponds to a Type I cover relation.
	The condition that $(\tau(i),i)$ is the sole $\bullet$ in the desired set is equivalent to the grey region in~\cref{eq:dsp_cover} for Type I moves being empty.
	For a Type II cover relation, we must have $\tau(i) < i < j < \tau(j)$.
	The condition that $(i,\tau(i))$ is the sole entry in the desired set is equivalent to the grey region in~\cref{eq:dsp_cover} for Type II moves being empty.
	This completes our proof.
\end{proof}

Note that Condition (I) in \cref{l:u_matchings} is equivalent to saying that $(\tau(i),i)$ is a $\preceq$-maximal $\bullet$ in $\Dsp_\tau \cap [j,k] \times [1,j]$.
This property is stronger than Condition (II), which allows for $\preceq$-greater $\bullet$'s in the rectangle $[j,k] \times [\tau(i),i]$.

\begin{figure}
\begin{center}
    \begin{tabular}{cccc}
         $\Dsp_\tau$ && $\Dsp_\rho$  \\
	\begin{tikzpicture}[scale=.5]
\draw (.5,7.5) -- (.5,-.5) -- (8.5,-.5) -- (.5,7.5);
\draw (.5,6.5) -- (1.5,6.5) -- (1.5,-.5);
\draw (.5,5.5) -- (2.5,5.5) -- (2.5,-.5);
\draw (.5,4.5) -- (3.5,4.5) -- (3.5,-.5);
\draw (.5,3.5) -- (4.5,3.5) -- (4.5,-.5);
\draw (.5,2.5) -- (5.5,2.5) -- (5.5,-.5);
\draw (.5,1.5) -- (6.5,1.5) -- (6.5,-.5);
\draw (.5,0.5) -- (7.5,0.5) -- (7.5,-.5);
\draw[thick,ForestGreen] (1,-.5) -- (1,2) -- (6,2) -- (6,-.5);
\draw[thick,ForestGreen] (2,-.5) -- (2,5) -- (3,5) -- (3,-.5);
\draw[thick,ForestGreen] (4,-.5) -- (4,0) -- (8,0) -- (8,-.5);
\draw[thick,ForestGreen] (5,-.5) -- (5,1) -- (7,1) -- (7,-.5);
\node at (1,2) {$\bullet$};
\node at (2,5) {$\bullet$};
\node at (4,0) {$\bullet$};
\node at (5,1) {$\bullet$};
\node at (1,-1) {1};
\node at (2,-1) {2};
\node at (3,-1) {3};
\node at (4,-1) {4};
\node at (5,-1) {5};
\node at (6,-1) {6};
\node at (7,-1) {7};
\node at (8,-1) {8};
\end{tikzpicture}  
&&
\begin{tikzpicture}[scale=.5]
\draw (.5,7.5) -- (.5,-.5) -- (8.5,-.5) -- (.5,7.5);
\draw (.5,6.5) -- (1.5,6.5) -- (1.5,-.5);
\draw (.5,5.5) -- (2.5,5.5) -- (2.5,-.5);
\draw (.5,4.5) -- (3.5,4.5) -- (3.5,-.5);
\draw (.5,3.5) -- (4.5,3.5) -- (4.5,-.5);
\draw (.5,2.5) -- (5.5,2.5) -- (5.5,-.5);
\draw (.5,1.5) -- (6.5,1.5) -- (6.5,-.5);
\draw (.5,0.5) -- (7.5,0.5) -- (7.5,-.5);
\draw[thick,ForestGreen] (1,-.5) -- (1,2) -- (6,2) -- (6,-.5);
\draw[thick,ForestGreen] (2,-.5) -- (2,5) -- (3,5) -- (3,-.5);
\draw[thick,ForestGreen] (5,-.5) -- (5,0) -- (8,0) -- (8,-.5);
\draw[thick,ForestGreen] (4,-.5) -- (4,1) -- (7,1) -- (7,-.5);
\node at (1,2) {$\bullet$};
\node at (2,5) {$\bullet$};
\node at (5,0) {$\bullet$};
\node at (4,1) {$\bullet$};
\node at (1,-1) {1};
\node at (2,-1) {2};
\node at (3,-1) {3};
\node at (4,-1) {4};
\node at (5,-1) {5};
\node at (6,-1) {6};
\node at (7,-1) {7};
\node at (8,-1) {8};
\end{tikzpicture}
    \end{tabular}
    \end{center}
\begin{center}
	\begin{tabular}{ccccc}
$\Dsp_{u_{14}(\rho)} $&& $\Dsp_{u_{34}(\rho)}$	&& $\Dsp_{u_{34}u_{14}(\rho)}$	\\
\begin{tikzpicture}[scale=.5]
\draw (.5,7.5) -- (.5,-.5) -- (8.5,-.5) -- (.5,7.5);
\draw (.5,6.5) -- (1.5,6.5) -- (1.5,-.5);
\draw (.5,5.5) -- (2.5,5.5) -- (2.5,-.5);
\draw (.5,4.5) -- (3.5,4.5) -- (3.5,-.5);
\draw (.5,3.5) -- (4.5,3.5) -- (4.5,-.5);
\draw (.5,2.5) -- (5.5,2.5) -- (5.5,-.5);
\draw (.5,1.5) -- (6.5,1.5) -- (6.5,-.5);
\draw (.5,0.5) -- (7.5,0.5) -- (7.5,-.5);
\draw[thick,ForestGreen] (1,-.5) -- (1,1) -- (7,1) -- (7,-.5);
\draw[thick,ForestGreen] (2,-.5) -- (2,5) -- (3,5) -- (3,-.5);
\draw[thick,ForestGreen] (5,-.5) -- (5,0) -- (8,0) -- (8,-.5);
\draw[thick,ForestGreen] (4,-.5) -- (4,2) -- (6,2) -- (6,-.5);
\node at (2,5) {$\bullet$};
\node at (4,2) {$\bullet$};
\node at (5,0) {$\bullet$};
\node at (1,1) {$\bullet$};
\node at (1,-1) {1};
\node at (2,-1) {2};
\node at (3,-1) {3};
\node at (4,-1) {4};
\node at (5,-1) {5};
\node at (6,-1) {6};
\node at (7,-1) {7};
\node at (8,-1) {8};
\end{tikzpicture}&&
\begin{tikzpicture}[scale=.5]
\draw (.5,7.5) -- (.5,-.5) -- (8.5,-.5) -- (.5,7.5);
\draw (.5,6.5) -- (1.5,6.5) -- (1.5,-.5);
\draw (.5,5.5) -- (2.5,5.5) -- (2.5,-.5);
\draw (.5,4.5) -- (3.5,4.5) -- (3.5,-.5);
\draw (.5,3.5) -- (4.5,3.5) -- (4.5,-.5);
\draw (.5,2.5) -- (5.5,2.5) -- (5.5,-.5);
\draw (.5,1.5) -- (6.5,1.5) -- (6.5,-.5);
\draw (.5,0.5) -- (7.5,0.5) -- (7.5,-.5);
\draw[thick,ForestGreen] (1,-.5) -- (1,2) -- (6,2) -- (6,-.5);

\draw[thick,ForestGreen] (2,-.5) -- (2,4) -- (4,4) -- (4,-.5);
\draw[thick,ForestGreen] (5,-.5) -- (5,0) -- (8,0) -- (8,-.5);
\draw[thick,ForestGreen] (3,-.5) -- (3,1) -- (7,1) -- (7,-.5);
\node at (1,2) {$\bullet$};
\node at (5,0) {$\bullet$};
\node at (3,1) {$\bullet$};
\node at (2,4) {$\bullet$};
\node at (1,-1) {1};
\node at (2,-1) {2};
\node at (3,-1) {3};
\node at (4,-1) {4};
\node at (5,-1) {5};
\node at (6,-1) {6};
\node at (7,-1) {7};
\node at (8,-1) {8};
\end{tikzpicture}
&&
\begin{tikzpicture}[scale=.5]
\draw (.5,7.5) -- (.5,-.5) -- (8.5,-.5) -- (.5,7.5);
\draw (.5,6.5) -- (1.5,6.5) -- (1.5,-.5);
\draw (.5,5.5) -- (2.5,5.5) -- (2.5,-.5);
\draw (.5,4.5) -- (3.5,4.5) -- (3.5,-.5);
\draw (.5,3.5) -- (4.5,3.5) -- (4.5,-.5);
\draw (.5,2.5) -- (5.5,2.5) -- (5.5,-.5);
\draw (.5,1.5) -- (6.5,1.5) -- (6.5,-.5);
\draw (.5,0.5) -- (7.5,0.5) -- (7.5,-.5);
\draw[thick,ForestGreen] (1,-.5) -- (1,1) -- (7,1) -- (7,-.5);
\draw[thick,ForestGreen] (2,-.5) -- (2,4) -- (4,4) -- (4,-.5);
\draw[thick,ForestGreen] (5,-.5) -- (5,0) -- (8,0) -- (8,-.5);
\draw[thick,ForestGreen] (3,-.5) -- (3,2) -- (6,2) -- (6,-.5);
\node at (2,4) {$\bullet$};
\node at (5,0) {$\bullet$};
\node at (1,1) {$\bullet$};
\node at (3,2) {$\bullet$};
\node at (1,-1) {1};
\node at (2,-1) {2};
\node at (3,-1) {3};
\node at (4,-1) {4};
\node at (5,-1) {5};
\node at (6,-1) {6};
\node at (7,-1) {7};
\node at (8,-1) {8};
\end{tikzpicture}
	\end{tabular}
\end{center}
\caption{Diagrams for $\tau = (1,6)(2,3)(4,8)(5,7)$ and its children in the symplectic LS tree.}
      \label{fig:transition-example}
\end{figure}
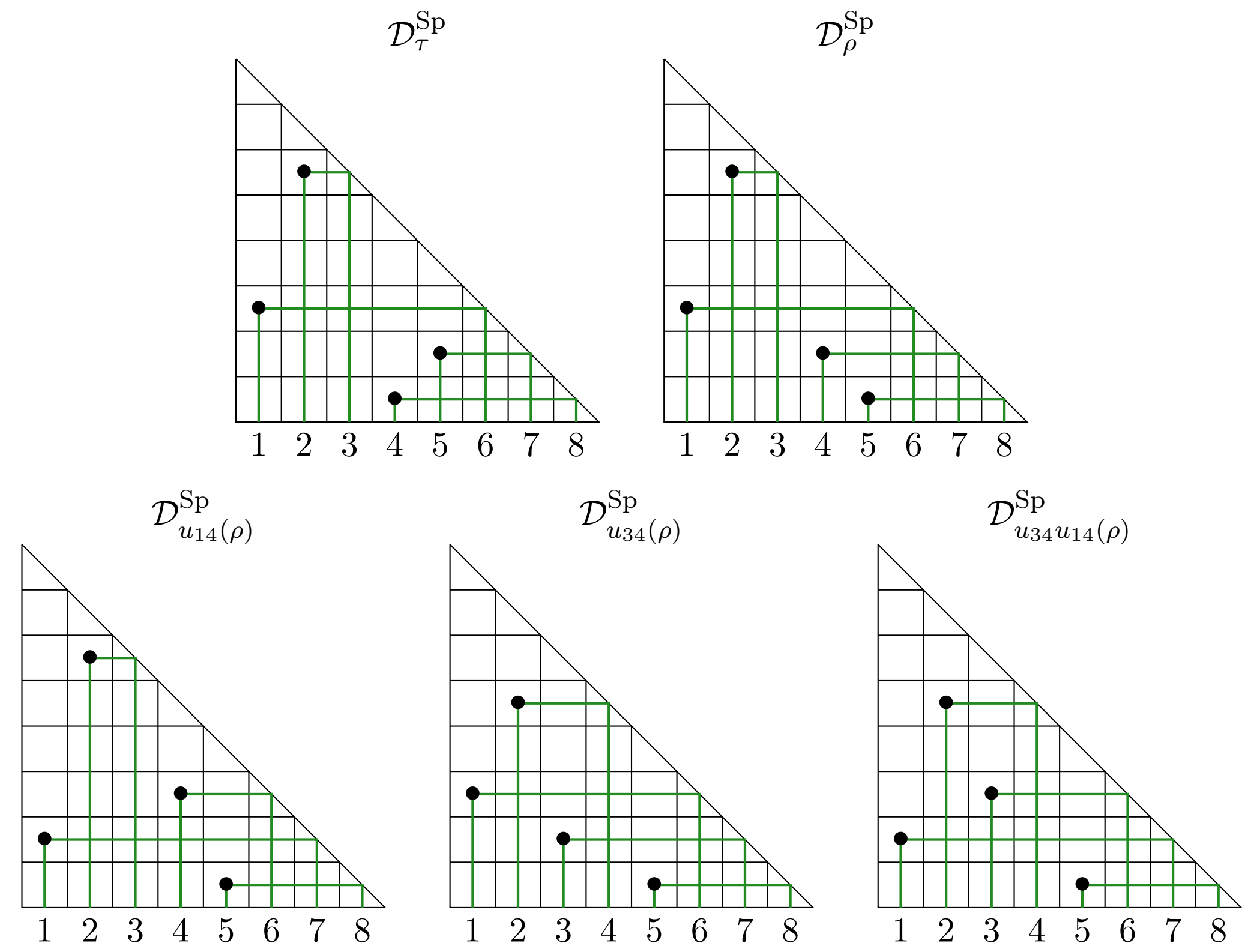

\begin{example}
	\label{ex:transition}
Consider the matching $\tau = (1,6)(2,3)(4,8)(5,7)$.
Note that $\Dsp_\tau$ has lexicographically maximal entry $(7,4)$, with $\rho = u_{7\tau(4)}(\tau) = (1,6)(2,3)(4,7)(5,8)$.
The diagrams of $\tau$ and $\rho$ appear in the first row of \cref{fig:transition-example}.
The maximal $\bullet$'s northwest of $(7,4)$ in $\Dsp_\tau$ are in positions $(6,1)$ and $(3,2)$.
The first satisfies Condition~(I) of \cref{l:u_matchings} with $i = 1$, since $\tau(1) = 6 > 4$; the second satisfies Condition~(II) with $i = 3$, since $\tau(3) = 2 < 3 < 4$.
Therefore $I_\Sp(\rho,4) = \{1,3\}$.
	Then
	\[
	\fkGsp_\tau = (x_4 \oplus x_7) \fkGsp_\rho + (1+\beta (x_4 \oplus x_7))\left(\fkGsp_{u_{14}(\rho)} + \fkGsp_{u_{34}(\rho)} + \beta \fkGsp_{u_{34} u_{14}(\rho)}\right).
	\]
As matchings, the three terms on the right-hand side are $u_{14}(\rho) = (1,7)(2,3)(4,6)(5,8)$, $u_{34}(\rho) = (1,6)(2,4)(3,7)(5,8)$, and $u_{34}u_{14}(\rho) = (1,7)(2,4)(3,6)(5,8)$.
    The diagrams for these matchings are in the second row of \cref{fig:transition-example}.
Note that $\rho \lessdot u_{14}(\rho)$ is a Type I cover relation and $\rho \lessdot u_{34}(\rho)$ is a Type II cover relation (compare with \cref{eq:dsp_cover}).
Additionally, $u_{14}u_{34}(\rho) = (1,2)(3,7)(4,6)(5,8) < u_{34}(\rho)$, so the order in which the operators are applied matters.
\end{example}

The \emph{symplectic Schubert polynomial} of $\tau \in \mf_{2n}$ is $\fkSsp_\tau = \fkGsp_\tau \big|_{\beta = 0}$.
Necessarily, the transition equations specialize as well.
For future reference, we state the specialization of~\cref{c:g_tree}, which first appeared in \cites{hamaker2018transition,hamaker2020fixed}.
\begin{corollary}
    \label{c:s-tree}
Let $\tau \in \mf_{2n}$ with $\Dsp_\tau \neq \varnothing$, and let $(k,j)$ be the lexicographically maximal entry in $\Dsp_\tau$.
For $\rho = u_{k\tau(j)}(\tau)$,
\[
\fkSsp_\tau = (x_j + x_k) \fkSsp_\rho + \sum_{i \in I_\Sp(\rho,j)} \fkSsp_{u_{ij}(\rho)}.
\]
\end{corollary}

\section{Symplectic ASMs}
\label{section:symplecticASM}

In this section, we explain how to adapt Lascoux's inflation procedure on ASMs to the symplectic setting.
As a result, we give a formula for symplectic Grothendieck polynomials in terms of certain ASMs, which we call \emph{symplectic}.

\subsection{Inflation on ASMs}

An \newword{alternating sign matrix (ASM)} is a square matrix with entries in $\{-1,0,1\}$ such that the entries in each row and column sum to $1$, and the nonzero entries alternate in sign within each row and column.  
Let $\asm(n)$ be the set of $n \times n$ alternating sign matrices.  
For our purposes, it is often useful to visualize $A\in \asm(n)$ by plotting within an $n\times n$ grid a $\bullet$ in position $(i,j)$ whenever $A_{i,j}=1$, and a $\whitecirc$ in position $(i,j)$ whenever $A_{i,j}=-1$.  
We call this the \newword{graph} of an ASM.  

\begin{figure}
$\begin{bmatrix}
0 & 1 & 0 & 0 & 0 & 0\\
0 & 0 & 0 & 0 & 1 & 0\\
1 & -1 & 0 & 1 & 0 & 0\\
0 & 0 & 1 & 0 & -1 & 1\\
0 & 1 & 0 & 0 & 0 & 0\\
0 & 0 & 0 & 0 & 1 & 0
\end{bmatrix}$
\hspace{3em}
\raisebox{-3.75em}{\begin{tikzpicture}[scale=.5]
 \draw (0,0) grid (6,6);
\draw[thick,ForestGreen] (.5,0) --(.5,3.5)--(1.5,3.5)--(1.5,3.5)--(1.5,5.5)--(6,5.5);
  \draw[thick,ForestGreen] (1.5,0) -- (1.5,1.5)--(6,1.5);
 \draw[thick,ForestGreen] (2.5,0)  -- (2.5,1.5) -- (2.5,2.5) -- (4.5,2.5) -- (4.5,4.5) -- (6,4.5);
  \draw[thick,ForestGreen] (3.5,0) -- (3.5,3.5)--(6,3.5);
  \draw[thick,ForestGreen] (4.5,0) -- (4.5,.5) -- (6,.5);
  \draw[thick,ForestGreen] (5.5,0) -- (5.5,2.5) -- (6,2.5);
\node at (1.5,5.5) {$\bullet$};
\node at (.5,3.5) {$\bullet$};
\node at (4.5,4.5) {$\bullet$};
\node at (1.5,3.5) {$\whitecirc$};
\node at (1.5,1.5) {$\bullet$};
\node at (2.5,2.5) {$\bullet$};
\node at (3.5,3.5) {$\bullet$};
\node at (4.5,2.5) {$\whitecirc$};
\node at (5.5,2.5) {$\bullet$};
\node at (4.5,.5) {$\bullet$};
\end{tikzpicture}   }
\caption{Pictured above on the left is an ASM $A$, and on the right is its Rothe diagram. We have $\rothe_A=\{(1,1),(2,1),(2,3),(2,4),(3,3),(4,2)\}$ and $\negset_A=\{(3,2),(4,5)\}$.}
    \label{fig:rothe}
\end{figure}
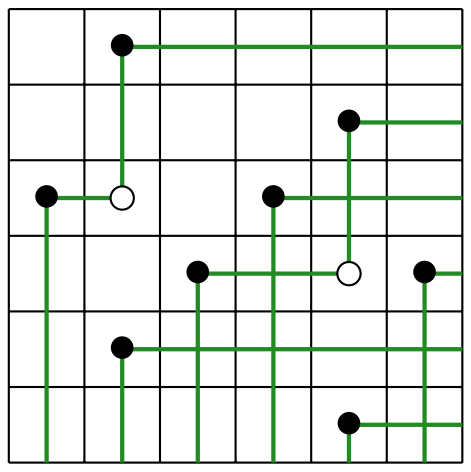

We may further decorate this picture by drawing rays that start at each $\bullet$, traveling horizontally to the right and vertically downward, and stopping whenever they encounter a $\whitecirc$ or the boundary of the $n \times n$ grid.  
We write $\negset_A$ for the set of coordinates of the $\whitecirc$'s, and $\rothe_A$ for the set of coordinates of the cells in this picture that do not contain any of the rays.  
We call both $\rothe_A$ and this visualization the \newword{Rothe diagram} of $A$.
The map $A \mapsto \rothe_A$, viewed as a visualization, places alternating sign matrices in transparent bijection with \newword{bumpless pipe dreams}: every cell of the $n \times n$ grid carries one of the six tiles $\nowire$, $\seelbow$, $\nwelbow$, $\cross$, $\hwire$, $\vwire$, with $\seelbow$ in the positions of the $1$'s of $A$, $\nwelbow$ in the positions of the $-1$'s, and $\nowire$ exactly on $\rothe_A$.
The cells of $\rothe_A$, viewed as a set, are called \newword{blank}.
See \cref{fig:rothe} for an example.

Given $(a,b)$ and an ASM $A$, the \newword{neighbors} $N_A(a,b)$ are the $1$-entries that are $\preceq$-maximal among all nonzero entries at positions $(i,j)\prec(a,b)$.
When $A_{ab} = -1$, there necessarily exist some $j< b$ and $i<a$ such that $(i,b)$ and $(a,j)$ are neighbors of $(a,b)$.  Call these the \newword{top} and \newword{left} neighbors of $(a,b)$, respectively.
For example, with $A$ as in \cref{fig:sharedneighbors}, we have $N_A(4,3) =\{(2,3), (3,2), (4,1)\}$,  with $(2,3)$ the top neighbor and $(4,1)$ the left neighbor.
Meanwhile, $N_A(4,4) =\varnothing$, since both $\preceq$-maximal nonzero entries northwest of $(4,4)$, at $(3,4)$ and $(4,3)$, are $-1$'s.

\begin{figure}
    \centering
    \begin{tikzpicture}[scale=.5]
    \draw[dashed] (4,1)--(4,5)--(2,5)--(2,4)--(1,4)--(1,2)--(0,2)--(0,1)--(4,1);
  \node at (4,5) {$\bullet$};
  \node at (2,4) {$\bullet$};
  \node at (1,2) {$\bullet$};
  \node at (0,1) {$\bullet$};
  \node at (4,1) {$\whitecirc$};
    \end{tikzpicture}
    \hspace{2em}
    \raisebox{3em}{$\mapsto$}
    \hspace{2em}
    \begin{tikzpicture}[scale=.5]
    \draw[dashed] (4,1)--(4,5)--(2,5)--(2,4)--(1,4)--(1,2)--(0,2)--(0,1)--(4,1);
  \node at (2,5) {$\bullet$};
  \node at (1,4) {$\bullet$};
  \node at (0,2) {$\bullet$};
  \node at (0,1){\,};
    \end{tikzpicture}
    \caption{An example of inflation on an ASM.}
    \label{fig:inflation}
\end{figure}
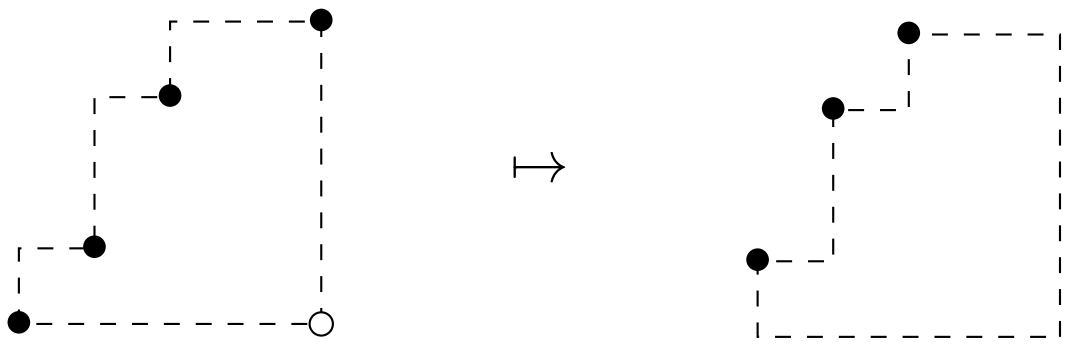

We now describe Lascoux's \cite{lascoux2002chern} \emph{inflation} operation on ASMs.
Fix $(a,b)$ so that $A_{a,b}=-1$ and $N_A(a,b) = \{(i_1,j_1),(i_2,j_2),\ldots,(i_k,j_k)\}$, ordered so that $i_1 < i_2 < \cdots < i_k = a$.
Then $(a,b)$ is \newword{removable} if $A_{ij} \neq -1$ for each $(i,j) \neq (a,b)$ in the union of intervals
\[\bigcup_{\ell=1}^{k-1}[(i_\ell,j_{\ell+1}),(a,b)]_\preceq.\]
Note that if $(i,j)$ lies in this union and $A_{ij} = 1$, then $(i,j) \in N_A(a,b)$, since otherwise the alternating condition would fail.
Given such a removable $-1$ in $A$, the \newword{inflation} operation sets $A_{a,b}=0$ and $A_{i_\ell,j_\ell}=0$ for all $\ell\in [k]$, then places $1$'s in the positions $(i_1,j_2),\ldots,(i_{k-1},j_k)$.
See \cref{fig:inflation} for a visualization of these replacements.

Note that applying inflation to a removable $-1$ produces a valid ASM which has one fewer $-1$ than the original ASM.
Every ASM with a $-1$ has a removable $-1$, namely the first $-1$ in reading order.
Therefore, we define the \newword{key} of an ASM as the permutation matrix obtained by applying inflation to removable $-1$'s until no $-1$'s remain.
By \cite{weigandt2021bumpless}*{Cor.~4.4}, the order in which one removes the removable $-1$'s does not matter.  In particular, inflation at any removable $-1$ preserves the key.
Write $\key(A)$ for the key of $A$, and for $w \in S_n$ let $\asm(w)$ be the set of ASMs with key $w$.

As we now explain, inflation is invertible.
A cell $(a,b) \in \rothe_A$ is \newword{addable} in $A$ if $N_A(a,b)$ is non-empty.
Note that because $(a,b)\in \rothe_A$, it has no neighbor to the left in row $a$, and also no neighbor above in column $b$.
For $(a,b)$ addable in $A$ and  $U = \{(i_1,j_1),\dots,(i_k,j_k)\} \subseteq N_A(a,b)$ non-empty, ordered so that $i_1 < \dots < i_k$, the \newword{deflation} operation at $(a,b)$ in $A$ for $U$ sets $A_{a,b} = -1$, $A_{a,j_k}=1$, and $A_{i_\ell,j_\ell} = 0$, $A_{i_\ell,j_{\ell-1}}=1$ for all $\ell \in [k]$, where $j_0 = b$.
Deflation and inflation at $(a,b)$ are inverse operations in the following sense.

\begin{lemma}
    \label{l:deflation}
    Let $A, A' \in \asm(n)$.
    \begin{enumerate}
        \item Assume $(a,b)$ is addable in $A$ so that $A'$ is the result of applying deflation at $(a,b)$ for a non-empty $U \subseteq N_A(a,b)$.
        Then $(a,b)$ is removable in $A'$, and applying inflation at $(a,b)$ results in $A$.
        \item Assume $(a,b)$ is removable in $A$ so that $A'$ is the result of applying inflation at $(a,b)$.
        Then $(a,b)$ is addable in $A'$, and there exists a unique $U \subseteq N_{A'}(a,b)$ so that applying deflation at $(a,b)$ for $U$ results in $A$.
    \end{enumerate}
\end{lemma}

\begin{proof}
To prove (1), first observe that the positions set equal to $1$ by deflation
    comprise all the neighbors of $(a,b)$ in $A'$.
Since elements of $U$ are neighbors of $(a,b)$ in $A$, each interval $[(i_\ell,j_\ell),(a,b)]_\preceq$ contains no other nonzero entry before deflation.
These are precisely the intervals defining removability afterward, so $(a,b)$ is removable and inflation recovers $A$.

    To prove (2), observe that the positions set equal to $1$ by inflation
    are neighbors of $(a,b)$ in $A'$.
    It is then straightforward to see that applying deflation at $(a,b)$ in $A'$ with $U$ consisting of the newly inserted positions results in $A$.
    Deflation for any other set $U'\subseteq N_{A'}(a,b)$ removes a different set of $1$'s, so uniqueness follows.
\end{proof}

\subsection{Inflation on skew ASMs}

From now on, we will focus on ASMs that are symmetric ($A=A^T$) and satisfy the additional condition that all entries on the main diagonal are $0$.  
We call such ASMs \newword{skew}.
Let $\asmsk(n)$ be the set of skew ASMs of size $n$.
For $A \in \asmsk(n)$, we call the entry at $(b,a)$ the \newword{partner} of the entry at $(a,b)$.

\begin{figure}
\centering
$\begin{bmatrix}
0 & 0 & 0 & 1 & 0 & 0 \\
0 & 0 & 1 & 0 & 0 & 0 \\
0 & 1 & 0 & -1 & 1 & 0 \\
1 & 0 & -1 & 0 & 0 & 1 \\
0 & 0 & 1 & 0 & 0 & 0 \\
0 & 0 & 0 & 1 & 0 & 0
\end{bmatrix}$
    \caption{Pictured above is a skew ASM $A$.
    There are removable $-1$'s at positions $(4,3)$ and $(3,4)$.  
    These entries share the neighbors $(2,3)$ and $(3,2)$, and thus $\keysp(A) = \varnothing$.}
    \label{fig:sharedneighbors}
\end{figure}
Our goal is to introduce a symplectic analogue of inflation on ASMs, with the corresponding symplectic key $\keysp$ defined as a composition of symplectic inflation moves.
We will show that there is one type of obstruction for symplectic inflation moves: when a removable $-1$ and its partner in $A$ share common neighbors (see \cref{fig:sharedneighbors} for an example).
When such a move is required, we say that $A$ has no symplectic key and write $\keysp(A) = \varnothing$.
Otherwise, we will show that the symplectic key agrees with the usual key defined by Lascoux.

Let $A\in \asmsk(2n)$, $a < b$, and $A_{b,a}=-1$.
We say that $(b,a)$ is \newword{$\Sp$-removable} if $(b,a)$ and $(a,b)$ are removable and have no common neighbors.
Note that by symmetry, if $(b,a)$ is a removable $-1$ in $A$, then so is $(a,b)$.
We call $(i,j) \in N_A(b,a)$ a \newword{lower neighbor} if $i > j$ and an \newword{upper neighbor} if $i < j$.

\begin{lemma}
\label{l:symplecticinflation}
Let $A \in \asmsk(2n)$, and suppose that $(b,a)$ is $\Sp$-removable in $A$.
After inflation at $(a,b)$, $(b,a)$ continues to be removable with identical neighbor set.
In particular, inflation at $(a,b)$, followed by inflation at $(b,a)$, modifies disjoint entries.
\end{lemma}

\begin{proof}
Since the neighbor sets $N_A(b,a)$ and $N_A(a,b)$ are disjoint, the union of intervals modified by inflation at $(b,a)$ does not intersect $N_A(a,b)$. Furthermore, inflation at $(b,a)$ can not introduce a new neighbor for $(a,b)$, as this would have implied the existence of a shared neighbor.  The result follows.
\end{proof}

For $A \in \asmsk(2n)$ such that $(b,a)$ is $\Sp$-removable in $A$, \newword{symplectic inflation} at $(b,a)$ is defined to be the result of applying inflation at $(a,b)$ and $(b,a)$, which, by \cref{l:symplecticinflation}, may be viewed as simultaneous.
The resulting ASM is skew: it is symmetric, and a modified diagonal entry would by symmetry be modified by both inflations.

By symmetry, for $(b,a)$ addable in $A\in\asmsk(2n)$ and
$\varnothing\ne U\subseteq N_A(b,a)$, consider the deflations at $(b,a)$ for $U$ and at $(a,b)$ for $U^T$.
If the resulting ASM is skew and
 $(b,a)$ is $\Sp$-removable, we call this operation
\newword{symplectic deflation}.
\cref{l:deflation,l:symplecticinflation} show that symplectic deflation and symplectic inflation are inverse operations.
\begin{lemma}
    \label{l:symplectic-deflation}
    Let $A\in\asmsk(2n)$, let $(b,a)$ be addable with $a<b$ and $\varnothing\ne U\subseteq N_A(b,a)$.
    Then symplectic deflation at $(b,a)$ for $U$ is possible if and only if every $(k,j)\in U\cap N_A(a,b)$ satisfies $k<j$.
\end{lemma}

We can visualize the statement by saying that $U$ is a subset of the grey region in the image below:
    
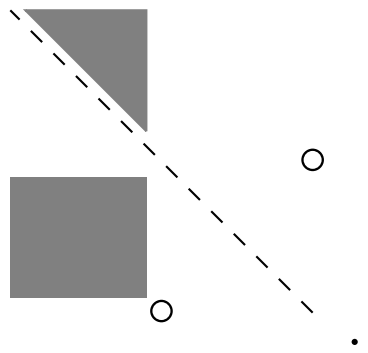
\[
    \begin{tikzpicture}[scale=.5]
        \draw[dashed] (2,-2) -- (-2,2);
        \fill [gray] (-2,-1.8) rectangle (-.2,-.2);
        \filldraw[gray] (-.2,.4) -- (-.2,2) -- (-1.8,2) -- (-.2,.4);
        \node at (2,0) {$\circ$};
        \node at (0,-2) {$\circ$};
        \end{tikzpicture}.
    \]
In light of \cref{l:symplectic-deflation}, for $(b,a)$ addable in $A \in \asmsk(2n)$, define
\[
N^\Sp_A(b,a) = \{(k,j) \in N_A(b,a) : k < j \ \mbox{or}\ k > a\}.
\]
That is, $N^\Sp_A(b,a)$ is obtained from $N_A(b,a)$ by discarding common neighbors below the diagonal.
For $U \subseteq N^\Sp_A(b,a)$ non-empty, we say that $(b,a),U$ is \newword{$\Sp$-addable}, since we can perform symplectic deflation at $(b,a)$ for $U$.

\begin{proof}[Proof of \cref{l:symplectic-deflation}]
    First, we prove the contrapositive of the forward direction.
    Unless the paired deflations are on disjoint
    entries and produce a skew ASM $A'$, we are done.
    Let $(k,j)\in U$ be a common neighbor of $(a,b)$ and $(b,a)$
    with $j<k\leq a<b$ minimizing $k-j$.
    Applying deflation at $(b,a)$ for $U$ results in a new $1$ in
    position $(k,\ell)$ with $\ell\leq a<b$.
    We claim that $(k,\ell)$ is a common neighbor for $(b,a)$
    and $(a,b)$ in $A'$.

    Clearly $(k,\ell)$ is northwest of $(a,b)$ and $(b,a)$.
    No unchanged nonzero entry can prevent it from being a common
    neighbor, since $(k,j)\prec(k,\ell)$ and $(k,j)$ was a common
    neighbor in $A$.
    By minimality of $k-j$, all selected neighbors preceding
    $(k,j)$ are upper neighbors.
    Consequently, apart from $(k,\ell)$ and its transpose, every
    newly inserted $1$ has a coordinate less than $\min(k,\ell)$.
    Its transpose is incomparable with it, since $A'$ has zero
    diagonal, so the claim follows.
    Then $(b,a)$ is not $\Sp$-removable, so
    $(b,a),U$ is not $\Sp$-addable.

    For the converse, assume that all
    $(j,k)\in U\cap N_A(a,b)$ satisfy $j<k$, so $U$ lies in the
    grey region.
    Therefore the paired deflations modify disjoint entries; let $A'$ be the
    resulting skew ASM.
    Any entry set equal to $1$ in the
    rectangular grey region lies below row $a$, so will not
    be a common neighbor of $(a,b)$ and $(b,a)$ in $A'$.
    Consider $(k,\ell)$ an entry set equal to $1$ in
    the triangular grey region.
    Necessarily, $U^T$ contains an entry below the diagonal, and
    the entry set equal to $1$ whose row index is that of the
    highest such neighbor has column greater than $a$ and lies
    southeast of $(k,\ell)$.
    Thus $(k,\ell)$ is not a common neighbor.
    By symmetry, $(b,a)$ and $(a,b)$ have no common neighbors in
    $A'$, while both remain removable so $(b,a)$ is $\Sp$-removable in $A'$.
    The result follows.
\end{proof}

We are now ready to define the \newword{symplectic key} $\keysp$ of an ASM.
We will do so by choosing a very particular sequence of symplectic inflations, but we will ultimately show that the symplectic key is independent of this choice (see \cref{c:symplectic-key}).

Let $A \in \asmsk(2n)$.
We define $\keysp(A)$ algorithmically:
\begin{enumerate}
    \item If $A$ has no $-1$'s, it is the permutation matrix for some $\mf$-involution $\tau$, and $\keysp(A) = \tau$.
    \item Otherwise, pick the lexicographically minimal $(b,a)$ with $b>a$ and $A_{ba}=-1$.
    \begin{enumerate}[(i)]
        \item If $(a,b)$ and $(b,a)$ share a common neighbor, then $\keysp(A) = \varnothing$.
        \item Otherwise, perform symplectic inflation at $(b,a)$ to obtain $A' \in \asmsk(2n)$.
        Set $\keysp(A) = \keysp(A')$ (computed recursively).
    \end{enumerate}
\end{enumerate}
Since symplectic inflation removes $-1$'s, the algorithm must terminate.
Let $\asmsp(2n)$ be the subset of $\asmsk(2n)$ with $\keysp(A) \neq \varnothing$, whose elements we call \newword{symplectic}, and for $\tau \in \mf_{2n}$ let $\asmsp(\tau)$ be the set of symplectic ASMs $A$ so that $\keysp(A )= \tau$.

\begin{remark}
    \label{rem:deflation-sequence}
    The same sequence of inflations used to construct $\keysp$ can be applied to any ASM, where we read in the same order but consider entries above the diagonal as well.
    This gives a canonical sequence of inflations from $A \in \asm(n)$ to the permutation matrix $w = \kappa(A)$. Inverting, this gives a \newword{deflation sequence} from $w$ to $A$.
    For $A \in \asmsp(2n)$, the algorithm to compute $\keysp(A)$ produces a sequence of symplectic inflations from $A$ to some permutation matrix $\tau$, or equivalently, a sequence of symplectic deflations from $\tau$ to $A$.
    We call the latter the \newword{symplectic deflation sequence} of $A$.
\end{remark}

\begin{lemma}
    Let $A\in\asmsk(2n)$.  If $\keysp(A) \neq \varnothing$ then $\keysp(A) = \key(A)$.
\end{lemma}
\begin{proof}
By \cref{l:symplecticinflation}, symplectic inflation at $(b,a)$ is the composition of the ordinary inflations at $(a,b)$ and $(b,a)$. Since the key of an ASM does not depend on the order in which removable $-1$'s are inflated~\cite{weigandt2021bumpless}*{Cor.~4.4}, the claim follows from the definitions of $\keysp$ and $\key$.
\end{proof}

We now demonstrate a skew ASM that is not symplectic, but whose $-1$'s are all $\Sp$-removable.
\begin{example}
\label{ex:sharedneighborsafteronestep}
    For the skew ASM $A$ below, both $(9,8)$ and $(10,7)$ are $\Sp$-removable, as they share no common neighbors with their partners.
    If we apply symplectic inflation to $(9,8)$, we obtain the skew ASM $A'$.
    \[
   A =  \left[\begin{array}{cccccccccccc}
        &&&&&\bullet\\
        &&&&&&\bullet\\
        &&&&&&&&&\bullet\\
        &&&&&&&\bullet\\
        &&&&&&&&\bullet\\
        \bullet\\
        &\bullet&&&&&&&&\circ&\bullet\\
        &&&\bullet&&&&&\circ&&&\bullet\\
        &&&&\bullet&&&\circ&&\bullet\\
        &&\bullet&&&&\circ&&\bullet\\
        &&&&&&\bullet\\
        &&&&&&&\bullet
    \end{array}\right]\,,
    \quad 
    A' = \left[\begin{array}{cccccccccccc}
        &&&&&\bullet\\
        &&&&&&\bullet\\
        &&&&&&&&&\bullet\\
        &&&&\bullet&&&&\\
        &&&\bullet&&&&&\\
        \bullet\\
        &\bullet&&&&&&&&\circ&\bullet\\
        &&&&&&&&&&&\bullet\\
        &&&&&&&&&\bullet\\
        &&\bullet&&&&\circ&&\bullet\\
        &&&&&&\bullet\\
        &&&&&&&\bullet
    \end{array}\right]
    \]
    We cannot apply symplectic inflation to $(10,7)$ in $A'$, as now $(10,7)$ and $(7,10)$ have common neighbors in positions $(4,5)$ and $(5,4)$.
    If we instead apply symplectic inflation to $(10,7)$ in $A$, we obtain
    \[
        A'' = \left[\begin{array}{cccccccccccc}
        &&&&&\bullet\\
        &&&&\bullet&&\\
        &&&&&&&&\bullet&\\
        &&&&&&&\bullet\\
        &\bullet&&&&&&&\\
        \bullet\\
        &&&&&&&&&&\bullet\\
        &&&\bullet&&&&&\circ&&&\bullet\\
        &&\bullet&&&&&\circ&&\bullet\\
        &&&&&&&&\bullet\\
        &&&&&&\bullet\\
        &&&&&&&\bullet
    \end{array}\right],
    \]
    after which we can no longer apply symplectic inflation to $(9,8)$.
\end{example}

\subsection{Transition on ASMs}
For $A \in \asmsp(2n)$, let $\negsetsp_A = \{(j,i):j >i, A_{ji} =-1\}$ and let $\Dsp_A = \{ (j,i) \in \rothe_A : j > i \}$, i.e., the set of cells in $\rothe_A$ that sit strictly below the main diagonal.
We define the symplectic weight of a symplectic ASM $A$ to be
\[\wt_{\Sp}(A)=\left(\prod_{(j,i)\in \Dsp_A}\beta (x_i\oplus x_j)\right)\left(\prod_{(j,i)\in \negsetsp_A}(1+\beta (x_i\oplus x_j))\right).\]

The main goal of the rest of this section is to prove the following.

\begin{theorem}
\label{thm:asmmain}
Given $\tau \in \mf_{2n}$,
    \begin{equation}
    \label{eq:asm-formula}
    \mathfrak G_\tau^{\Sp}=\beta^{-\ellfpf(\tau)}\sum_{A\in \asmsp(\tau)} \wt_{\Sp}(A).
    \end{equation}
\end{theorem}

We will do so by showing the right-hand side of \cref{eq:asm-formula} satisfies the transition recurrence for symplectic Grothendieck polynomials.
First, we recall Lascoux's transition map on ASMs (see \cite{lascoux2002chern}).
We refer the reader to \cite{weigandt2021bumpless}*{\S 5} for background.

Fix $w\in S_n$, and let $(a,b)$ be the lexicographically last diagram cell of $\rothe_w$. 
Let $U \subseteq N_w(a,b)$, where $U = \{(i_1,j_1),\dots,(i_h,j_h)\}$ with row indices $i_1 < \dots < i_h$.
Writing $t_{ij}=(i,j)$, we define 
\[
w_U = 
\begin{cases}
  wt_{a,w^{-1}(b)}, & U=\varnothing, \\
  wt_{a,w^{-1}(b)}t_{a,i_1}t_{a,i_2} \cdots t_{a,i_h}, & U\neq\varnothing.
\end{cases}
\]
\begin{lemma}
\label{lemma:normalASMtransitionmap}
Let $w\in S_n$, and let $(a,b)$ be the lexicographically last cell in $\rothe_w$.
Form an $n\times n$ matrix $M_{w,(a,b)}$ that has $1$'s in positions $(a,b)$ and $(w^{-1}(b),w(a))$, and $-1$'s in positions $(w^{-1}(b),b)$ and $(a,w(a))$.  
The map $A\mapsto A+M_{w,(a,b)}$ defines a bijection from $\asm(w)$ to 
\[
\bigsqcup_{U\subseteq N_w(a,b)} \asm(w_U).
\]
Furthermore, the preimage of $\asm(w_U)$ is the set of ASMs in $\asm(w)$ for which the first step in the deflation sequence is at $(a,b)$ for $U$ (when $U = \varnothing$, no deflation occurs at $(a,b)$).
\end{lemma}

We note that the map $A\mapsto A+M_{w,(a,b)}$ can also be visualized by making the replacements pictured below, according to the two cases where $A_{a,b} = 0$ or $A_{a,b} = -1$.
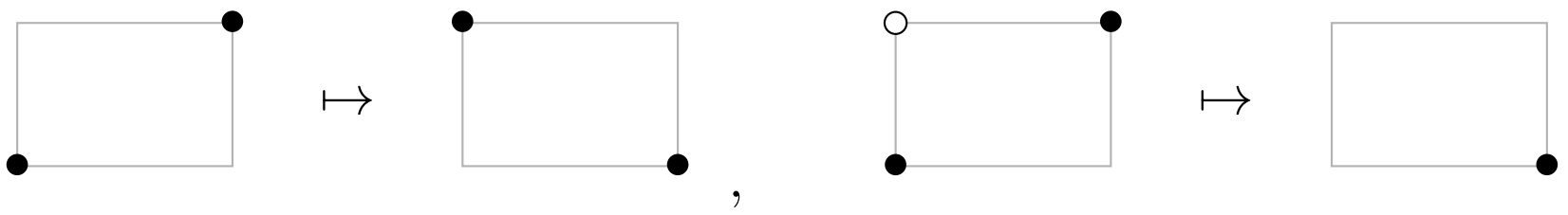
\begin{equation}
    \label{eq:asm-map}
\begin{tikzpicture}[scale=.5]
  \draw[gray!60, thin] (0,0) rectangle (3,-2);
  \node[] (A) at (0,0) {\,};
  \node (B) at (0,-2) {$\bullet$};
  \node (C) at (3,0) {$\bullet$};
  \node (D) at (3,-2) {\,};
\end{tikzpicture}
\hspace{1em}
\raisebox{1.5em}{$\mapsto$}
\hspace{1em}
\begin{tikzpicture}[scale=.5]
  \draw[gray!60, thin] (0,0) rectangle (3,-2);
  \node[] (A) at (0,0) {$\bullet$};
  \node (B) at (0,-2) {\,};
  \node (C) at (3,0) {\,};
  \node (D) at (3,-2) {$\bullet$};
\end{tikzpicture}\ ,
\hspace{2em}
\begin{tikzpicture}[scale=.5]
  \draw[gray!60, thin] (0,0) rectangle (3,-2);
  \node[] (A) at (0,0) {$\whitecirc$};
  \node (B) at (0,-2) {$\bullet$};
  \node (C) at (3,0) {$\bullet$};
  \node (D) at (3,-2) {\,};
\end{tikzpicture}
\hspace{1em}
\raisebox{1.5em}{$\mapsto$}
\hspace{1em}
\begin{tikzpicture}[scale=.5]
  \draw[gray!60, thin] (0,0) rectangle (3,-2);
  \node[] (A) at (0,0) {\,};
  \node (B) at (0,-2) {\,};
  \node (C) at (3,0) {\,};
  \node (D) at (3,-2) {$\bullet$};
\end{tikzpicture}
\end{equation}
In both cases, the replacement occurs on the rectangle $[a,w^{-1}(b)]\times[b,w(a)]$.

\begin{lemma}
\label{l:I-neighbors}
Let $\tau\in\mf_{2n}$ with lexicographically maximal cell $(b,a)$ in $\Dsp_\tau$.
Then for $\rho=u_{b,\tau(a)}(\tau)$,
\[
I_\Sp(\rho,a) = \{i: (j,i) \in N^\Sp_\tau(b,a)\}.
\]
\end{lemma}

Let $\tau \in \mf_{2n}$ with lexicographically maximal cell $(b,a)$ in $\Dsp_\tau$.
For $U\subseteq N^\Sp_\tau(b,a)$, write
\[
U=\{(r_1,c_1),\ldots,(r_h,c_h)\}
\quad \mbox{with} \quad  c_1<\cdots<c_h,
\]
so $r_1>\cdots>r_h$ and $r_j=\tau(c_j)$.
Define
\[
\tau_U=u_{c_h,a}\cdots u_{c_1,a}\bigl(u_{b,\tau(a)}(\tau)\bigr).
\]
In particular, $\tau_\varnothing=\rho$ from \cref{c:g_tree}, and \cref{l:I-neighbors} identifies the column indices of $U$ with elements of $I_\Sp(\rho,a)$.
Translating \cref{lemma:normalASMtransitionmap} to symplectic ASMs, we have:
\begin{proposition}
\label{prop:sympasmbijection}
Let $\tau\in\mf_{2n}$ with $(b,a)$ the lexicographically maximal element of $\Dsp_\tau$.
The map $A\mapsto A+M_{\tau,(b,a)}+M_{\tau,(a,b)}$ defines a bijection
\[
\asmsp(\tau)\longrightarrow
\bigsqcup_{U\subseteq N^\Sp_\tau(b,a)}\asmsp(\tau_U).
\]
For $U\neq\varnothing$, the preimage of $\asmsp(\tau_U)$ consists of those ASMs whose first symplectic deflation is at $(b,a)$ for $U$.
For $U=\varnothing$, it consists of those for which no deflation occurs at $(b,a)$.
\end{proposition}

\begin{proof}
    The key to our proof is the following observation: fix $A\in \asmsp(\tau)$, and let $A' = A + M_{\tau,(b,a)} + M_{\tau,(a,b)} \in \asmsk(2n)$. 
    If $(j,i)$ is strictly west or north of $(b,a)$, then $(j,i),U$ is $\Sp$-addable in $A$ if and only if it is $\Sp$-addable in $A'$.
    With this in mind, we first establish that the given map defines a bijection from those $A\in \asmsp(\tau)$ with $(b,a)\in \Dsp_A$ to $\asmsp(\tau_\varnothing)$.
    Here, we apply the first move from~\cref{eq:asm-map} symmetrically.
    This move maps $\tau$ (as a matrix) to $\tau_\varnothing$, which corresponds to an $\mf$-Bruhat cover.
    For other such $A$'s, any deflation in its symplectic deflation sequence will be north or west of $(b,a)$, so the same sequence can be applied equally to $\tau$ and $\tau_\varnothing$.
    The result follows.

    For $U = \{(i_1,j_1),\dots,(i_h,j_h)\} \subseteq N^\Sp_\tau(b,a)$ non-empty with $j_1 < \dots < j_h$, let $A_U$ be the result of applying symplectic deflation to $\tau$ at $(b,a)$ for $U$.
    Since symmetric application of the second move from~\cref{eq:asm-map} removes the only two $-1$'s in $A_U$, the resulting ASM is in fact an $\mf$-involution.
    We claim this is $\tau_U$.
    We argue by induction, with our base case being $U=\varnothing$.
    The induction hypothesis shows the result holds for $U' = \{(i_1,j_1),\dots,(i_{h-1},j_{h-1})\}$.
    Applying $u_{j_h,a}$ to $\tau_{U'}$ is equivalent to swapping rows $j_h$ and $a$ and columns $j_h$ and $a$.
    Direct inspection shows this is the difference between deflation at $(b,a)$ for $U$ and for $U'$.
    By our initial observation, the rest of the argument proceeds as in the case where $U = \varnothing$.
    By \cref{l:children}, every deflation in a symplectic
deflation sequence from $\tau_U$ precedes $(b,a)$,
so the same argument applies in reverse.
\end{proof}

\begin{remark}
\label{rem:asm-tree}
\cref{prop:sympasmbijection} allows us to interpret the symplectic deflation sequence of $A \in \asmsp(\tau)$ as a maximal path in the symplectic LS tree.
For $\tau_U$ as defined in \cref{prop:sympasmbijection},
in each application of \cref{c:g_tree}, symplectic deflation at $(b,a)$ for $U$ assigns $A$ to a child $\tau_U$.
Iterating, we assign $A$ to one of $\tau_U$'s children, and so on.
When no further symplectic deflations will be applied, all further choices will be $U = \varnothing$, continuing until the resulting child is $\mf$-dominant.
\end{remark}

We are now ready to prove \cref{thm:asmmain}.
\begin{proof}[Proof of \cref{thm:asmmain}]
Both sides are $1$ for $\tau=\Theta_{2n}$.
Otherwise, let $(b,a)$ be the last cell of $\Dsp_\tau$ in reading order.
By \cref{l:children}, induction applies to every $\tau_U$.
    Furthermore, note that $\ellfpf(\tau_{U})=\ellfpf(\tau)+|U|-1$ since each $u$-operator corresponds to an $\mf$-Bruhat cover by the $\prec$-maximality of neighbors.
    For $A \in \asmsp(\tau)$, if no symplectic deflation occurs at $(b,a)$ (corresponding to $U = \varnothing$), the map in \cref{prop:sympasmbijection} removes a factor of $\beta (x_a\oplus x_b)$ from $\wt_\Sp(A)$, since $(b,a)$ is no longer a diagram cell.
    Otherwise, it removes a factor of $1 + \beta (x_a \oplus x_b)$ as $(b,a)$ is no longer a $-1$.
    Therefore,
\begin{align*}
    \beta^{\ellfpf(\tau)}\mathfrak G_{\tau}^{\Sp}&=\beta(x_a \oplus x_b )\beta^{\ellfpf(\tau_\varnothing)}\fkGsp_{\tau_\varnothing} + (1+\beta (x_b \oplus x_a)) \left[\sum_{\varnothing\neq U\subseteq N^{\Sp}_\tau(b,a)}\beta^{\ellfpf(\tau_U)}\fkGsp_{\tau_U} \right]\\
    &=\beta(x_a \oplus x_b )\sum_{A\in \asmsp(\tau_\varnothing)} \wt_{\Sp}(A)+(1+\beta (x_b \oplus x_a)) \left[\sum_{\varnothing\neq U\subseteq N^{\Sp}_\tau(b,a)}\sum_{A\in \asmsp(\tau_U)} \wt_{\Sp}(A) \right]\\
    &=\sum_{A\in \asmsp(\tau)} \wt_{\Sp}(A),
\end{align*}
where the first equality is the transition formula, the second follows by induction, and the third is an application of \cref{prop:sympasmbijection}.\end{proof}

\section{Symplectic bumpless pipe dreams}
\label{section:symplecticbpd}

In this section we introduce symplectic bumpless pipe dreams, which we will show are the restrictions of Rothe diagrams of symplectic ASMs to $\lrtri_{2n}$.
In order to establish this fact, we require several operations on these objects.
We then use this correspondence to show the symplectic key $\keysp(A)$ is independent of the choice of symplectic inflations used to identify it.
This allows us to show that $\bpdsp(\tau)$ is connected by droop moves, symplectic droop moves, and their $K$-theoretic analogues.  This is analogous to the corresponding result for ordinary BPDs in~\cite{weigandt2021bumpless}.

\subsection{Symplectic BPDs}
\label{subsection:sympBPD}
A \newword{triangular bumpless pipe dream} is a filling of the diagram $\lrtri_{2n}$ together with the triangular half-cells along its diagonal, using the tiles
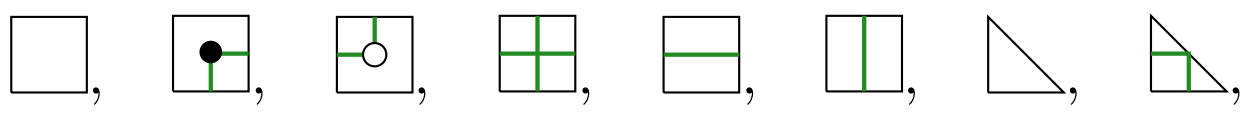
\[
\begin{tikzpicture}[scale=.5]
	\draw (0,0) -- (1,0) -- (1,1) -- (0,1) -- (0,0);
\end{tikzpicture},
\quad
\begin{tikzpicture}[scale=.5]
	\draw (0,0) -- (1,0) -- (1,1) -- (0,1) -- (0,0);
	\draw[thick,ForestGreen] (.5,0) -- (.5,.5) -- (1,.5);
	\node at (.5,.5) {$\bullet$};
\end{tikzpicture},
\quad
\begin{tikzpicture}[scale=.5]
	\draw (0,0) -- (1,0) -- (1,1) -- (0,1) -- (0,0);
	\draw[thick,ForestGreen] (0,.5) -- (.5,.5) -- (.5,1);
	\node at (.5,.5) {\whitecirc};
\end{tikzpicture},
\quad
\begin{tikzpicture}[scale=.5]
	\draw (0,0) -- (1,0) -- (1,1) -- (0,1) -- (0,0);
	\draw[thick,ForestGreen] (.5,0) -- (.5,1);
	\draw[thick,ForestGreen] (0,.5) -- (1,.5);
\end{tikzpicture},
\quad
\begin{tikzpicture}[scale=.5]
	\draw (0,0) -- (1,0) -- (1,1) -- (0,1) -- (0,0);
	\draw[thick,ForestGreen] (0,.5) -- (1,.5);
\end{tikzpicture},
\quad
\begin{tikzpicture}[scale=.5]
	\draw (0,0) -- (1,0) -- (1,1) -- (0,1) -- (0,0);
	\draw[thick,ForestGreen] (.5,0) -- (.5,1);
\end{tikzpicture},
\quad
\begin{tikzpicture}[scale=.5]
	\draw (0,0) -- (1,0) -- (0,1) -- (0,0);
\end{tikzpicture},
\quad
\begin{tikzpicture}[scale=.5]
	\draw (0,0) -- (1,0) -- (0,1) -- (0,0);
	\draw[thick,ForestGreen] (.5,0) -- (.5,.5) -- (0,.5);
\end{tikzpicture},
\]
with the first six tiles in square cells and the last two in diagonal half-cells, satisfying:
\begin{enumerate}
    \item pipes agree along shared edges of adjacent cells;
    \item each cell in the bottom row has a pipe entering it from below;
    \item no cell in the leftmost column has a pipe entering it from the left.
\end{enumerate}
The tiles $\cross$ and $\yestri$ are called \newword{cross} and \newword{reflect}, respectively. For a triangular bumpless pipe dream $B$ and a tile $t$, let $B(t)$ be the set of positions $(j,i)$ with $B(j,i) = t$.

\begin{example}
\label{ex:triangular-bpd-pictures}
Below are three triangular bumpless pipe dreams:
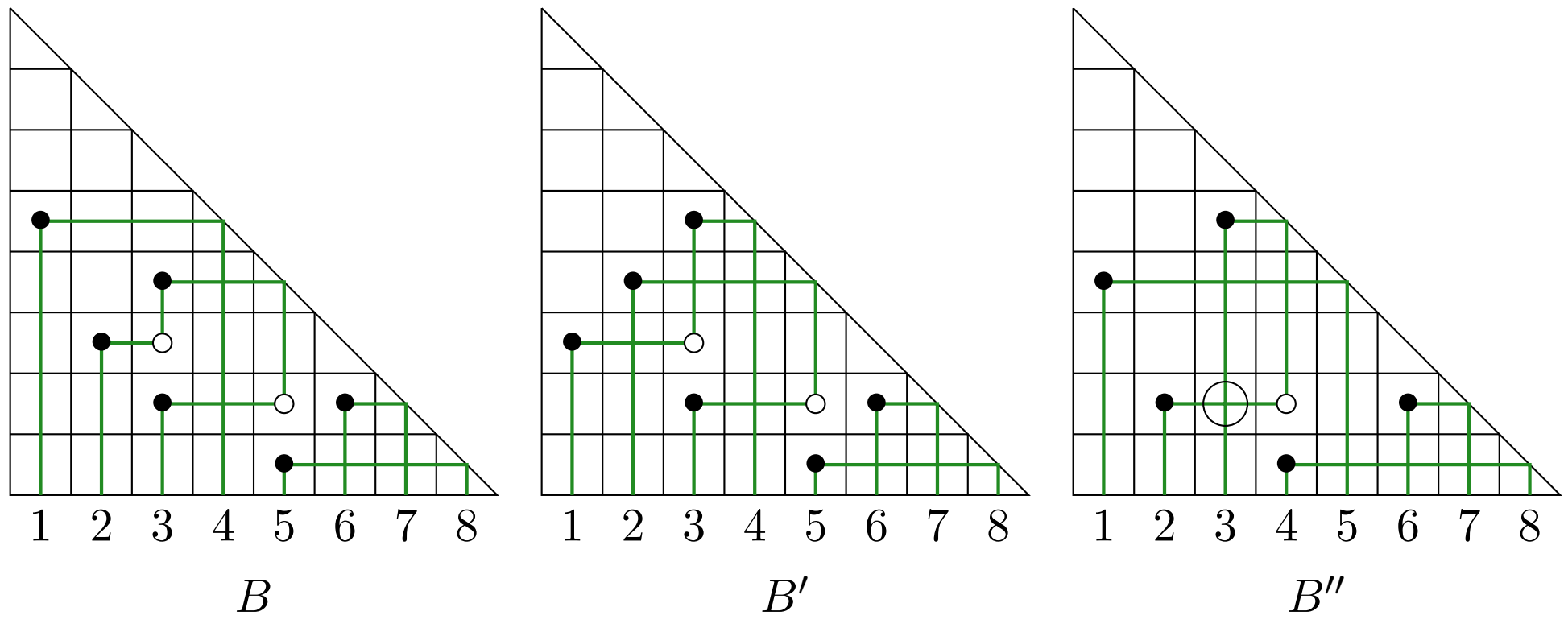
\[	\begin{tabular}{ccc}
\begin{tikzpicture}[scale=.5]
\draw (.5,7.5) -- (.5,-.5) -- (8.5,-.5) -- (.5,7.5);
\draw (.5,6.5) -- (1.5,6.5) -- (1.5,-.5);
\draw (.5,5.5) -- (2.5,5.5) -- (2.5,-.5);
\draw (.5,4.5) -- (3.5,4.5) -- (3.5,-.5);
\draw (.5,3.5) -- (4.5,3.5) -- (4.5,-.5);
\draw (.5,2.5) -- (5.5,2.5) -- (5.5,-.5);
\draw (.5,1.5) -- (6.5,1.5) -- (6.5,-.5);
\draw (.5,0.5) -- (7.5,0.5) -- (7.5,-.5);
\draw[thick,ForestGreen] (1,-.5) -- (1,4)  -- (4,4)--(4,-.5);
\draw[thick,ForestGreen] (2,-.5) --(2,2)--(3,2)--(3,3)--(5,3)--(5,1)--(3,1)--(3,-.5);
\draw[thick,ForestGreen] (5,-.5) --(5,0)--(8,0)--(8,-.5);
\draw[thick,ForestGreen] (6,-.5) --(6,1)--(7,1)--(7,-.5);
\node at (1,4) {$\bullet$};
\node at (2,2) {$\bullet$};
\node at (3,1) {$\bullet$};
\node at (3,2) {$\whitecirc$};
\node at (5,1) {$\whitecirc$};
\node at (3,3) {$\bullet$};
\node at (5,0) {$\bullet$};
\node at (6,1) {$\bullet$};
\node at (1,-1) {1};
\node at (2,-1) {2};
\node at (3,-1) {3};
\node at (4,-1) {4};
\node at (5,-1) {5};
\node at (6,-1) {6};
\node at (7,-1) {7};
\node at (8,-1) {8};
\end{tikzpicture} 
&
\begin{tikzpicture}[scale=.5]
\draw (.5,7.5) -- (.5,-.5) -- (8.5,-.5) -- (.5,7.5);
\draw (.5,6.5) -- (1.5,6.5) -- (1.5,-.5);
\draw (.5,5.5) -- (2.5,5.5) -- (2.5,-.5);
\draw (.5,4.5) -- (3.5,4.5) -- (3.5,-.5);
\draw (.5,3.5) -- (4.5,3.5) -- (4.5,-.5);
\draw (.5,2.5) -- (5.5,2.5) -- (5.5,-.5);
\draw (.5,1.5) -- (6.5,1.5) -- (6.5,-.5);
\draw (.5,0.5) -- (7.5,0.5) -- (7.5,-.5);

\draw[thick,ForestGreen] (1,-.5) -- (1,2)  -- (3,2)--(3,4)--(4,4)--(4,-.5);
\draw[thick,ForestGreen] (2,-.5) --(2,3)--(5,3)--(5,1)--(3,1)--(3,-.5);
\draw[thick,ForestGreen] (5,-.5) --(5,0)--(8,0)--(8,-.5);
\draw[thick,ForestGreen] (6,-.5) --(6,1)--(7,1)--(7,-.5);
\node at (1,2) {$\bullet$};
\node at (2,3) {$\bullet$};
\node at (3,1) {$\bullet$};
\node at (3,2) {$\whitecirc$};
\node at (5,1) {$\whitecirc$};
\node at (3,4) {$\bullet$};
\node at (5,0) {$\bullet$};
\node at (6,1) {$\bullet$};
\node at (1,-1) {1};
\node at (2,-1) {2};
\node at (3,-1) {3};
\node at (4,-1) {4};
\node at (5,-1) {5};
\node at (6,-1) {6};
\node at (7,-1) {7};
\node at (8,-1) {8};
\end{tikzpicture} 
& \begin{tikzpicture}[scale=.5]
\draw (.5,7.5) -- (.5,-.5) -- (8.5,-.5) -- (.5,7.5);
\draw (.5,6.5) -- (1.5,6.5) -- (1.5,-.5);
\draw (.5,5.5) -- (2.5,5.5) -- (2.5,-.5);
\draw (.5,4.5) -- (3.5,4.5) -- (3.5,-.5);
\draw (.5,3.5) -- (4.5,3.5) -- (4.5,-.5);
\draw (.5,2.5) -- (5.5,2.5) -- (5.5,-.5);
\draw (.5,1.5) -- (6.5,1.5) -- (6.5,-.5);
\draw (.5,0.5) -- (7.5,0.5) -- (7.5,-.5);

\draw[thick,ForestGreen] (1,-.5) -- (1,3)  -- (2,3)--(5,3)--(5,-.5);
\draw[thick,ForestGreen] (2,-.5) --(2,1)--(4,1)--(4,4)--(3,4)--(3,-.5);
\draw[thick,ForestGreen] (4,-.5) --(4,0)--(8,0)--(8,-.5);
\draw[thick,ForestGreen] (6,-.5) --(6,1)--(7,1)--(7,-.5);
\node at (1,3) {$\bullet$};
\node at (2,1) {$\bullet$};
\node at (4,1) {$\whitecirc$};
\node at (3,4) {$\bullet$};
\node at (4,0) {$\bullet$};
\node at (6,1) {$\bullet$};
\node at (1,-1) {1};
\node at (2,-1) {2};
\node at (3,-1) {3};
\node at (4,-1) {4};
\node at (5,-1) {5};
\node at (6,-1) {6};
\node at (7,-1) {7};
\node at (8,-1) {8};

\node[draw,circle] at (3,1) {};
\end{tikzpicture}
\\
$B$ & $B'$ & $B''$
\end{tabular}
\]
\end{example}

For $A\in\asmsk(2n)$, let $\Psi(A)$ be the triangular filling obtained by restricting the BPD associated to $A$ to the lower triangle, including the diagonal half-cells.
Its set of blank cells is $\Dsp_A=\rothe_A\cap\lrtri_{2n}$.
As with $\Dsp_\tau$, we also use $\Dsp_A$ for the corresponding tile filling.

\begin{proposition}
    \label{p:triangular-bpd}
    The map $\Psi$  is a bijection from $\asmsk(2n)$ to triangular BPDs on $\lrtri_{2n}$.
\end{proposition}
\begin{proof}
The usual correspondence from (ordinary) ASMs to BPDs is a bijection.
This map immediately restricts to a bijection between symmetric ASMs and symmetric BPDs.
Further restricting to skew ASMs imposes the condition that there are only $0$'s on the diagonal, which means the BPD of $A$ has no $\nwelbow$'s or $\seelbow$'s on the diagonal. 
Thus, the diagonal of $\rothe_A$ consists of $\cross$ and $\nowire$ tiles, and so restricting the filling to the lower triangle, including its diagonal half-cells, gives a triangular BPD.

On the other hand, starting with a triangular BPD $B$ on $\lrtri_{2n}$, it is clear that taking the union of $B$ with its reflection across the diagonal produces a full BPD that is symmetric.  Furthermore, this BPD has no $\nwelbow$ or $\seelbow$ on its diagonal, thus it is the image of a skew ASM.
\end{proof}

To each triangular bumpless pipe dream $B$ with $k$ crossings, associate the word $\mathbf{a}_B=(a_1,\ldots,a_k)$ as follows.
First, order the crossings of $B$ from right to left within each row, from top to bottom.
Then $a_i$ is the number of pipes strictly to the left of the $i$th crossing, i.e., pick a point immediately below and to the left of the crossing, then record how many pipes are to the left of this point.
If $\mathbf{a}_B$ is a symplectic Hecke word for $\tau \in \mf_{2n}$, then we say that $B$ is a \newword{symplectic BPD} and write $\demprodfpf(B)=\tau$.
A symplectic bumpless pipe dream $B$ is \newword{reduced} if $\mathbf{a}_B$ is a reduced symplectic Hecke word.
Given $\tau \in \mf_{2n}$, let $\bpdsp(\tau)$ be the set of symplectic bumpless pipe dreams of size $2n$ whose words are symplectic Hecke words for $\tau$, and let $\bpdsp_\mr(\tau)$ be the subset of reduced elements of $\bpdsp(\tau)$.
Note that $\Dsp_\tau \in \bpdsp_\mr(\tau)$ for all $\tau \in \mf_{2n}$.

It is useful to be able to visually read off the matching associated to a symplectic BPD by reading along its pipes.
To do so, fix $B\in\bpdsp(\tau)$ with $k$ crossings, ordered as before, and associated word $\mathbf{a}_B=(a_1,\ldots,a_k)$.
We obtain a sequence of matchings as follows.
Let $\tau^{(0)}=\Theta_{2n}$.
Then we define $\tau^{(i)}=\tau^{(i-1)} \iprod s_{a_i}$ for all $i\in [k]$.
Note by definition that $\tau^{(k)} = \demprodfpf(B) = \tau$.
For all $i\in[k]$, if $\tau^{(i)}>\tau^{(i-1)}$, then we call crossing $i$ an \newword{active cross}.
Otherwise, we must have $\tau^{(i)}=\tau^{(i-1)}$ and we say that crossing $i$ is a \newword{passive cross}.

\begin{lemma}
    \label{l:reading-word}
    For $\tau \in \mf_{2n}$ and $B \in \bpdsp(\tau)$, replace all passive crosses with bumping tiles $\bump$ and label the pipes along the bottom edge of the grid with the numbers $1,2,\ldots, 2n$.
When two pipes cross, propagate their labels as pictured below:
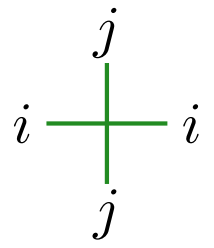
\[\begin{tikzpicture}[scale=.4]
  \draw[thick,color=ForestGreen] (-1,0) -- (1,0);
  \node at (-1.4,0) {$i$};
  \node at (1.4,0) {$i$};
  \draw[thick,color=ForestGreen] (0,-1) -- (0,1);
  \node at (0,1.5) {$j$};
  \node at (0,-1.5) {$j$};
\end{tikzpicture}.\]
Then for all $(j,i) \in \Cyc(\tau)$, the pipes $i$ and $j$ meet in a reflect tile $\yestri$.
\end{lemma}

\begin{proof}
For $(k,j) \in \lrtri_{2n}$, consider the restriction $B'$ of $B$ to the rows above $k$ as well as the $k$th row weakly to the right of the $j$th column.
Label the pipes at the southeast boundary of this filling (including each $\yestri$ and $\notri$ tile) with $1,2,\ldots,2n$ from top left to bottom right.
Treating each passive cross as a $\bump$, the pairs of pipes meeting at the diagonal determine a matching.
If $B(k,j)$ is not a cross, it is clear that removing this tile does not change the matching.
Therefore, we reduce to the case where $(k,j)$ is the $i$th cross of $B$.
We claim the matching for $B'$ is $\tau^{(i)}$.
In particular, the matching for all of $B$ will be $\tau$ as claimed.

Our proof is by induction, the base case being when the diagram is unfilled so the matching is $\tau^{(0)} = \Theta_{2n}$ with each cycle corresponding to one of the $n$ $\yestri$ tiles on the diagonal.
Inductively, we may assume the matching of  $B' - (k,j)$ is $\tau^{(i-1)}$.
If the cross in $(k,j)$ is passive, direct inspection shows the matching for $B'$ is $\tau^{(i-1)} = \tau^{(i)}$.
Otherwise, it swaps the pipes labeled $a_i$ and $a_i{+}1$, so the matching is $s_{a_i} \tau^{(i-1)} s_{a_i} = \tau^{(i-1)} \iprod s_{a_i} = \tau^{(i)}$ as desired.
\end{proof}

Note that when $B$ is reduced, there are no passive crossings, so we may directly read the matching associated to $B$ by tracing along pipes.
In particular, reduced symplectic bumpless pipe dreams for $\tau$ have $\ellfpf(\tau)$ cross tiles.

\begin{proposition}
    \label{p:reduced-doubling}
    For $\tau \in \mf_{2n}$, the map $B \mapsto B \cup B^T$ is a bijection from $\bpdsp_\mr(\tau)$ to the reduced symmetric BPDs for $\tau$.
\end{proposition}
\begin{proof}
For any triangular BPD  $B$ with reading word $\bfa = (a_1,\ldots,a_k)$,
 we observe that $B \cup B^T$ has $2k+n$ crosses, including the $n$ on the diagonal.
Tracing its pipes as in \cref{l:reading-word} gives the permutation
\[
w=s_{a_k}\cdots s_{a_1}\Theta_{2n}s_{a_1}\cdots s_{a_k} = \demprodfpf(\bfa),
\]
with the final equality since $\bfa$ is reduced as $B \in \bpdsp(\tau)$.
Therefore, $B \cup B^T \in \bpd(\tau)$ and is reduced since each conjugation is non-trivial.
Conversely, a symmetric reduced BPD for $\tau\in\mf_{2n}$ has no diagonal elbows: a pipe through such an elbow is fixed by reflection and would give a fixed point of $\tau$.
It therefore restricts to a triangular BPD, and restriction is inverse to our map.
\end{proof}

We now characterize $\Dsp_\tau$, viewed as a triangular bumpless pipe dream.

\begin{lemma}
    \label{l:diagram-bpd}
    For $\tau \in \mf_{2n}$, $\Dsp_\tau$ is the unique element of $\bpdsp(\tau)$ with no $\nwelbow$'s.
\end{lemma}

\begin{proof}
If a triangular bumpless pipe dream $B$ has no $\nwelbow$'s, then $\Psi^{-1}(B)$ is the permutation matrix of some $\tau \in \mf_{2n}$.
Therefore $B = \Dsp_\tau$ for some $\tau\in \mf_{2n}$, so the result follows.
\end{proof}

The following is an immediate consequence of~\cite{weigandt2021bumpless}*{Lem.~3.7~(1)}.

\begin{lemma}
    \label{l:reduced-blanks}
    Let $\tau \in \mf_{2n}$ and $B \in \bpdsp(\tau)$.
    Then $B$ has as many $\nowire$ tiles as $\cross$ tiles.
    In particular $B$ is reduced if and only if $|B(\nowire)| = \ellfpf(\tau)$, and otherwise $|B(\nowire)| > \ellfpf(\tau)$.
\end{lemma}

\begin{example}
In \cref{ex:triangular-bpd-pictures}, $B$ is reduced, $B'$ is symplectic but not reduced, and $B''$ is not a symplectic bumpless pipe dream.
By tracing pipes, we see $\demprodfpf(B) = (1,4)(2,3)(5,8)(6,7)$.
To see why $B''$ is not symplectic, note its crosses, in the order used to define $\mathbf{a}_B$, are in positions
\[
(5,4),(5,3),(7,3),(8,7),(8,6),(8,5).
\]
Recall that to each cross we associate the number of pipes to the left of the point just below and left of the crossing.
Therefore $\bfa_{B''} = (2,1,2,6,5,4)$.
However, $\Theta_{2n} \iprod s_{2}\iprod s_1\iprod s_2=\varnothing$, so $\mathbf a_{B''}$ is not a symplectic Hecke word for any $\tau\in \mf_{2n}$.
This is witnessed in $B''$ by the fact that the third (circled) cross is between two pipes that meet on the diagonal.

The crossings for $B'$, in the same order, are in positions
\[
(5,4),(5,3),(6,2),(7,4),(8,7),(8,6),
\]
so $\bfa_{B'} = (2,1,1,3,6,5)$.
Here the third and fourth crossings are passive.
Replacing the passive crossings with $\bump$ tiles as in \cref{l:reading-word}, we see $\demprodfpf(B') = (1,4)(2,3)(5,8)(6,7)$:
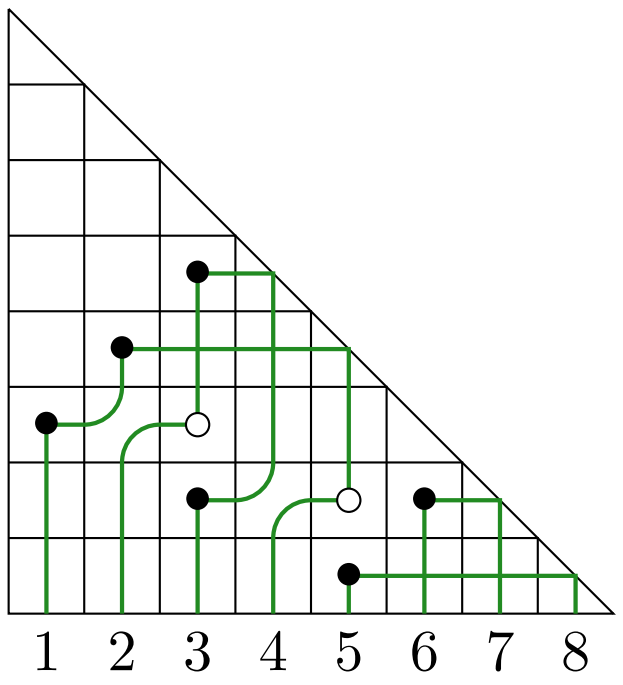
\[
\begin{tikzpicture}[scale=.5]
\draw (.5,7.5) -- (.5,-.5) -- (8.5,-.5) -- (.5,7.5);
\draw (.5,6.5) -- (1.5,6.5) -- (1.5,-.5);
\draw (.5,5.5) -- (2.5,5.5) -- (2.5,-.5);
\draw (.5,4.5) -- (3.5,4.5) -- (3.5,-.5);
\draw (.5,3.5) -- (4.5,3.5) -- (4.5,-.5);
\draw (.5,2.5) -- (5.5,2.5) -- (5.5,-.5);
\draw (.5,1.5) -- (6.5,1.5) -- (6.5,-.5);
\draw (.5,0.5) -- (7.5,0.5) -- (7.5,-.5);

\draw[thick,ForestGreen] (1,-.5) -- (1,2) -- (1.5,2) to [bend right=45] (2,2.5) -- (2,3)--(5,3)--(5,1) -- (4.5,1) to [bend right=45] (4,.5) -- (4,-.5);
\draw[thick,ForestGreen] (2,-.5) -- (2,1.5) to [bend left=45] (2.5,2) -- (3,2)-- (3,4)--(4,4)--(4,1.5) to [bend left=45] (3.5,1) --(3,1)--(3,-.5);
\draw[thick,ForestGreen] (5,-.5) --(5,0)--(8,0)--(8,-.5);
\draw[thick,ForestGreen] (6,-.5) --(6,1)--(7,1)--(7,-.5);
\node at (1,2) {$\bullet$};
\node at (2,3) {$\bullet$};
\node at (3,1) {$\bullet$};
\node at (3,2) {$\whitecirc$};
\node at (5,1) {$\whitecirc$};
\node at (3,4) {$\bullet$};
\node at (5,0) {$\bullet$};
\node at (6,1) {$\bullet$};
\node at (1,-1) {1};
\node at (2,-1) {2};
\node at (3,-1) {3};
\node at (4,-1) {4};
\node at (5,-1) {5};
\node at (6,-1) {6};
\node at (7,-1) {7};
\node at (8,-1) {8};
\end{tikzpicture} \qedhere
\]
\end{example}

\subsection{Triangular planar histories}

  Fix a triangular region in the plane with vertices at positions $(0,0), (m,0), (0,m)$ for some positive integer $m$.  A \newword{triangular planar history} is a filling of this triangle with pseudolines that satisfies the following conditions:
\begin{enumerate}
    \item each pseudoline starts on the bottom of the triangle,
    \item each segment of a pseudoline travels north or east,
    \item distinct pseudolines do not travel concurrently (i.e., they only overlap at points),
    \item each pseudoline ends at the diagonal edge of the triangle, and
    \item each pseudoline is paired with a distinct pseudoline that it meets on the diagonal.  
\end{enumerate}  
We assume finitely many pseudolines and intersections, distinct bottom endpoints, and no intersections involving three or more pseudolines.

We interpret triangular BPDs as special cases of  triangular planar histories.  
Also, in keeping with our visual representation of Rothe diagrams of ASMs, we draw turns in pseudolines with $\seelbow$ and $\nwelbow$.  
It is possible for two pseudolines to meet but not cross; we represent this situation by drawing $\bump$.  Furthermore, in keeping with the literature on pipe dream combinatorics, we refer to pseudolines as \newword{pipes}.

  By interpreting each $\yestri$ on the diagonal as a place where two pipes almost meet, but do not cross, we may readily apply results about \emph{rectangular planar histories} to triangular planar histories.   See \cite{Weigandt.Bases}*{\S 3} for background on rectangular planar histories.
  
  Given a triangular planar history $T$, we assign a reading word $\mathbf a_T$ in the same way as we did for triangular BPDs.  
  We define $\demprodfpf(T)=\demprodfpf(\mathbf a_T)$.  Reversing the word $\mathbf a_T$ produces the word one would normally assign to a rectangular planar history.  
  If $\demprodfpf(T)\neq \varnothing$, we say $T$ is a \newword{symplectic} planar history.

\begin{example}
\label{example:triangularnotsymplectic}
Consider the triangular planar history $T$ pictured below.
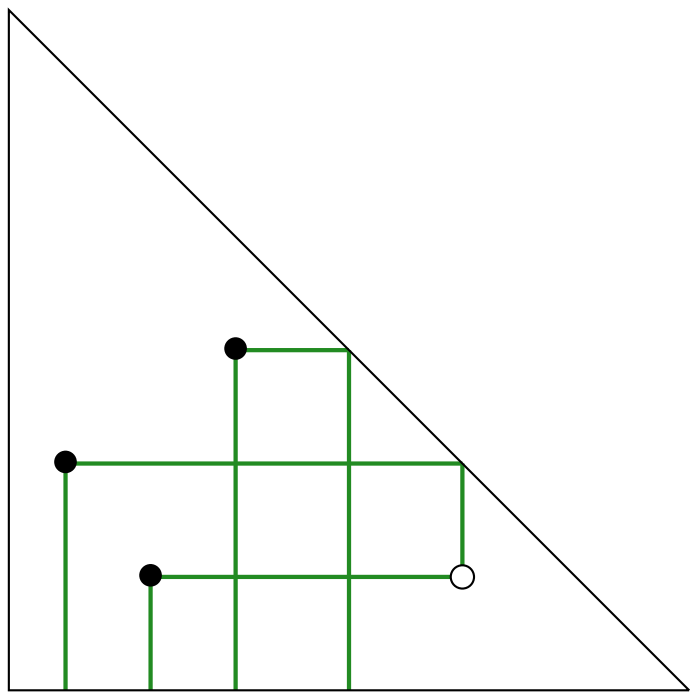
\[
    \begin{tikzpicture}[scale=.75]
		\draw[ForestGreen,thick] (.5,0)--(.5,2)--(4,2);
        \draw[ForestGreen,thick] (4,2)--(4,1)--(1.25,1)--(1.25,0);
		\draw[ForestGreen,thick] (3,0)--(3,3);
        \draw[ForestGreen,thick] (3,3)--(2,3)--(2,0);
		\draw (6,0) -- (0,6)--(0,0)--(6,0);
        \node at (.5,2) {$\bullet$};
        \node at (4,1) {$\whitecirc$};
        \node at (1.25,1) {$\bullet$};
        \node at (2,3) {$\bullet$};
	\end{tikzpicture}\]
    We have $\mathbf {a}_T=(2,1,3,2)$.
   Because $\mathbf {a}_T=(2,1,3,2)\sim_\Sp (2,3,3,2)\sim_\Sp (2,3,2)\sim_\Sp (3,2,3)$, and the last word starts with an odd letter, we know that this triangular planar history is associated to the null matching.
\end{example}

There are other reading orders for the crossings of $T$ that produce a Hecke equivalent word.

\begin{lemma}
\label{lemma:flexiblereadingword}
    Let $T$ be a triangular planar history with $p$ crossings.  
    Label these crossings with the integers $1,2,\ldots, p$ so that, when reading along each pipe starting from the diagonal and working towards the bottom of the triangle, the labels are in increasing order.
    Let $b_i$ be the number of pipes that pass left of a point immediately below and left of the $i$th crossing.
    Then $\bfb=(b_1,\ldots,b_p) \sim \bfa_T$.
\end{lemma}
\begin{proof}
    This follows from \cite{Weigandt.Bases}*{Lem.~3.2}.
\end{proof}

There are certain fundamental local moves on planar histories that preserve equivalence classes of Hecke words (and hence symplectic Hecke words).  

\begin{lemma}
\label{lemma:basiclocalreplacements}
    Fix triangular planar histories $T$ and $T'$.  If $T\mapsto T'$ by making any of the local replacements pictured below, then $\bfa_T \sim \bfa_{T'}$.  In particular, $\demprodfpf(T)=\demprodfpf(T')$.
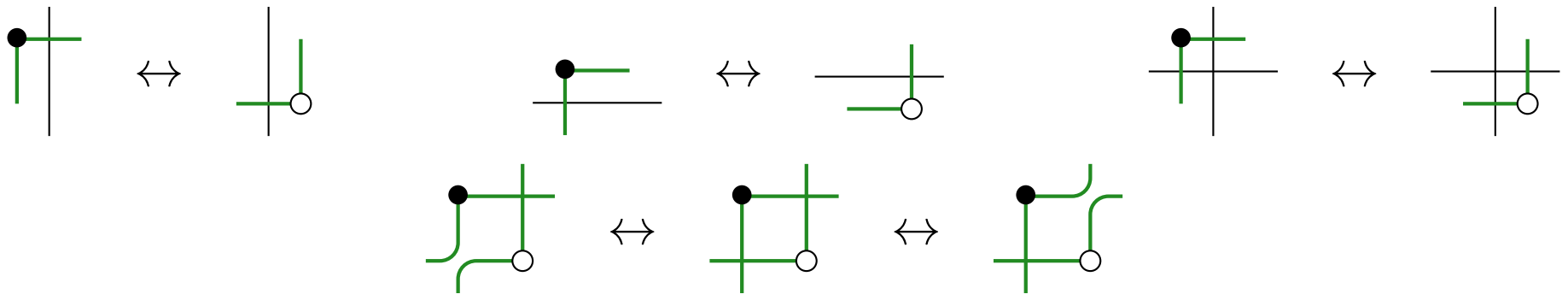
    \[\begin{tikzpicture}[scale=1]
  \draw[] (0,-0.5) -- (0,0.5);
  \draw[thick,ForestGreen] (-.25,-.25)--(-0.25,.25) -- (.25,.25);
  \node at (-.25,.25) {$\bullet$};
\end{tikzpicture}
  \hspace{1em}
\raisebox{1em}{$\leftrightarrow$}
\hspace{1em}
\begin{tikzpicture}[scale=1]
  \draw[] (0,-0.5) -- (0,0.5);
  \draw[thick,ForestGreen] (-.25,-.25)--(0.25,-.25) -- (.25,.25);
  \node at (.25,-.25) {$\whitecirc$};
\end{tikzpicture}
\hspace{4em}
\begin{tikzpicture}[scale=1]
  \draw[] (-0.5,0) -- (0.5,0);
  \draw[thick,ForestGreen] (-.25,-.25)--(-0.25,.25) -- (.25,.25);
  \node at (-.25,.25) {$\bullet$};
\end{tikzpicture}
  \hspace{1em}
\raisebox{1em}{$\leftrightarrow$}
\hspace{1em}
\begin{tikzpicture}[scale=1]
  \draw[] (-0.5,0) -- (0.5,0);
  \draw[thick,ForestGreen] (-.25,-.25)--(0.25,-.25) -- (.25,.25);
  \node at (.25,-.25) {$\whitecirc$};
\end{tikzpicture}
\hspace{4em}
\begin{tikzpicture}[scale=1]
  \draw[] (-0.5,0) -- (0.5,0);
    \draw[] (0,-0.5) -- (0,0.5);
  \draw[thick,ForestGreen] (-.25,-.25)--(-0.25,.25) -- (.25,.25);
  \node at (-.25,.25) {$\bullet$};
\end{tikzpicture}
  \hspace{1em}
\raisebox{1em}{$\leftrightarrow$}
\hspace{1em}
\begin{tikzpicture}[scale=1]
  \draw[] (-0.5,0) -- (0.5,0);
    \draw[] (0,-0.5) -- (0,0.5);
  \draw[thick,ForestGreen] (-.25,-.25)--(0.25,-.25) -- (.25,.25);
  \node at (.25,-.25) {$\whitecirc$};
\end{tikzpicture}
\]
\[\begin{tikzpicture}[scale=1]
  \draw[thick,ForestGreen,rounded corners](-.5,-.25)--(-.25,-.25)[sharp corners]-- (-0.25,.25)  -- (.25,.25)--(.5,.25);
\draw[thick,ForestGreen,rounded corners] (-.25,-.5)--(-.25,-.25)[sharp corners]--(0.25,-.25) -- (.25,.5);
  \node at (-.25,.25) {$\bullet$};
    \node at (.25,-.25) {$\whitecirc$};
\end{tikzpicture}
  \hspace{1em}
\raisebox{1em}{$\leftrightarrow$}
\hspace{1em}
\begin{tikzpicture}[scale=1]
  \draw[thick,ForestGreen] (-.25,-.5)--(-0.25,.25) -- (.25,.25)--(.5,.25);
\draw[thick,ForestGreen] (-.5,-.25)--(0.25,-.25) -- (.25,.5);
  \node at (-.25,.25) {$\bullet$};
    \node at (.25,-.25) {$\whitecirc$};
\end{tikzpicture}
  \hspace{1em}
\raisebox{1em}{$\leftrightarrow$}
\hspace{1em}
\begin{tikzpicture}[scale=1]
  \draw[thick,ForestGreen,rounded corners](.5,.25)--(.25,.25)[sharp corners]-- (.25,-.25)  -- (-.25,-.25)--(-.5,-.25);
  \draw[thick,ForestGreen,rounded corners] (.25,.5)--(.25,.25)[sharp corners]--(-.25,.25) -- (-.25,-.5);
  \node at (.25,-.25) {$\whitecirc$};
  \node at (-.25,.25) {$\bullet$};
\end{tikzpicture}
\]
\end{lemma}

\begin{proof}
    The given replacements on rectangular planar histories preserve equivalence classes of Hecke words (see, e.g., \cite{weigandt2021bumpless}*{Lem.~2.1}).
    Thus, these replacements also preserve equivalence classes of symplectic Hecke words for triangular planar histories.
\end{proof}

The moves from \cref{lemma:basiclocalreplacements} can be composed to produce more complicated modifications that also preserve equivalence classes of Hecke words, such as the droop moves of \cite{lam2021back}. 
We now introduce an additional move on triangular planar histories, called a \emph{skew droop move}, that preserves equivalence classes of symplectic Hecke words.
Droop moves on triangular bumpless pipe dreams, which are more restrictive, will be defined in \cref{subsection:droops}.

A \newword{skew droop} is a local replacement of the form pictured below.
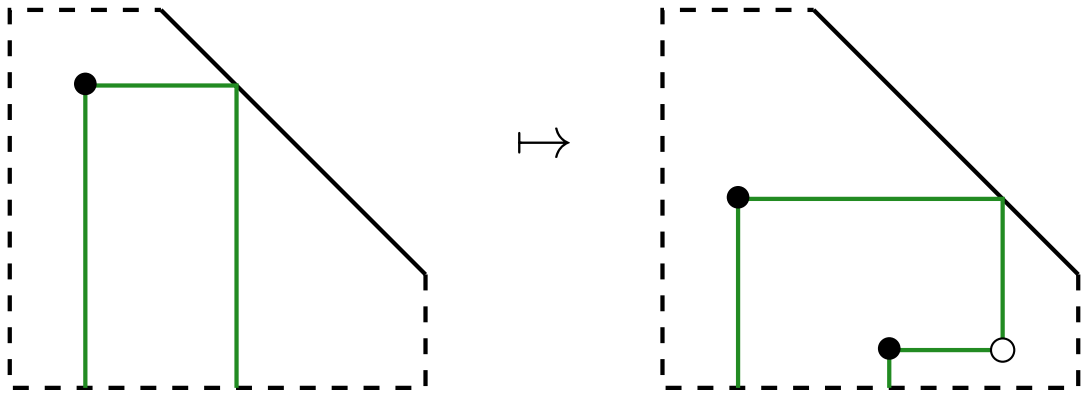
\[
\begin{tikzpicture}[scale=1]
  \draw[thick] (-.5,3.5) -- (1.25,1.75);
  \draw[thick,dashed](1.25,1.75) --(1.25,1)--(-1.5,1)--(-1.5,3.5)--(-.5,3.5);
  \draw[thick, ForestGreen]
    (-1,1) -- (-1,3) -- (0,3) -- (0,1);
  \node at (-1,3) {$\bullet$};
\end{tikzpicture}
  \hspace{1.5em}
\raisebox{4em}{$\mapsto$}
\hspace{1.5em}
\begin{tikzpicture}[scale=1]
  \draw[thick] (-.5,3.5) -- (1.25,1.75);
  \draw[thick,dashed](1.25,1.75) --(1.25,1)--(-1.5,1)--(-1.5,3.5)--(-.5,3.5);
  \draw[thick, ForestGreen]
    (-1,1) -- (-1,2.25) -- (.75,2.25)--(.75,1.25) --(0,1.25)-- (0,1);
  \node at (-1,2.25) {$\bullet$};
  \node at (0,1.25) {$\bullet$};
  \node at (.75,1.25) {$\whitecirc$};
\end{tikzpicture}
\]

Within this region, there may be other unpictured pipes, though none are allowed to contain bends of the form $\seelbow$ or $\nwelbow$. 
However, the green pipe may move past some number of $\yestri$ tiles on the diagonal.  See the image below for an example.
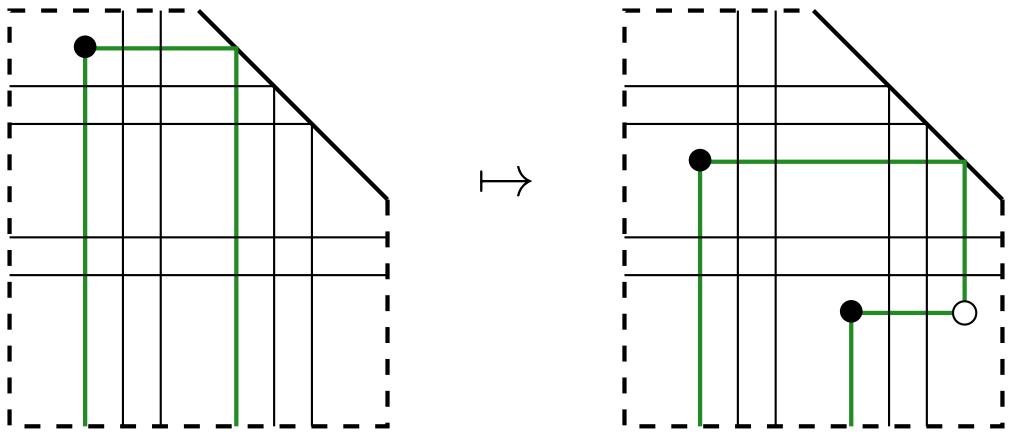
\[
\begin{tikzpicture}[scale=1]
  \draw[thick] (-.25,3.25) -- (1,2);
  \draw[thick,dashed] (1,2) -- (1,0.5) -- (-1.5,0.5) -- (-1.5,3.25) -- (-.25,3.25);
  \draw[thick, ForestGreen] (-1,0.5) -- (-1,3) -- (0,3) -- (0,0.5);
  \draw[] (-1.5,2.5) -- (.5,2.5) -- (.5,.5);     
  \draw[] (-1.5,2.75) -- (.25,2.75) -- (.25,.5); 
  \draw (-1.5,1.5) -- (1,1.5);
  \draw (-1.5,1.75) -- (1,1.75);
  \draw (-.75,3.25) -- (-.75,.5);   
  \draw (-.5,3.25) -- (-.5,.5);    
  \node at (-1,3) {$\bullet$};
\end{tikzpicture}
\hspace{1.5em}
\raisebox{4em}{$\mapsto$}
\hspace{1.5em}
\begin{tikzpicture}[scale=1]
  \draw[thick] (-.25,3.25) -- (1,2);
  \draw[thick,dashed] (1,2) -- (1,0.5) -- (-1.5,0.5) -- (-1.5,3.25) -- (-.25,3.25);
  \draw[thick, ForestGreen] (-1,0.5) -- (-1,2.25) -- (.75,2.25) -- (.75,1.25) -- (0,1.25) -- (0,0.5);
  \draw[] (-1.5,2.5) -- (.5,2.5) -- (.5,.5);     
  \draw[] (-1.5,2.75) -- (.25,2.75) -- (.25,.5);
  \draw (-1.5,1.5) -- (1,1.5);
  \draw (-1.5,1.75) -- (1,1.75);
  \draw (-.75,3.25) -- (-.75,.5);   
  \draw (-.5,3.25) -- (-.5,.5);    
  \node at (-1,2.25) {$\bullet$};
  \node at (0,1.25) {$\bullet$};
  \node at (.75,1.25) {$\whitecirc$};
\end{tikzpicture}
\]
If $T\mapsto T'$ is a skew droop, we say that $T'\mapsto T$ is an \newword{inverse skew droop}.

\begin{lemma}
    Let $T,T'$ be triangular planar histories.  If $T\mapsto T'$ by a skew droop move, then $\demprodfpf(T)=\demprodfpf(T')$.
\end{lemma}
\begin{proof}
    If the skew droop does not involve any crossings, then the result automatically follows. 
    Otherwise, assume first that we perform a skew droop that passes one $\yestri$ on the diagonal.
    If there are no other pipes in the pictured region, then by \cref{lemma:flexiblereadingword}, we may choose reading words $\mathbf b$ for $T$ and $\mathbf b'$ for $T'$ so that $\mathbf b$ starts with $(2k,2k-1)$, $\mathbf b'$ starts with $(2k,2k+1)$, and $\mathbf b$ and $\mathbf b'$ agree after the first two entries.
    In this case $\mathbf b \sim_{\Sp}\mathbf b'$, and thus $\demprodfpf(T)=\demprodfpf(T')$.
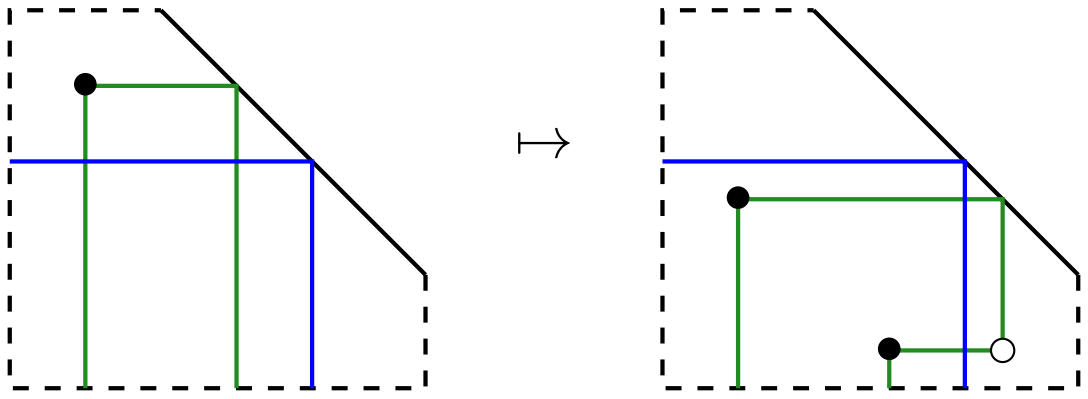
\[
\begin{tikzpicture}[scale=1]
  \draw[thick] (-.5,3.5) -- (1.25,1.75);
  \draw[thick,dashed] (1.25,1.75) -- (1.25,1) -- (-1.5,1) -- (-1.5,3.5) -- (-.5,3.5);
  \draw[thick, ForestGreen] (-1,1) -- (-1,3) -- (0,3) -- (0,1);
  \draw[thick,blue](-1.5,2.5)--(.5,2.5)--(.5,1);
  \node at (-1,3) {$\bullet$};
\end{tikzpicture}
  \hspace{1.5em}
\raisebox{4em}{$\mapsto$}
  \hspace{1.5em}
\begin{tikzpicture}[scale=1]
  \draw[thick] (-.5,3.5) -- (1.25,1.75);
  \draw[thick,dashed] (1.25,1.75) -- (1.25,1) -- (-1.5,1) -- (-1.5,3.5) -- (-.5,3.5);
  \draw[thick, ForestGreen] (-1,1) -- (-1,2.25) -- (.75,2.25)--(.75,1.25) --(0,1.25)-- (0,1);
  \draw[thick,blue](-1.5,2.5)--(.5,2.5)--(.5,1);
  \node at (-1,2.25) {$\bullet$};
  \node at (0,1.25) {$\bullet$};
  \node at (.75,1.25) {$\whitecirc$};
\end{tikzpicture}
\]

Now, if there are additional horizontal and vertical pipes that cross through the region, we may apply a sequence of the moves from \cref{lemma:basiclocalreplacements}, as well as a skew droop that involves no horizontal or vertical pipes, for which the moves compose to the desired skew droop.  See the images below.
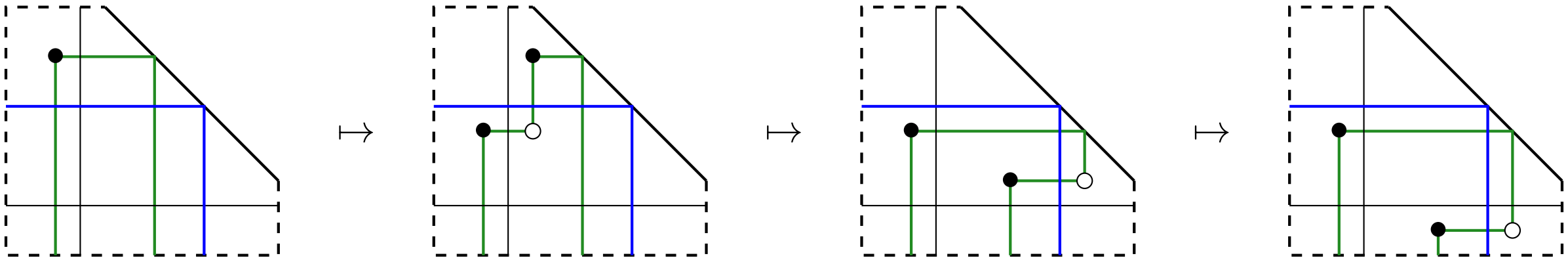
\[
\begin{tikzpicture}[scale=1]
  \draw[thick] (-.5,3.5) -- (1.25,1.75);
  \draw[thick,dashed] (1.25,1.75) -- (1.25,1) -- (-1.5,1) -- (-1.5,3.5) -- (-.5,3.5);
  \draw[thick, ForestGreen] (-1,1) -- (-1,3) -- (0,3) -- (0,1);
  \draw(-1.5,1.5)--(1.25,1.5); 
  \draw(-.75,3.5)--(-.75,1);
  \draw[thick,blue](-1.5,2.5)--(.5,2.5)--(.5,1);
  \node at (-1,3) {$\bullet$};
\end{tikzpicture}
  \hspace{1.5em}
\raisebox{3em}{$\mapsto$}
\hspace{1.5em}
\begin{tikzpicture}[scale=1]
  \draw[thick] (-.5,3.5) -- (1.25,1.75);
  \draw[thick,dashed] (1.25,1.75) -- (1.25,1) -- (-1.5,1) -- (-1.5,3.5) -- (-.5,3.5);
  \draw[thick, ForestGreen] (-1,1) -- (-1,2.25) -- (-.5,2.25)--(-.5,3) --(0,3)-- (0,1);
  \draw(-1.5,1.5)--(1.25,1.5);
  \draw(-.75,3.5)--(-.75,1);
  \draw[thick,blue](-1.5,2.5)--(.5,2.5)--(.5,1);
  \node at (-1,2.25) {$\bullet$};
  \node at (-.5,3) {$\bullet$};
  \node at (-.5,2.25) {$\whitecirc$};
\end{tikzpicture}
  \hspace{1.5em}
\raisebox{3em}{$\mapsto$}
\hspace{1.5em}
\begin{tikzpicture}[scale=1]
  \draw[thick] (-.5,3.5) -- (1.25,1.75);
  \draw[thick,dashed] (1.25,1.75) -- (1.25,1) -- (-1.5,1) -- (-1.5,3.5) -- (-.5,3.5);
  \draw[thick, ForestGreen] (-1,1) -- (-1,2.25) -- (.75,2.25)--(.75,1.75) --(0,1.75)-- (0,1);
  \draw(-1.5,1.5)--(1.25,1.5); 
  \draw(-.75,3.5)--(-.75,1);
  \draw[thick,blue](-1.5,2.5)--(.5,2.5)--(.5,1);
  \node at (-1,2.25) {$\bullet$};
  \node at (0,1.75) {$\bullet$};
  \node at (.75,1.75) {$\whitecirc$};
\end{tikzpicture}
  \hspace{1.5em}
\raisebox{3em}{$\mapsto$}
\hspace{1.5em}
\begin{tikzpicture}[scale=1]
  \draw[thick] (-.5,3.5) -- (1.25,1.75);
  \draw[thick,dashed] (1.25,1.75) -- (1.25,1) -- (-1.5,1) -- (-1.5,3.5) -- (-.5,3.5);
  \draw[thick, ForestGreen] (-1,1) -- (-1,2.25) -- (.75,2.25)--(.75,1.25) --(0,1.25)-- (0,1);
  \draw(-1.5,1.5)--(1.25,1.5); 
  \draw(-.75,3.5)--(-.75,1);
  \draw[thick,blue](-1.5,2.5)--(.5,2.5)--(.5,1);
  \node at (-1,2.25) {$\bullet$};
  \node at (0,1.25) {$\bullet$};
  \node at (.75,1.25) {$\whitecirc$};
\end{tikzpicture}
\]

Otherwise, if we are performing a skew droop that passes by multiple $\yestri$'s, it is enough to note that this move can be interpreted as the composition of moves where we do a skew droop past one $\yestri$ at a time, along with ordinary droops, as pictured below.
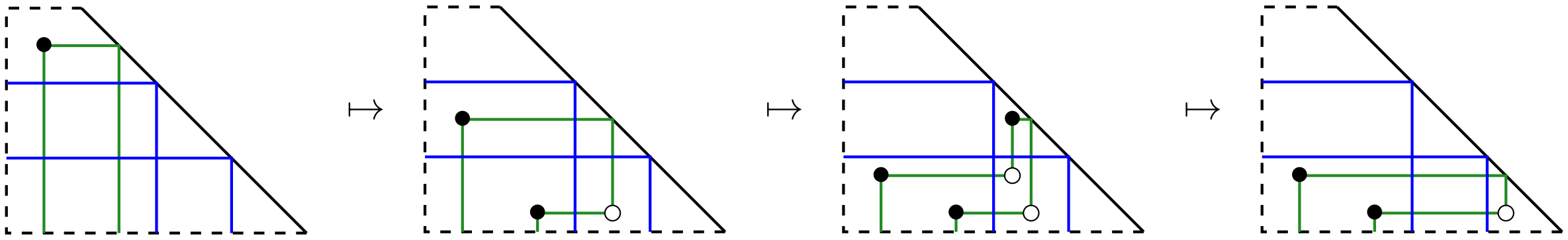
\[
\begin{tikzpicture}[scale=.75]
  \draw[thick] (-.5,3.5) -- (2.5,.5);
  \draw[thick,dashed] (2.5,.5) -- (-1.5,.5) -- (-1.5,3.5) -- (-.5,3.5);
  \draw[thick, ForestGreen] (-1,0.5) -- (-1,3) -- (0,3) -- (0,0.5);
  \draw[thick,blue](-1.5,2.5)--(.5,2.5)--(.5,.5);
  \draw[thick,blue](-1.5,1.5)--(1.5,1.5)--(1.5,.5);
  \node at (-1,3) {$\bullet$};
\end{tikzpicture}
  \hspace{1em}
\raisebox{3em}{$\mapsto$}
\hspace{1em}
\begin{tikzpicture}[scale=.75]
  \draw[thick] (-.5,3.5) -- (2.5,.5); 
  \draw[thick,dashed] (2.5,.5) -- (-1.5,.5) -- (-1.5,3.5) -- (-.5,3.5);
  \draw[thick, ForestGreen] (-1,0.5) -- (-1,2) -- (1,2)--(1,.75) --(0,.75)-- (0,0.5);
  \draw[thick,blue](-1.5,2.5)--(.5,2.5)--(.5,.5);
  \draw[thick,blue](-1.5,1.5)--(1.5,1.5)--(1.5,.5);
  \node at (-1,2) {$\bullet$};
  \node at (0,.75) {$\bullet$};
  \node at (1,.75) {$\whitecirc$};
\end{tikzpicture}
  \hspace{1em}
\raisebox{3em}{$\mapsto$}
\hspace{1em}
\begin{tikzpicture}[scale=.75]
  
  \draw[thick] (-.5,3.5) -- (2.5,.5); 
  \draw[thick,dashed] (2.5,.5) -- (-1.5,.5) -- (-1.5,3.5) -- (-.5,3.5);

  \draw[thick, ForestGreen] (-1,0.5) -- (-1,1.25) --(.75,1.25)--(.75,2) -- (1,2)--(1,.75) --(0,.75)-- (0,0.5);
  \draw[thick,blue](-1.5,2.5)--(.5,2.5)--(.5,.5);
  \draw[thick,blue](-1.5,1.5)--(1.5,1.5)--(1.5,.5);
  \node at (-1,1.25) {$\bullet$};
  \node at (.75,1.25) {$\whitecirc$};
  \node at (.75,2) {$\bullet$};
  \node at (0,.75) {$\bullet$};
  \node at (1,.75) {$\whitecirc$};
\end{tikzpicture}
\hspace{1em}
\raisebox{3em}{$\mapsto$}
\hspace{1em}
\begin{tikzpicture}[scale=.75]
  
  \draw[thick] (-.5,3.5) -- (2.5,.5); 
  \draw[thick,dashed] (2.5,.5) -- (-1.5,.5) -- (-1.5,3.5) -- (-.5,3.5);
  \draw[thick, ForestGreen] (-1,0.5) -- (-1,1.25) -- (1.75,1.25)--(1.75,.75) --(0,.75)-- (0,0.5);
  \draw[thick,blue](-1.5,2.5)--(.5,2.5)--(.5,.5);
  \draw[thick,blue](-1.5,1.5)--(1.5,1.5)--(1.5,.5);
  \node at (-1,1.25) {$\bullet$};
  \node at (0,.75) {$\bullet$};
  \node at (1.75,.75) {$\whitecirc$};
\end{tikzpicture}
\]
All of these moves preserve equivalence classes of symplectic Hecke words, so the result follows.    
\end{proof}

By interpreting $\seelbow$'s as $1$'s and $\nwelbow$'s as $-1$'s, we may readily import terms from ASMs for use in the setting of triangular planar histories.  In particular, we define neighbors of $-1$'s in the same way. It is possible in this setting that a $-1$ will not have a top neighbor.  When it does, if the $-1$ is removable, then we may perform inflation in the usual way directly on a triangular planar history, and this produces a new triangular planar history.

\begin{lemma}
\label{lemma:easy_symplectic_inflation}
For $T,T'$ triangular planar histories,  if $T\mapsto T'$ by inflation then $\demprodfpf(T)=\demprodfpf(T')$.
\end{lemma}
\begin{proof}
This follows from \cite{weigandt2021bumpless}*{Lem.~4.2}, which applies unchanged to triangular planar histories.
\end{proof}

The same applies for symplectic inflation.

\begin{proposition}
\label{prop:symplecticinflationpreservesdemprod}
    Let $A,A'\in \asmsk(2n)$ so that $A\mapsto A'$ by symplectic inflation at $(b,a)$ with $a<b$.  Then $\demprodfpf(\Psi(A))=\demprodfpf(\Psi(A'))$.
\end{proposition}
\begin{proof}
    If all neighbors of $(b,a)$ sit below the diagonal, we are in the case of \cref{lemma:easy_symplectic_inflation}, so assume not.
    In particular, the top neighbor of $(b,a)$ sits above the diagonal.
    Since $(b,a)$ and $(a,b)$ have no common neighbors, every lower neighbor of $(b,a)$ lies in a row greater than $a$.

    We first focus on the lower neighbors of $(b,a)$.  By doing sequences of the moves from \cref{lemma:basiclocalreplacements}, we may reduce to the case where $(b,a)$ has exactly one lower neighbor as illustrated:
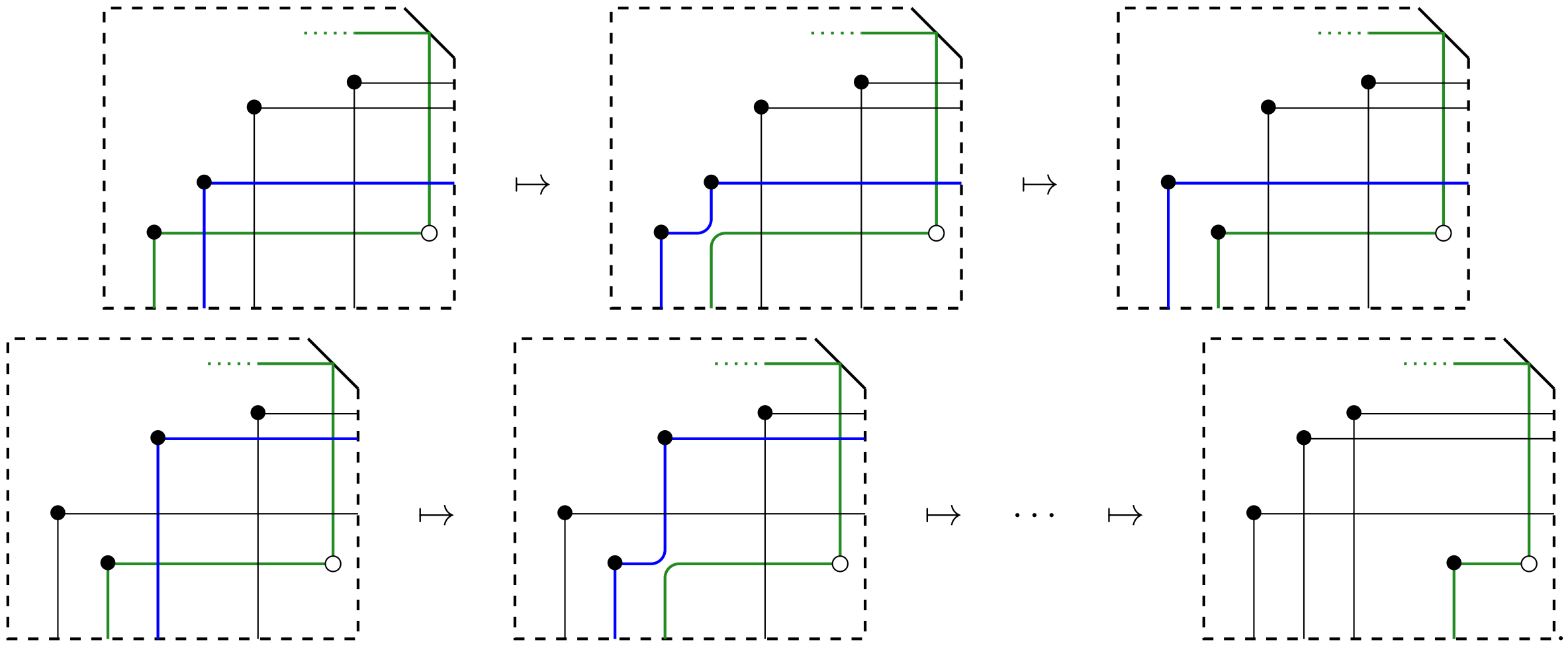
\[
\begin{tikzpicture}[scale=1]
  \draw[thick] (.5,3.5) -- (1,3);
  \draw[thick,dashed] (1,3) -- (1,.5) -- (-2.5,.5) -- (-2.5,3.5) -- (.5,3.5);
   \draw[thick, ForestGreen,dotted](-.5,3.25)--(0,3.25);
  \draw[thick, ForestGreen]
   (0,3.25)-- (.75,3.25)--(.75,1.25) --(-2,1.25)-- (-2,0.5);
     \draw[](0,.5)--(0,2.75)--(1,2.75);
     \draw[](-1,.5)--(-1,2.5)--(1,2.5);
     \draw[thick, blue](-1.5,.5)--(-1.5,1.75)--(1,1.75);
  \node at (0,2.75) {$\bullet$};
  \node at (-1,2.5) {$\bullet$};
  \node at (-1.5,1.75) {$\bullet$};
\node at (-2,1.25) {$\bullet$};
\node at (.75,1.25) {$\whitecirc$};
\end{tikzpicture}
\hspace{1.5em}
\raisebox{3em}{$\mapsto$}
\hspace{1.5em}
\begin{tikzpicture}[scale=1]
  \draw[thick] (.5,3.5) -- (1,3);
  \draw[thick,dashed] (1,3) -- (1,.5) -- (-2.5,.5) -- (-2.5,3.5) -- (.5,3.5);
   \draw[thick, ForestGreen,dotted](-.5,3.25)--(0,3.25);
  \draw[thick, ForestGreen,rounded corners]
     (-1.5,0.5) -- (-1.5,1.25)[sharp corners] -- (.75,1.25) -- (.75,3.25)--(0,3.25);
     \draw[](0,.5)--(0,2.75)--(1,2.75);
     \draw[](-1,.5)--(-1,2.5)--(1,2.5);
     \draw[thick, blue](-2,.5)--(-2,1.25)[rounded corners]--(-1.5,1.25)[sharp corners]--(-1.5,1.75)--(1,1.75);
  \node at (0,2.75) {$\bullet$};
  \node at (-1,2.5) {$\bullet$};
  \node at (-1.5,1.75) {$\bullet$};
\node at (-2,1.25) {$\bullet$};
\node at (.75,1.25) {$\whitecirc$};
\end{tikzpicture}
\hspace{1.5em}
\raisebox{3em}{$\mapsto$}
\hspace{1.5em}
\begin{tikzpicture}[scale=1]
  \draw[thick] (.5,3.5) -- (1,3);
  \draw[thick,dashed] (1,3) -- (1,.5) -- (-2.5,.5) -- (-2.5,3.5) -- (.5,3.5);
   \draw[thick, ForestGreen,dotted](-.5,3.25)--(0,3.25);
  \draw[thick, ForestGreen]
     (-1.5,0.5) -- (-1.5,1.25)[sharp corners] -- (.75,1.25) -- (.75,3.25)--(0,3.25);
     \draw[](0,.5)--(0,2.75)--(1,2.75);
     \draw[](-1,.5)--(-1,2.5)--(1,2.5);
     \draw[thick, blue](-2,.5)--(-2,1.75)--(1,1.75);
  \node at (0,2.75) {$\bullet$};
  \node at (-1,2.5) {$\bullet$};
  \node at (-1.5,1.25) {$\bullet$};
\node at (-2,1.75) {$\bullet$};
\node at (.75,1.25) {$\whitecirc$};
\end{tikzpicture}
\]
\[
\begin{tikzpicture}[scale=1]
  \draw[thick] (.5,3.5) -- (1,3);
  \draw[thick,dashed] (1,3) -- (1,.5) -- (-2.5,.5) -- (-2.5,3.5) -- (.5,3.5);
   \draw[thick, ForestGreen,dotted](-.5,3.25)--(0,3.25);
  \draw[thick, ForestGreen]
     (-1.5,0.5) -- (-1.5,1.25)[sharp corners] -- (.75,1.25) -- (.75,3.25)--(0,3.25);
     \draw[](0,.5)--(0,2.75)--(1,2.75);
     \draw[thick, blue](-1,.5)--(-1,2.5)--(1,2.5);
     \draw[](-2,.5)--(-2,1.75)--(1,1.75);
  \node at (0,2.75) {$\bullet$};
  \node at (-1,2.5) {$\bullet$};
  \node at (-1.5,1.25) {$\bullet$};
\node at (-2,1.75) {$\bullet$};
\node at (.75,1.25) {$\whitecirc$};
\end{tikzpicture}
\hspace{1.5em}
\raisebox{3em}{$\mapsto$}
\hspace{1.5em}
\begin{tikzpicture}[scale=1]
  \draw[thick] (.5,3.5) -- (1,3);
  \draw[thick,dashed] (1,3) -- (1,.5) -- (-2.5,.5) -- (-2.5,3.5) -- (.5,3.5);
   \draw[thick, ForestGreen,dotted](-.5,3.25)--(0,3.25);
  \draw[thick, ForestGreen,rounded corners]
     (-1,0.5) -- (-1,1.25)[sharp corners] -- (.75,1.25) -- (.75,3.25)--(0,3.25);
     \draw[](0,.5)--(0,2.75)--(1,2.75);
     \draw[thick, blue](-1.5,.5)--(-1.5,1.25)[rounded corners]--(-1,1.25)[sharp corners]--(-1,2.5)--(1,2.5);
     \draw[](-2,.5)--(-2,1.75)--(1,1.75);
  \node at (0,2.75) {$\bullet$};
  \node at (-1,2.5) {$\bullet$};
  \node at (-1.5,1.25) {$\bullet$};
\node at (-2,1.75) {$\bullet$};
\node at (.75,1.25) {$\whitecirc$};
\end{tikzpicture}
\hspace{1.5em}
\raisebox{3em}{$\mapsto \quad \cdots \quad \mapsto$}
\hspace{1.5em}
\begin{tikzpicture}[scale=1]
  \draw[thick] (.5,3.5) -- (1,3);
  \draw[thick,dashed] (1,3) -- (1,.5) -- (-2.5,.5) -- (-2.5,3.5) -- (.5,3.5);
   \draw[thick, ForestGreen,dotted](-.5,3.25)--(0,3.25);
  \draw[thick, ForestGreen]
     (0,3.25)--(.75,3.25)--(.75,1.25) --(0,1.25)-- (0,0.5);
     \draw[](-1,.5)--(-1,2.75)--(1,2.75);
     \draw[](-1.5,.5)--(-1.5,2.5)--(1,2.5);
     \draw[](-2,.5)--(-2,1.75)--(1,1.75);
  \node at (0,1.25) {$\bullet$};
  \node at (-1,2.75) {$\bullet$};
  \node at (-1.5,2.5) {$\bullet$};
\node at (-2,1.75) {$\bullet$};
\node at (.75,1.25) {$\whitecirc$};
\end{tikzpicture}.
\]

Now we are in the case where $(b,a)$ has a single lower neighbor; since it is then the left neighbor of $(b,a)$, it lies in row $b$, say at $(b,c)$.
Here $c$ is the column of the topmost lower neighbor in the original $A$.
Because $(b,a)$ did not share neighbors with $(a,b)$ in $A$, the partners of the upper neighbors occur strictly to the left of column $c$.
As such, we may perform an inverse skew droop to modify the pipe that passes through $(b,a)$ in such a way that we only modify pipes to the right of the partners of the upper neighbors of $(b,a)$.
Specifically, we reroute this pipe so that it passes vertically in column $c$, until it reflects at $(c,c)$, and the pipe that it meets at this reflection point travels horizontally in row $c$ to the left, then bends downwards at $(c,c-\epsilon)$ for some $\epsilon>0$ sufficiently small.
In our figure, any choice $(c,c-\epsilon)$ in the small dashed box suffices.
 
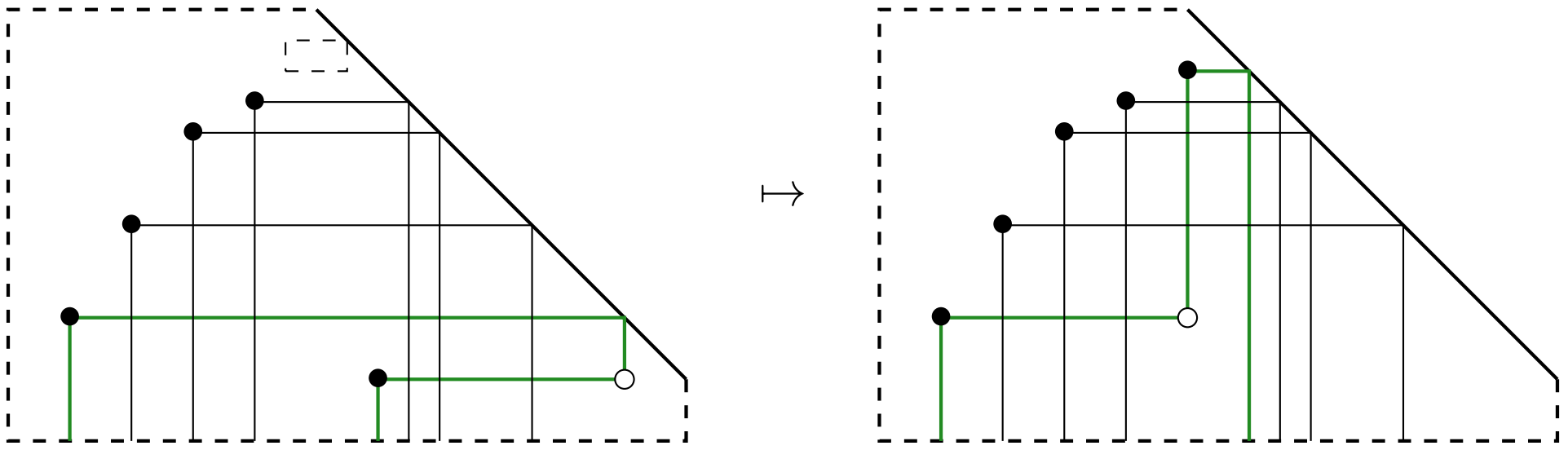
\[\begin{tikzpicture}[scale=1]
  \draw[thick] (-.5,3.5) -- (2.5,.5);
  \draw[thick,dashed] (2.5,.5) -- (2.5,0) -- (-3,0) -- (-3,3.5) -- (-.5,3.5);
  \draw[thick, ForestGreen]
     (-2.5,0)--(-2.5,1)--(2,1)--(2,.5) --(0,.5)-- (0,0);
  \draw[](-1,0)--(-1,2.75)--(.25,2.75)--(.25,0);
  \draw[](-1.5,0)--(-1.5,2.5)--(.5,2.5)--(.5,0);
  \draw[](-2,0)--(-2,1.75)--(1.25,1.75)--(1.25,0);
  \node at (0,.5) {$\bullet$};
  \draw[dashed] (-.75,3) -- (-.25,3) -- (-.25,3.25) -- (-.75,3.25) -- (-.75,3);
  \node at (-1,2.75) {$\bullet$};
  \node at (-1.5,2.5) {$\bullet$};
  \node at (-2,1.75) {$\bullet$};
  \node at (-2.5,1) {$\bullet$};
  \node at (2,.5) {$\whitecirc$};
\end{tikzpicture}
\hspace{1.5em}
\raisebox{5em}{$\mapsto$}
\hspace{1.5em}
\begin{tikzpicture}[scale=1]
  \draw[thick] (-.5,3.5) -- (2.5,.5);
  \draw[thick,dashed] (2.5,.5) -- (2.5,0) -- (-3,0) -- (-3,3.5) -- (-.5,3.5);
  \draw[thick, ForestGreen]
     (-2.5,0)--(-2.5,1)--(-.5,1)--(-.5,3)--(0,3)-- (0,0);
  \draw[](-1,0)--(-1,2.75)--(.25,2.75)--(.25,0);
  \draw[](-1.5,0)--(-1.5,2.5)--(.5,2.5)--(.5,0);
  \draw[](-2,0)--(-2,1.75)--(1.25,1.75)--(1.25,0);
  \node at (-.5,3) {$\bullet$};
  \node at (-1,2.75) {$\bullet$};
  \node at (-1.5,2.5) {$\bullet$};
  \node at (-2,1.75) {$\bullet$};
  \node at (-2.5,1) {$\bullet$};
  \node at (-.5,1) {$\whitecirc$};
\end{tikzpicture}
\]
The new $\nwelbow$ is in position $(a,c-\epsilon)$, and is removable.
Thus, we may perform inflation at $(a,c-\epsilon)$, resulting in the replacements pictured below.
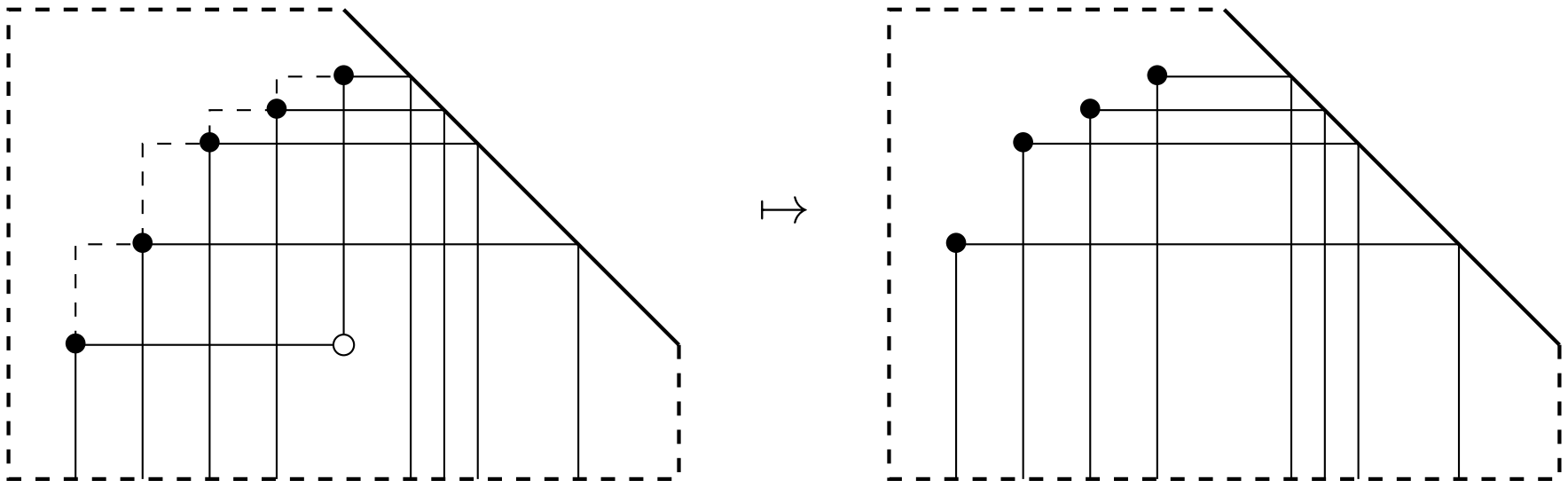
\[
\begin{tikzpicture}[scale=1]
  \draw[thick] (-.5,3.5) -- (2,1);
  \draw[thick,dashed] (2,1) -- (2,0) -- (-3,0) -- (-3,3.5) -- (-.5,3.5);
  \draw[]
     (-2.5,0)--(-2.5,1)--(-.5,1)--(-.5,3)--(0,3)-- (0,0);
  \draw[](-1,0)--(-1,2.75)--(.25,2.75)--(.25,0);
  \draw[](-1.5,0)--(-1.5,2.5)--(.5,2.5)--(.5,0);
  \draw[](-2,0)--(-2,1.75)--(1.25,1.75)--(1.25,0);
  \draw[dashed](-2.5,1)--(-2.5,1.75)--(-2,1.75)--(-2,2.5)--(-1.5,2.5)--(-1.5,2.75)--(-1,2.75)--(-1,3)--(-.5,3);
  \node at (-.5,3) {$\bullet$};
  \node at (-1,2.75) {$\bullet$};
  \node at (-1.5,2.5) {$\bullet$};
  \node at (-2,1.75) {$\bullet$};
  \node at (-2.5,1) {$\bullet$};
  \node at (-.5,1) {$\whitecirc$};
\end{tikzpicture}
\hspace{1.5em}
\raisebox{5em}{$\mapsto$}
\hspace{1.5em}
\begin{tikzpicture}[scale=1]
  \draw[thick] (-.5,3.5) -- (2,1);
  \draw[thick,dashed] (2,1) -- (2,0) -- (-3,0) -- (-3,3.5) -- (-.5,3.5);
  \draw[]
     (-1,0)--(-1,3)--(0,3)-- (0,0);
  \draw[](-1.5,0)--(-1.5,2.75)--(.25,2.75)--(.25,0);
  \draw[](-2,0)--(-2,2.5)--(.5,2.5)--(.5,0);
  \draw[](-2.5,0)--(-2.5,1.75)--(1.25,1.75)--(1.25,0);
  \node at (-1,3) {$\bullet$};
  \node at (-1.5,2.75) {$\bullet$};
  \node at (-2,2.5) {$\bullet$};
  \node at (-2.5,1.75) {$\bullet$};
\end{tikzpicture}\]
Throughout this whole process, we have made a sequence of replacements within triangular planar histories that ultimately result in $\Psi(A')$.
Because each of these replacements preserves equivalence classes of symplectic Hecke words, the result follows.
\end{proof}

We now introduce two configurations in triangular planar histories that prevent them from being symplectic:
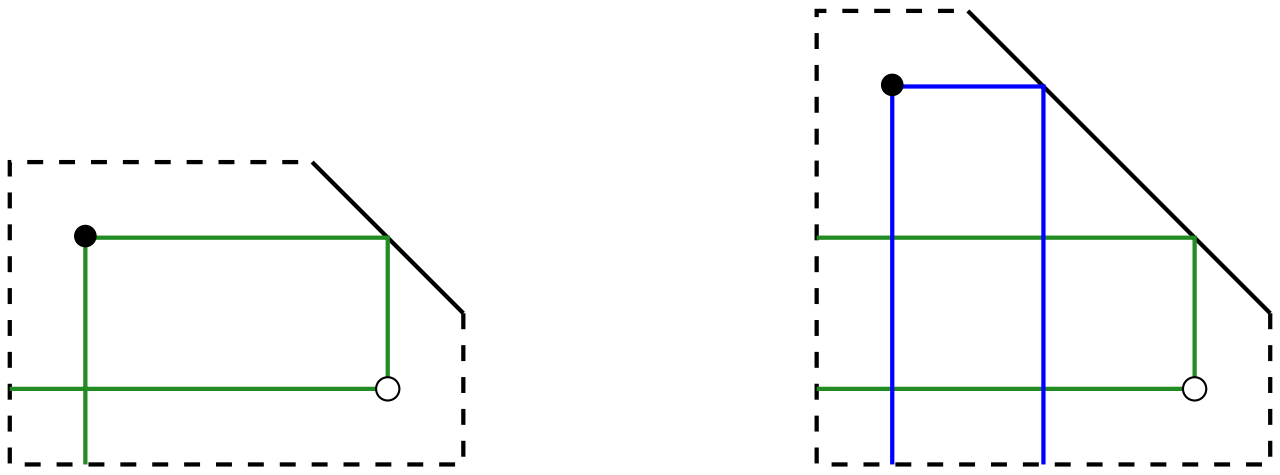
\[
\begin{tikzpicture}[scale=1]
  \draw[thick] (.5,2.5) -- (1.5,1.5);
  \draw[thick,dashed](1.5,1.5)--(1.5,.5)--(-1.5,.5)--(-1.5,2.5)--(.5,2.5);
\draw[thick, ForestGreen]
  (-1,.5) -- (-1,2) -- (1,2) -- (1,1) -- (-1.5,1);
\node at (1,1) {$\whitecirc$};
    \node at (-1,2) {$\bullet$};
\end{tikzpicture}
\hspace{6em}
\begin{tikzpicture}[scale=1]
  \draw[thick] (-.5,3.5) -- (1.5,1.5);
  \draw[thick,dashed](1.5,1.5)--(1.5,.5)--(-1.5,.5)--(-1.5,3.5)--(-.5,3.5);
  \draw[thick, ForestGreen]
    (-1.5,2) -- (1,2) -- (1,1) -- (-1.5,1);
  \draw[thick, blue]
    (-1,0.5) -- (-1,3) -- (0,3) -- (0,0.5);
  \node at (1,1) {$\whitecirc$};
  \node at (-1,3) {$\bullet$};
\end{tikzpicture}
\]
Call the configuration on the left a \newword{null cycle} and the configuration on the right \newword{overlapping horseshoes}.
Within the dashed regions, there may be other unpictured pipes, though none are allowed to contain bends of the form $\seelbow$ or $\nwelbow$.

\begin{lemma}
\label{lemma:basicbadconfigurations}
Let $T$ be a triangular planar history.
If $T$ contains a null cycle or overlapping horseshoes, then $\demprodfpf(T)=\varnothing$.
\end{lemma}

\begin{proof}
    Assume $T$ contains a null cycle.  Suppose first that there are no other pipes passing through the green rectangle.  
    The letter associated to the pictured crossing must be odd.  
    Because of this, using \cref{lemma:flexiblereadingword}, we can choose a reading word $\mathbf b$ for $T$ starting with an odd letter.  
    Thus $\mathbf{a}_T\sim\mathbf b$ and $\demprodfpf(\mathbf b)=\varnothing$, since $\mathbf b$ begins with an odd letter.
    Hence $\demprodfpf(T)=\varnothing$.
    
    Now, suppose that there are pipes that pass through the green rectangle, and these pipes do not bend within the rectangle. 
    By applying a sequence of the moves from \cref{lemma:basiclocalreplacements}, we may reduce to the case where no other pipes pass through the given rectangle.  
    See the images below.
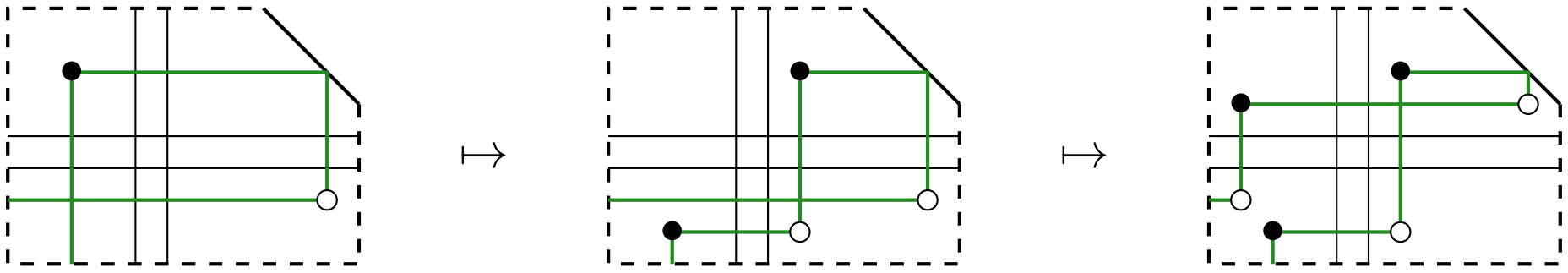
\[
\begin{tikzpicture}[scale=1]
  \draw[thick] (.5,2.5) -- (1.25,1.75);
  \draw[thick,dashed] (1.25,1.75)--(1.25,.5)--(-1.5,.5)--(-1.5,2.5)--(.5,2.5);
  \draw (-1.5,1.25)--(1.25,1.25);
  \draw (-1.5,1.5)--(1.25,1.5);
  \draw (-.5,2.5)--(-.5,.5);
  \draw (-.25,2.5)--(-.25,.5);
    \draw[thick, ForestGreen] (-1,.5) -- (-1,2) -- (1,2) -- (1,1) -- (-1.5,1);
    \node at (1,1) {$\whitecirc$};
  \node at (-1,2) {$\bullet$};
\end{tikzpicture}
\hspace{2em}
\raisebox{2em}{$\mapsto$}
\hspace{2em}
\begin{tikzpicture}[scale=1]
  \draw[thick] (.5,2.5) -- (1.25,1.75);
  \draw[thick,dashed] (1.25,1.75)--(1.25,.5)--(-1.5,.5)--(-1.5,2.5)--(.5,2.5);
    \draw (-1.5,1.25)--(1.25,1.25);
  \draw (-1.5,1.5)--(1.25,1.5);
  \draw (-.5,2.5)--(-.5,.5);
  \draw (-.25,2.5)--(-.25,.5);
    \draw[thick, ForestGreen] (-1,.5) -- (-1,.75)--(0,.75)--(0,2) -- (1,2) -- (1,1) -- (-1.5,1);
    \node at (1,1) {$\whitecirc$};
  \node at (-1,.75) {$\bullet$};
  \node at (0,.75) {$\whitecirc$};
  \node at (0,2) {$\bullet$};
\end{tikzpicture}
\hspace{2em}
\raisebox{2em}{$\mapsto$}
\hspace{2em}
\begin{tikzpicture}[scale=1]

  \draw[thick] (.5,2.5) -- (1.25,1.75);
  \draw[thick,dashed] (1.25,1.75)--(1.25,.5)--(-1.5,.5)--(-1.5,2.5)--(.5,2.5);

  \draw (-1.5,1.25)--(1.25,1.25);
  \draw (-1.5,1.5)--(1.25,1.5);
  \draw (-.5,2.5)--(-.5,.5);
  \draw (-.25,2.5)--(-.25,.5);
  \draw[thick, ForestGreen] (-1,.5) -- (-1,.75)--(0,.75)--(0,2) -- (1,2) -- (1,1.75) -- (-1.25,1.75)--(-1.25,1)--(-1.5,1);
  \node at (1,1.75) {$\whitecirc$};
  \node at (-1.25,1) {$\whitecirc$};
  \node at (-1.25,1.75) {$\bullet$};
  \node at (-1,.75) {$\bullet$};
  \node at (0,.75) {$\whitecirc$};
  \node at (0,2) {$\bullet$};
\end{tikzpicture}
\]

Now suppose $T$ contains overlapping horseshoes.
If no other pipes pass through the pictured region, we can choose a reading order on $T$ giving a word $\mathbf b$ with initial segment  $(2k,2k-1,2k+1,2k)$.
Applying the relations to this initial segment gives
\[
\begin{aligned}
(2k,2k-1,2k+1,2k)&\sim_\Sp(2k,2k+1,2k+1,2k)\\
&\sim_\Sp(2k,2k+1,2k)\sim_\Sp(2k+1,2k,2k+1).
\end{aligned}
\]
The resulting word begins with an odd letter, so its symplectic Demazure product is null.   Since reading words for $T$ in any admissible order are Hecke equivalent by \cref{lemma:flexiblereadingword}, $\demprodfpf(T)=\varnothing$, as well.

Now suppose we have some number of horizontal pipes passing through the green horseshoe, vertical pipes passing through the blue horseshoe, and some pipes which reflect between the blue reflection point and the green reflection point.  See the leftmost picture below. 
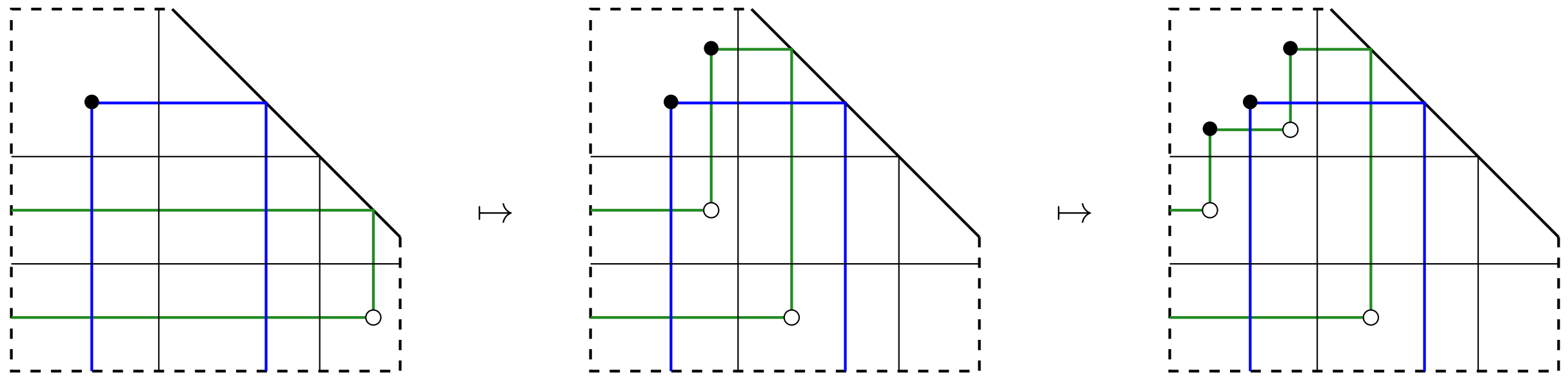
\[\begin{tikzpicture}[scale=.55]
  \draw[thick] (-.75,6.75) -- (3.5,2.5);
  \draw[thick,dashed] (3.5,2.5) -- (3.5,0) -- (-3.75,0) -- (-3.75,6.75) -- (-.75,6.75);
  \draw[thick, ForestGreen]
    (-3.75,3) -- (3,3) -- (3,1) -- (-3.75,1);
  \draw[thick, blue]
    (-2.25,0) -- (-2.25,5) -- (1,5) -- (1,0);
    \draw(-3.75,2)--(3.5,2);
    \draw(-1,6.75)--(-1,0);
    \draw(-3.75,4)--(2,4)--(2,0);
  \node at (3,1) {$\whitecirc$};
  \node at (-2.25,5) {$\bullet$};
\end{tikzpicture}
\hspace{2em}
\raisebox{4em}{$\mapsto$}
\hspace{2em}
\begin{tikzpicture}[scale=.55]
  \draw[thick] (-.75,6.75) -- (3.5,2.5);
  \draw[thick,dashed] (3.5,2.5) -- (3.5,0) -- (-3.75,0) -- (-3.75,6.75) -- (-.75,6.75);
  \draw[thick, ForestGreen]
    (-3.75,3)--(-1.5,3) -- (-1.5,6) --(0,6) -- (0,1) -- (-3.75,1);
  \draw[thick, blue]
    (-2.25,0) -- (-2.25,5) -- (1,5) -- (1,0);
    \draw(-3.75,2)--(3.5,2);
    \draw(-1,6.75)--(-1,0);
    \draw(-3.75,4)--(2,4)--(2,0);
  \node at (0,1) {$\whitecirc$};
    \node at (-1.5,3) {$\whitecirc$};
  \node at (-2.25,5) {$\bullet$};
    \node at (-1.5,6) {$\bullet$};
\end{tikzpicture}
\hspace{2em}
\raisebox{4em}{$\mapsto$}
\hspace{2em}
\begin{tikzpicture}[scale=.55]
  \draw[thick] (-.75,6.75) -- (3.5,2.5);
  \draw[thick,dashed] (3.5,2.5) -- (3.5,0) -- (-3.75,0) -- (-3.75,6.75) -- (-.75,6.75);
  \draw[thick, ForestGreen]
    (-3.75,3)--(-3,3)--(-3,4.5) --(-1.5,4.5)-- (-1.5,6) --(0,6) -- (0,1) -- (-3.75,1);
  \draw[thick, blue]
    (-2.25,0) -- (-2.25,5) -- (1,5) -- (1,0);
    \draw(-3.75,2)--(3.5,2);
    \draw(-1,6.75)--(-1,0);
    \draw(-3.75,4)--(2,4)--(2,0);
  \node at (0,1) {$\whitecirc$};
    \node at (-3,3) {$\whitecirc$};
  \node at (-2.25,5) {$\bullet$};
    \node at (-3,4.5) {$\bullet$};
  \node at (-1.5,4.5) {$\whitecirc$};
  \node at (-1.5,6) {$\bullet$};
\end{tikzpicture}
\]
Doing inflation at the topmost pictured $\whitecirc$ produces the configuration below, which contains a null cycle (see the green pipe).
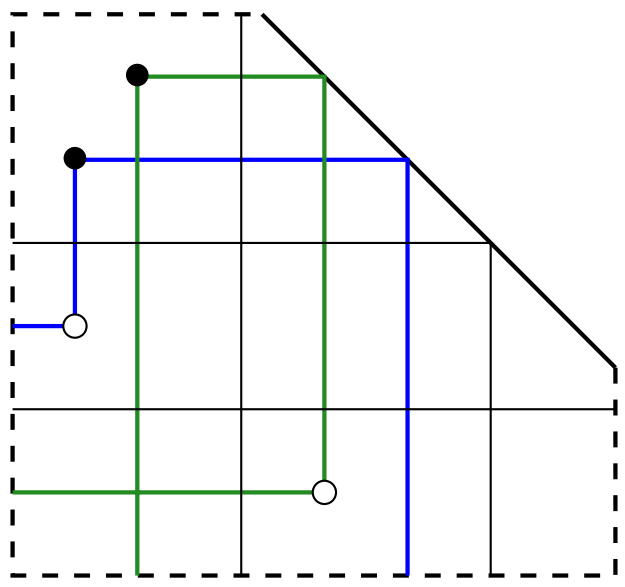
\[\begin{tikzpicture}[scale=.55]
  \draw[thick] (-.75,6.75) -- (3.5,2.5);
  \draw[thick,dashed] (3.5,2.5) -- (3.5,0) -- (-3.75,0) -- (-3.75,6.75) -- (-.75,6.75);
  \draw[thick, blue]
    (-3.75,3)--(-3,3)--(-3,5)-- (1,5) -- (1,0)  ;
  \draw[thick, ForestGreen]
   (-2.25,0) -- (-2.25,6) --(0,6) -- (0,1) -- (-3.75,1);
    \draw(-3.75,2)--(3.5,2);
    \draw(-1,6.75)--(-1,0);
    \draw(-3.75,4)--(2,4)--(2,0);
  \node at (0,1) {$\whitecirc$};
    \node at (-3,3) {$\whitecirc$};
  \node at (-2.25,6) {$\bullet$};
    \node at (-3,5) {$\bullet$};
\end{tikzpicture}\]
    Hence, $\demprodfpf(T)=\varnothing$.
\end{proof}

We now show that if a removable $-1$ in $A$ shares a neighbor with its partner, then $\Psi(A)$ is not symplectic.

\begin{lemma}
\label{lemma:sharedneighborbaddemprod}
    Let $A\in \asmsk(2n)$ and $(b,a)$ with $a<b$ so that $(b,a)$ is a removable $-1$ in $A$ that is not $\Sp$-removable,
    that is, $(b,a)$ and $(a,b)$ have common neighbors.
    Then $\demprodfpf(\Psi(A))=\varnothing$.
\end{lemma}

\begin{proof}
    First suppose that the left neighbor of $(a,b)$ is also a neighbor of $(b,a)$.  Then we have a configuration of the form 
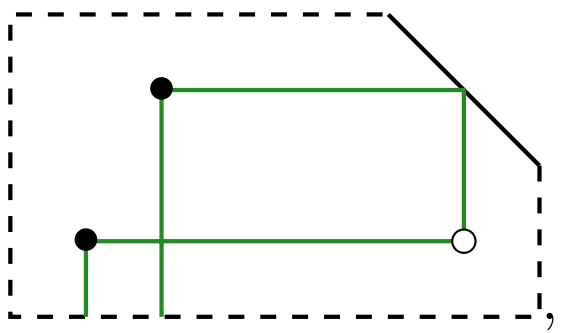
\[\begin{tikzpicture}[scale=1]
  \draw[thick] (.5,2.5) -- (1.5,1.5);
  \draw[thick,dashed](1.5,1.5)--(1.5,.5)--(-2,.5)--(-2,2.5)--(.5,2.5);
  \draw[thick, ForestGreen]
    (-1,.5) -- (-1,2) -- (1,2) -- (1,1) -- (-1.5,1)--(-1.5,.5);
  \node at (1,1) {$\whitecirc$};
  \node at (-1,2) {$\bullet$};
  \node at (-1.5,1) {$\bullet$};
\end{tikzpicture},\]
which, by definition of neighbor, has no $\seelbow$ or $\nwelbow$ within the green rectangle (though there may be other unpictured pipes and neighbors).
Thus $\Psi(A)$ contains a null cycle, so by \cref{lemma:basicbadconfigurations} $\demprodfpf(\Psi(A))=\varnothing$. 

Now suppose the left neighbor of $(a,b)$ is not a neighbor of $(b,a)$.
By assumption, $(a,b)$ and $(b,a)$ share some neighbor.
Choosing its partner if necessary, we have a common neighbor that sits below the diagonal.
This case leads to overlapping horseshoes, as pictured below.
\[\begin{tikzpicture}[scale=1]
  \draw[thick] (-.5,3.5) -- (1.5,1.5);
  \draw[thick,dashed](1.5,1.5)--(1.5,.5)--(-3,.5)--(-3,3.5)--(-.5,3.5);
  \draw[thick, ForestGreen]
    (-2.5,.5)--(-2.5,2) -- (1,2) -- (1,1) -- (-1.5,1)--(-1.5,.5);
  \draw[thick, blue]
    (-1,0.5) -- (-1,3) -- (0,3) -- (0,0.5);
  \node at (1,1) {$\whitecirc$};
  \node at (-1,3) {$\bullet$};
  \node at (-2.5,2) {$\bullet$};
  \node at (-1.5,1) {$\bullet$};
\end{tikzpicture}\]
Thus, by \cref{lemma:basicbadconfigurations}, $\demprodfpf(\Psi(A))=\varnothing$. 
\end{proof}

\begin{example}
    Let $A$ and $A'$ be as in \cref{ex:sharedneighborsafteronestep}.  Pictured below are $\Dsp_A$ and $\Dsp_{A'}$.
\[
\begin{tikzpicture}[scale=.45]
\draw (.5,11.5) -- (.5,-.5) -- (12.5,-.5) -- (.5,11.5);
\foreach \i in {1,...,11}{
  \draw (.5,{11.5-\i}) -- ({.5+\i},{11.5-\i}) -- ({.5+\i},-.5);
}
\draw[thick,ForestGreen] (1,-.5) -- (1,6) -- (6,6) -- (6,-.5);
\draw[thick,ForestGreen] (2,-.5) -- (2,5) -- (7,5) --(7,2)--(3,2) --(3,-.5);
\draw[thick,ForestGreen] (4,-.5) -- (4,4) -- (8,4) --(8,3)-- (5,3)--(5,-.5);
\draw[thick,ForestGreen] (7,-.5) -- (7,1) -- (11,1) -- (11,-.5);
\draw[thick,ForestGreen] (8,-.5) -- (8,0) -- (12,0) -- (12,-.5);
\draw[thick,ForestGreen] (9,-.5) -- (9,2) -- (10,2) -- (10,-.5);
\node at (1,6) {$\bullet$};
\node at (2,5) {$\bullet$};
\node at (3,2) {$\bullet$};
\node at (4,4) {$\bullet$};
\node at (5,3) {$\bullet$};
\node at (7,1) {$\bullet$};
\node at (8,0) {$\bullet$};
\node at (9,2) {$\bullet$};
\node at (7,2) {$\whitecirc$};
\node at (8,3) {$\whitecirc$};
\foreach \i in {1,...,12}{
  \node at (\i,-1) {\i};
}
\end{tikzpicture}
\hspace{4em}
\begin{tikzpicture}[scale=.45]
\draw (.5,11.5) -- (.5,-.5) -- (12.5,-.5) -- (.5,11.5);
\foreach \i in {1,...,11}{
  \draw (.5,{11.5-\i}) -- ({.5+\i},{11.5-\i}) -- ({.5+\i},-.5);
}
\draw[thick,ForestGreen] (1,-.5) -- (1,6) -- (6,6) -- (6,-.5);
\draw[very thick,blue] (2,-.5) -- (2,5) -- (7,5) --(7,2)--(3,2) --(3,-.5);
\draw[very thick,blue] (4,-.5) -- (4,7) -- (5,7) --(5,-.5);
\draw[thick,ForestGreen] (7,-.5) -- (7,1) -- (11,1) -- (11,-.5);
\draw[thick,ForestGreen] (8,-.5) -- (8,0) -- (12,0) -- (12,-.5);
\draw[thick,ForestGreen] (9,-.5) -- (9,2) -- (10,2) -- (10,-.5);
\node at (1,6) {$\bullet$};
\node at (2,5) {$\bullet$};
\node at (3,2) {$\bullet$};
\node at (4,7) {$\bullet$};
\node at (7,1) {$\bullet$};
\node at (8,0) {$\bullet$};
\node at (9,2) {$\bullet$};
\node at (7,2) {$\whitecirc$};
\foreach \i in {1,...,12}{
  \node at (\i,-1) {\i};
}
\end{tikzpicture}
\]
Notice that $\Dsp_A$ does not have any null cycles or overlapping horseshoes.
The overlapping horseshoes in $\Dsp_{A'}$ have been highlighted in blue.
Thus we know $\demprodfpf(\Dsp_{A'})=\varnothing$.
Because $A'$ is obtained from $A$ by symplectic inflation, \cref{prop:symplecticinflationpreservesdemprod} gives $\demprodfpf(\Dsp_A)=\varnothing$ as well.
\end{example}

\begin{proposition}
\label{prop:symplectickeyagreeswithdemprod}
Let $A\in \asmsk(2n)$.
Then $\keysp(A)=\demprodfpf(\Psi(A))$.
\end{proposition}
\begin{proof}
We may induct on the number of $-1$'s below the diagonal in $A$.
If there are none, $A$ is the permutation matrix of some $\tau \in \mf_{2n}$, and $\keysp(A)=\tau=\demprodfpf(\Psi(\tau))$.

Now assume $A$ has some $-1$ below the diagonal.
Let $(b,a)$ with $a<b$ be the lexicographically minimal such position, as in the algorithm defining $\keysp$.
This $-1$ is removable in $A$, since any other negative entry northwest of it, or that entry's partner, would be an earlier negative entry below the diagonal.
If it is not $\Sp$-removable (i.e., $(b,a)$ shares neighbors with $(a,b)$), then by definition $\keysp(A)=\varnothing$.
Applying \cref{lemma:sharedneighborbaddemprod} also gives $\demprodfpf(\Psi(A))=\varnothing$.

Now suppose $(b,a)$ is $\Sp$-removable.
Then we may apply symplectic inflation to get $A'$.
Thus,
\begin{align*}
    \keysp(A)&=\keysp(A')  &\text{(by definition of symplectic key)}\\
    &=\demprodfpf(\Psi(A'))&\text{(by the inductive hypothesis)}\\
    &=\demprodfpf(\Psi(A)) & \text{(by \cref{prop:symplecticinflationpreservesdemprod}),}
\end{align*}
as desired.
\end{proof}

As an immediate consequence, we have:
\begin{corollary}
    \label{c:symplectic-key}
    Let $A\in \asmsk(2n)$, and suppose we obtain $A'$ from $A$ by performing symplectic inflation at an $\Sp$-removable $-1$ in $A$.  Then $\keysp(A) = \keysp(A')$. 
\end{corollary}

Therefore, $\keysp$ is well defined regardless of the sequence of symplectic inflations applied.

We are now prepared to show symplectic bumpless pipe dreams and symplectic alternating sign matrices are in correspondence.

\begin{theorem}
    \label{t:symplectic-equality}
    For $\tau \in \mf_{2n}$, the restriction $\Psi:\asmsp(\tau) \to \bpdsp(\tau)$ is a bijection.
\end{theorem}

\begin{proof}
 By \cref{p:triangular-bpd}, $\Psi$ is a bijection from $\asmsk(2n)$ to triangular BPDs supported on $\lrtri_{2n}$.  By \cref{prop:symplectickeyagreeswithdemprod}, given $A\in\asmsk(2n)$, $\keysp(A)=\demprodfpf(\Psi(A))$.  In particular, $A\in \asmsp(\tau)$ if and only if $\Psi(A)\in\bpdsp(\tau)$.
\end{proof}

For $B$ a symplectic bumpless pipe dream, define
\[
\wt(B) = \left(\prod_{(j,i) \in B(\nowire)}  \beta (x_i \oplus x_j)\right) \left(  \prod_{(j,i) \in B(\nwelbow)} (1 + \beta (x_i \oplus x_j))\right).
\]
It is straightforward to see $\wt(\Psi(A)) = \wt_\Sp(A)$ for $A \in \asmsp(2n)$.
One immediate consequence is the main result of this paper.

\begin{corollary}
    \label{c:sbpd-formula}
    For $\tau \in \mf_{2n}$,
\[
\fkGsp_\tau = \beta^{-\ellfpf(\tau)} \sum_{B \in \bpdsp(\tau)} \wt(B).
\]
\end{corollary}

\cref{c:sbpd-formula} has a simple reinterpretation.
A \newword{marked symplectic bumpless pipe dream} is a pair $(B,M)$ where $B$ is a symplectic bumpless pipe dream and $M$ is a subset of its $\nwelbow$ tiles.
Let $\bpdsp_\mm(\tau)$ be the set of marked symplectic bumpless pipe dreams for $\tau$.
Define
\[
\wt_\mm(B,M) = \prod_{(j,i) \in B(\nowire) \cup M} \beta (x_i \oplus x_j).
\]
\begin{corollary}
    \label{c:marked}
    For $\tau \in \mf_{2n}$,
    \[
    \fkGsp_\tau =  \beta^{-\ellfpf(\tau)} \sum_{(B,M) \in \bpdsp_\mm(\tau)} \wt_\mm(B,M).
    \]
\end{corollary}

Next, we will use the $\beta = 0$ specialization of \cref{c:sbpd-formula} to derive a formula for the symplectic Schubert polynomial $\fkSsp_\tau$.
Towards that end, define
\[
\wt_\fkS(B) = \prod_{(j,i) \in B(\nowire)} (x_i + x_j).
\]

\begin{corollary}
    \label{c:schubert-bpd}
For $\tau \in \mf_{2n}$,
\[
\fkSsp_\tau = \sum_{B \in \bpdsp_\mr(\tau)} \wt_\fkS(B).
\]
\end{corollary}
\begin{proof}
    Given $B \in \bpdsp(\tau)$, \cref{l:reduced-blanks} shows that $|B(\nowire)| = \ellfpf(\tau)$ if and only if $B$ is reduced.
    Therefore the result follows by setting $\beta = 0$.
\end{proof}

\subsection{\texorpdfstring{Droops and $K$-theoretic droops}{Droops and K-theoretic droops}}
\label{subsection:droops}

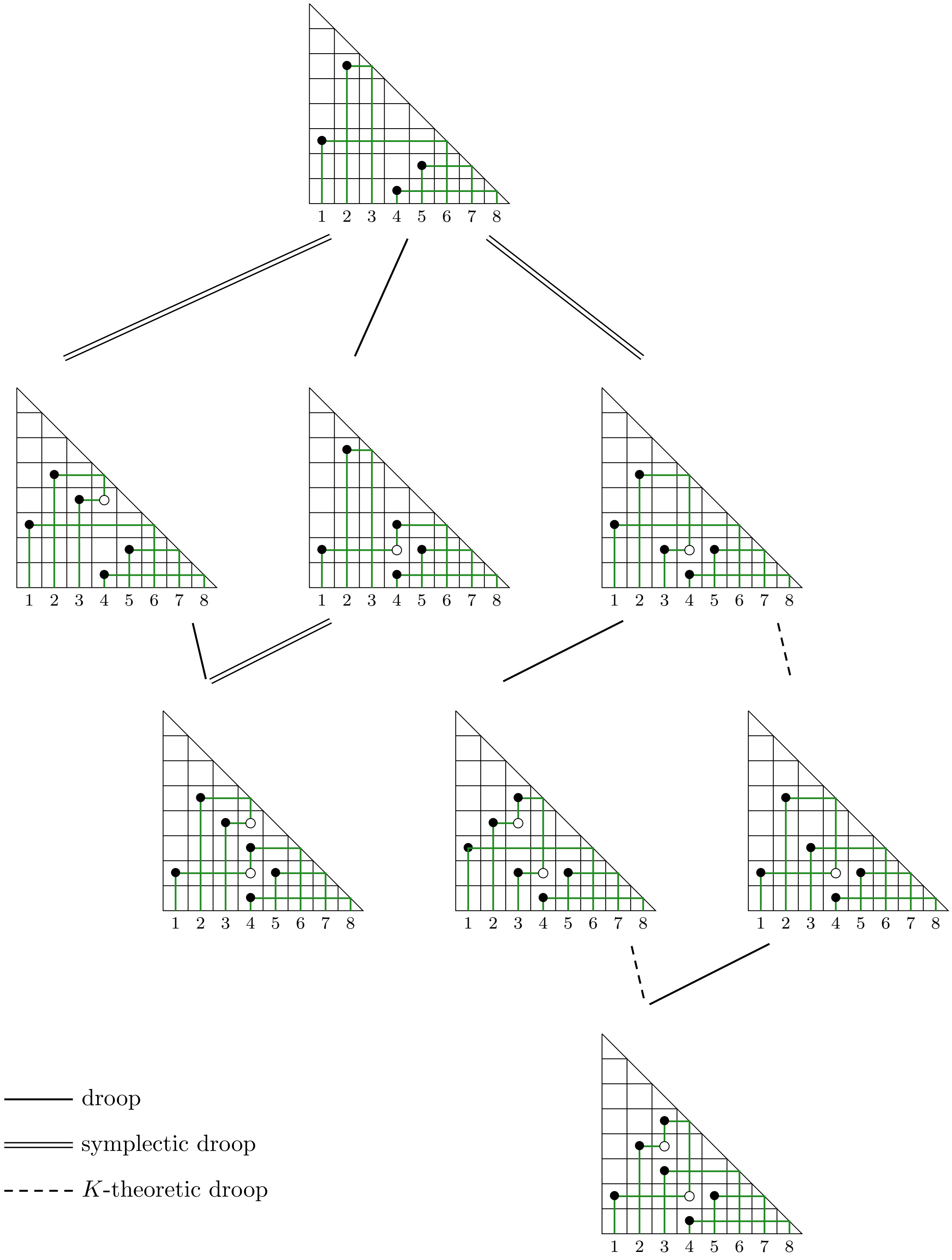
\begin{figure}
\centering
\begin{tikzpicture}[
    droop/.style={line width=.9pt},
    sdroop/.style={double,line width=.6pt,double distance=1.6pt},
    kdroop/.style={line width=.9pt,dash pattern=on 4pt off 3pt},
    skdroop/.style={double,line width=.6pt,double distance=1.6pt,dash pattern=on 4pt off 3pt},
    move/.style={shorten <=2pt, shorten >=2pt}]
\def\bpdgrid{%
  \draw (.5,7.5) -- (.5,-.5) -- (8.5,-.5) -- (.5,7.5);
  \foreach \rr in {1,...,7} {\draw (.5,7.5-\rr) -- (\rr+.5,7.5-\rr) -- (\rr+.5,-.5);}
  \foreach \cc in {1,...,8} {\node at (\cc,-1) {\scriptsize\cc};}}
\begin{scope}[shift={(7.2,1)},scale=.41]
  \bpdgrid
  \draw[thick,ForestGreen] (1,-.5) -- (1,2) -- (6,2) -- (6,-.5);
  \draw[thick,ForestGreen] (2,-.5) -- (2,5) -- (3,5) -- (3,-.5);
  \draw[thick,ForestGreen] (4,-.5) -- (4,0) -- (8,0) -- (8,-.5);
  \draw[thick,ForestGreen] (5,-.5) -- (5,1) -- (7,1) -- (7,-.5);
  \node at (1,2) {$\bullet$};
  \node at (2,5) {$\bullet$};
  \node at (4,0) {$\bullet$};
  \node at (5,1) {$\bullet$};
\end{scope}
\begin{scope}[shift={(2.4,-5.3)},scale=.41]
  \bpdgrid
  \draw[thick,ForestGreen] (1,-.5) -- (1,2) -- (6,2) -- (6,-.5);
  \draw[thick,ForestGreen] (2,-.5) -- (2,4) -- (4,4) -- (4,3) -- (3,3) -- (3,-.5);
  \draw[thick,ForestGreen] (4,-.5) -- (4,0) -- (8,0) -- (8,-.5);
  \draw[thick,ForestGreen] (5,-.5) -- (5,1) -- (7,1) -- (7,-.5);
  \node at (1,2) {$\bullet$};
  \node at (2,4) {$\bullet$};
  \node at (4,3) {$\whitecirc$};
  \node at (3,3) {$\bullet$};
  \node at (4,0) {$\bullet$};
  \node at (5,1) {$\bullet$};
\end{scope}
\begin{scope}[shift={(12,-5.3)},scale=.41]
  \bpdgrid
  \draw[thick,ForestGreen] (1,-.5) -- (1,2) -- (6,2) -- (6,-.5);
  \draw[thick,ForestGreen] (2,-.5) -- (2,4) -- (4,4) -- (4,1) -- (3,1) -- (3,-.5);
  \draw[thick,ForestGreen] (4,-.5) -- (4,0) -- (8,0) -- (8,-.5);
  \draw[thick,ForestGreen] (5,-.5) -- (5,1) -- (7,1) -- (7,-.5);
  \node at (1,2) {$\bullet$};
  \node at (2,4) {$\bullet$};
  \node at (4,1) {$\whitecirc$};
  \node at (3,1) {$\bullet$};
  \node at (4,0) {$\bullet$};
  \node at (5,1) {$\bullet$};
\end{scope}
\begin{scope}[shift={(7.2,-5.3)},scale=.41]
  \bpdgrid
  \draw[thick,ForestGreen] (1,-.5) -- (1,1) -- (4,1) -- (4,2) -- (6,2) -- (6,-.5);
  \draw[thick,ForestGreen] (2,-.5) -- (2,5) -- (3,5) -- (3,-.5);
  \draw[thick,ForestGreen] (4,-.5) -- (4,0) -- (8,0) -- (8,-.5);
  \draw[thick,ForestGreen] (5,-.5) -- (5,1) -- (7,1) -- (7,-.5);
  \node at (1,1) {$\bullet$};
  \node at (4,1) {$\whitecirc$};
  \node at (4,2) {$\bullet$};
  \node at (2,5) {$\bullet$};
  \node at (4,0) {$\bullet$};
  \node at (5,1) {$\bullet$};
\end{scope}
\begin{scope}[shift={(4.8,-10.6)},scale=.41]
  \bpdgrid
  \draw[thick,ForestGreen] (1,-.5) -- (1,1) -- (4,1) -- (4,2) -- (6,2) -- (6,-.5);
  \draw[thick,ForestGreen] (2,-.5) -- (2,4) -- (4,4) -- (4,3) -- (3,3) -- (3,-.5);
  \draw[thick,ForestGreen] (4,-.5) -- (4,0) -- (8,0) -- (8,-.5);
  \draw[thick,ForestGreen] (5,-.5) -- (5,1) -- (7,1) -- (7,-.5);
  \node at (1,1) {$\bullet$};
  \node at (4,1) {$\whitecirc$};
  \node at (4,2) {$\bullet$};
  \node at (2,4) {$\bullet$};
  \node at (4,3) {$\whitecirc$};
  \node at (3,3) {$\bullet$};
  \node at (4,0) {$\bullet$};
  \node at (5,1) {$\bullet$};
\end{scope}
\begin{scope}[shift={(9.6,-10.6)},scale=.41]
  \bpdgrid
  \node at (1,2) {$\bullet$};
  \draw[thick,ForestGreen] (1,-.5) -- (1,2) -- (6,2) -- (6,-.5);
  \draw[thick,ForestGreen] (2,-.5) -- (2,3) -- (3,3) -- (3,4) -- (4,4) -- (4,1) -- (3,1) -- (3,-.5);
  \draw[thick,ForestGreen] (4,-.5) -- (4,0) -- (8,0) -- (8,-.5);
  \draw[thick,ForestGreen] (5,-.5) -- (5,1) -- (7,1) -- (7,-.5);
  \node at (2,3) {$\bullet$};
  \node at (3,3) {$\whitecirc$};
  \node at (3,4) {$\bullet$};
  \node at (4,1) {$\whitecirc$};
  \node at (3,1) {$\bullet$};
  \node at (4,0) {$\bullet$};
  \node at (5,1) {$\bullet$};
\end{scope}
\begin{scope}[shift={(14.4,-10.6)},scale=.41]
  \bpdgrid
  \draw[thick,ForestGreen] (1,-.5) -- (1,1) -- (4,1) -- (4,2) -- (6,2) -- (6,-.5);
  \draw[thick,ForestGreen] (2,-.5) -- (2,4) -- (4,4) -- (4,2) -- (3,2) -- (3,-.5);
  \draw[thick,ForestGreen] (4,-.5) -- (4,0) -- (8,0) -- (8,-.5);
  \draw[thick,ForestGreen] (5,-.5) -- (5,1) -- (7,1) -- (7,-.5);
  \node at (1,1) {$\bullet$};
  \node at (4,1) {$\whitecirc$};
  \node at (2,4) {$\bullet$};
  \node at (4,0) {$\bullet$};
  \node at (5,1) {$\bullet$};
  \node at (3,2) {$\bullet$};
\end{scope}
\begin{scope}[shift={(12,-15.9)},scale=.41]
  \bpdgrid
  \draw[thick,ForestGreen] (1,-.5) -- (1,1) -- (4,1) -- (4,2) -- (6,2) -- (6,-.5);
  \draw[thick,ForestGreen] (2,-.5) -- (2,3) -- (3,3) -- (3,4) -- (4,4) -- (4,2) -- (3,2) -- (3,-.5);
  \draw[thick,ForestGreen] (4,-.5) -- (4,0) -- (8,0) -- (8,-.5);
  \draw[thick,ForestGreen] (5,-.5) -- (5,1) -- (7,1) -- (7,-.5);
  \node at (1,1) {$\bullet$};
  \node at (4,1) {$\whitecirc$};
  \node at (2,3) {$\bullet$};
  \node at (3,3) {$\whitecirc$};
  \node at (3,4) {$\bullet$};
  \node at (3,2) {$\bullet$};
  \node at (5,1) {$\bullet$};
  \node at (4,0) {$\bullet$};
\end{scope}
\draw[move,sdroop] (7.815,.283) -- (3.322,-1.774); 
\draw[move,droop] (9.045,.283) -- (8.122,-1.774); 
\draw[move,sdroop] (10.275,.283) -- (12.922,-1.774); 
\draw[move,droop] (5.475,-6.018) -- (5.723,-7.074); 
\draw[move,sdroop] (7.815,-6.018) -- (5.723,-7.074); 
\draw[move,droop] (12.615,-6.018) -- (10.522,-7.074); 
\draw[move,kdroop] (15.075,-6.018) -- (15.322,-7.074); 
\draw[move,kdroop] (12.675,-11.317) -- (12.922,-12.374); 
\draw[move,droop] (15.015,-11.317) -- (12.922,-12.374); 
\begin{scope}[shift={(2.4,-13.9)}]
  \draw[droop] (0,-0) -- ++(32pt,0) node[right] {droop};
  \draw[sdroop] (0,-.75) -- ++(32pt,0) node[right] {symplectic droop};
  \draw[kdroop] (0,-1.5) -- ++(32pt,0) node[right] {$K$-theoretic droop};
\end{scope}
\end{tikzpicture}
\caption{For $\tau = (1,6)(2,3)(4,8)(5,7)$, the eight elements of $\bpdsp(\tau)$, connected by the moves shown in the legend.
Each move goes downward.
No $K$-theoretic symplectic droop applies to these BPDs.
\label{f:bpd}}
\end{figure}

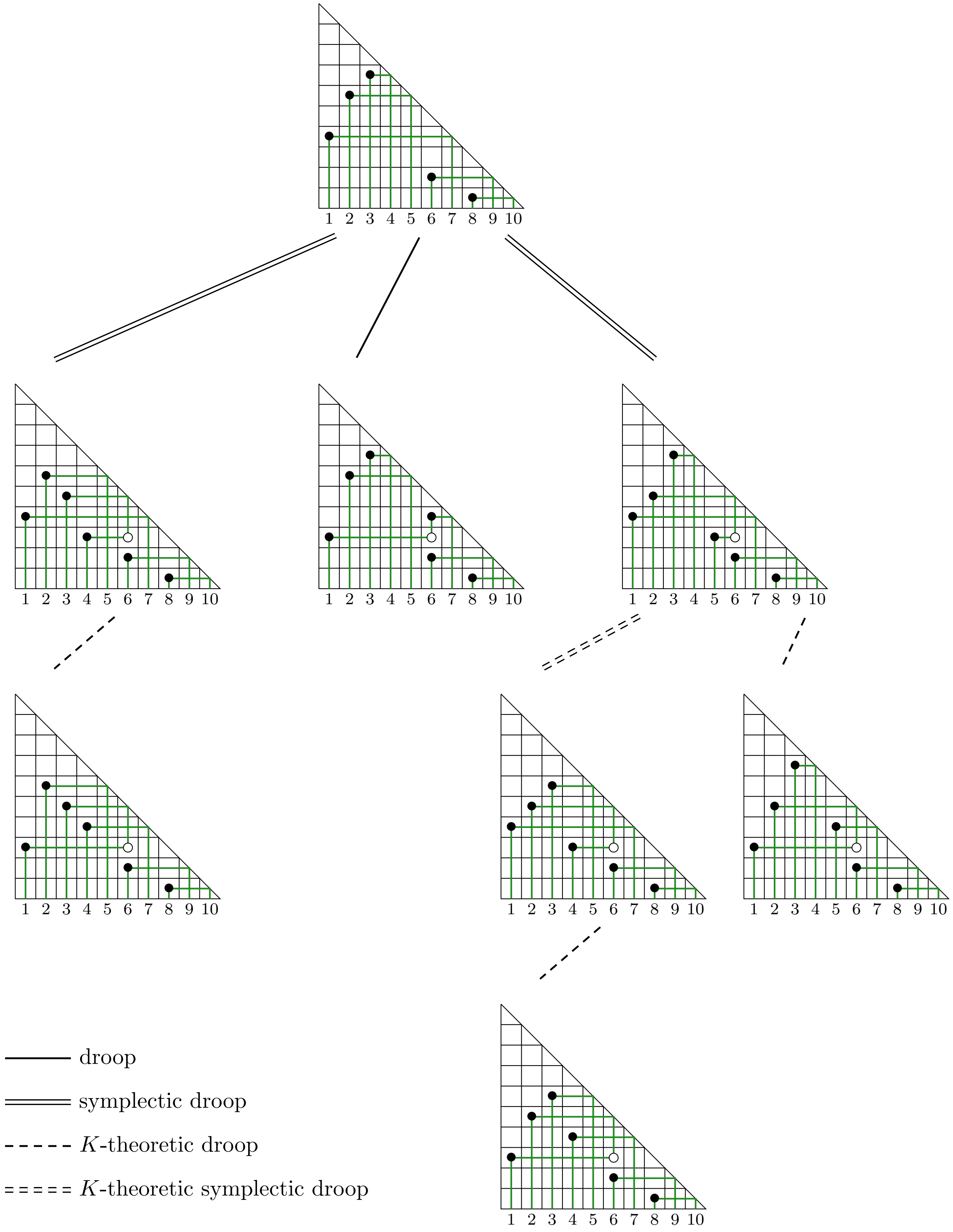
\begin{figure}
\centering
\begin{tikzpicture}[
    droop/.style={line width=.9pt},
    sdroop/.style={double,line width=.6pt,double distance=1.6pt},
    kdroop/.style={line width=.9pt,dash pattern=on 4pt off 3pt},
    skdroop/.style={double,line width=.6pt,double distance=1.6pt,dash pattern=on 4pt off 3pt},
    move/.style={shorten <=2pt, shorten >=2pt}]
\def\bpdgridten{%
  \draw (.5,9.5) -- (.5,-.5) -- (10.5,-.5) -- (.5,9.5);
  \foreach \rr in {1,...,9} {\draw (.5,9.5-\rr) -- (\rr+.5,9.5-\rr) -- (\rr+.5,-.5);}
  \foreach \cc in {1,...,10} {\node at (\cc,-1) {\scriptsize\cc};}}
\begin{scope}[shift={(5.188,1.2)},scale=.35]
  \bpdgridten
  \draw[thick,ForestGreen] (1,-.5) -- (1,3) -- (7,3) -- (7,-.5);
  \draw[thick,ForestGreen] (2,-.5) -- (2,5) -- (5,5) -- (5,-.5);
  \draw[thick,ForestGreen] (3,-.5) -- (3,6) -- (4,6) -- (4,-.5);
  \draw[thick,ForestGreen] (6,-.5) -- (6,1) -- (9,1) -- (9,-.5);
  \draw[thick,ForestGreen] (8,-.5) -- (8,0) -- (10,0) -- (10,-.5);
  \node at (3,6) {$\bullet$};
  \node at (2,5) {$\bullet$};
  \node at (1,3) {$\bullet$};
  \node at (6,1) {$\bullet$};
  \node at (8,0) {$\bullet$};
\end{scope}
\begin{scope}[shift={(0,-5.3)},scale=.35]
  \bpdgridten
  \draw[thick,ForestGreen] (1,-.5) -- (1,3) -- (7,3) -- (7,-.5);
  \draw[thick,ForestGreen] (2,-.5) -- (2,5) -- (5,5) -- (5,-.5);
  \draw[thick,ForestGreen] (3,-.5) -- (3,4) -- (6,4) -- (6,2) -- (4,2) -- (4,-.5);
  \draw[thick,ForestGreen] (6,-.5) -- (6,1) -- (9,1) -- (9,-.5);
  \draw[thick,ForestGreen] (8,-.5) -- (8,0) -- (10,0) -- (10,-.5);
  \node at (2,5) {$\bullet$};
  \node at (3,4) {$\bullet$};
  \node at (1,3) {$\bullet$};
  \node at (4,2) {$\bullet$};
  \node at (6,2) {$\whitecirc$};
  \node at (6,1) {$\bullet$};
  \node at (8,0) {$\bullet$};
\end{scope}
\begin{scope}[shift={(10.375,-5.3)},scale=.35]
  \bpdgridten
  \draw[thick,ForestGreen] (1,-.5) -- (1,3) -- (7,3) -- (7,-.5);
  \draw[thick,ForestGreen] (2,-.5) -- (2,4) -- (6,4) -- (6,2) -- (5,2) -- (5,-.5);
  \draw[thick,ForestGreen] (3,-.5) -- (3,6) -- (4,6) -- (4,-.5);
  \draw[thick,ForestGreen] (6,-.5) -- (6,1) -- (9,1) -- (9,-.5);
  \draw[thick,ForestGreen] (8,-.5) -- (8,0) -- (10,0) -- (10,-.5);
  \node at (3,6) {$\bullet$};
  \node at (2,4) {$\bullet$};
  \node at (1,3) {$\bullet$};
  \node at (5,2) {$\bullet$};
  \node at (6,2) {$\whitecirc$};
  \node at (6,1) {$\bullet$};
  \node at (8,0) {$\bullet$};
\end{scope}
\begin{scope}[shift={(5.188,-5.3)},scale=.35]
  \bpdgridten
  \draw[thick,ForestGreen] (1,-.5) -- (1,2) -- (6,2) -- (6,3) -- (7,3) -- (7,-.5);
  \draw[thick,ForestGreen] (2,-.5) -- (2,5) -- (5,5) -- (5,-.5);
  \draw[thick,ForestGreen] (3,-.5) -- (3,6) -- (4,6) -- (4,-.5);
  \draw[thick,ForestGreen] (6,-.5) -- (6,1) -- (9,1) -- (9,-.5);
  \draw[thick,ForestGreen] (8,-.5) -- (8,0) -- (10,0) -- (10,-.5);
  \node at (3,6) {$\bullet$};
  \node at (2,5) {$\bullet$};
  \node at (6,3) {$\bullet$};
  \node at (1,2) {$\bullet$};
  \node at (6,2) {$\whitecirc$};
  \node at (6,1) {$\bullet$};
  \node at (8,0) {$\bullet$};
\end{scope}
\begin{scope}[shift={(0,-10.6)},scale=.35]
  \bpdgridten
  \draw[thick,ForestGreen] (1,-.5) -- (1,2) -- (6,2) -- (6,4) -- (3,4) -- (3,-.5);
  \draw[thick,ForestGreen] (2,-.5) -- (2,5) -- (5,5) -- (5,-.5);
  \draw[thick,ForestGreen] (4,-.5) -- (4,3) -- (7,3) -- (7,-.5);
  \draw[thick,ForestGreen] (6,-.5) -- (6,1) -- (9,1) -- (9,-.5);
  \draw[thick,ForestGreen] (8,-.5) -- (8,0) -- (10,0) -- (10,-.5);
  \node at (2,5) {$\bullet$};
  \node at (3,4) {$\bullet$};
  \node at (4,3) {$\bullet$};
  \node at (1,2) {$\bullet$};
  \node at (6,2) {$\whitecirc$};
  \node at (6,1) {$\bullet$};
  \node at (8,0) {$\bullet$};
\end{scope}
\begin{scope}[shift={(8.3,-10.6)},scale=.35]
  \bpdgridten
  \draw[thick,ForestGreen] (1,-.5) -- (1,3) -- (7,3) -- (7,-.5);
  \draw[thick,ForestGreen] (2,-.5) -- (2,4) -- (6,4) -- (6,2) -- (4,2) -- (4,-.5);
  \draw[thick,ForestGreen] (3,-.5) -- (3,5) -- (5,5) -- (5,-.5);
  \draw[thick,ForestGreen] (6,-.5) -- (6,1) -- (9,1) -- (9,-.5);
  \draw[thick,ForestGreen] (8,-.5) -- (8,0) -- (10,0) -- (10,-.5);
  \node at (3,5) {$\bullet$};
  \node at (2,4) {$\bullet$};
  \node at (1,3) {$\bullet$};
  \node at (4,2) {$\bullet$};
  \node at (6,2) {$\whitecirc$};
  \node at (6,1) {$\bullet$};
  \node at (8,0) {$\bullet$};
\end{scope}
\begin{scope}[shift={(12.45,-10.6)},scale=.35]
  \bpdgridten
  \draw[thick,ForestGreen] (1,-.5) -- (1,2) -- (6,2) -- (6,4) -- (2,4) -- (2,-.5);
  \draw[thick,ForestGreen] (3,-.5) -- (3,6) -- (4,6) -- (4,-.5);
  \draw[thick,ForestGreen] (5,-.5) -- (5,3) -- (7,3) -- (7,-.5);
  \draw[thick,ForestGreen] (6,-.5) -- (6,1) -- (9,1) -- (9,-.5);
  \draw[thick,ForestGreen] (8,-.5) -- (8,0) -- (10,0) -- (10,-.5);
  \node at (3,6) {$\bullet$};
  \node at (2,4) {$\bullet$};
  \node at (5,3) {$\bullet$};
  \node at (1,2) {$\bullet$};
  \node at (6,2) {$\whitecirc$};
  \node at (6,1) {$\bullet$};
  \node at (8,0) {$\bullet$};
\end{scope}
\begin{scope}[shift={(8.3,-15.9)},scale=.35]
  \bpdgridten
  \draw[thick,ForestGreen] (1,-.5) -- (1,2) -- (6,2) -- (6,4) -- (2,4) -- (2,-.5);
  \draw[thick,ForestGreen] (3,-.5) -- (3,5) -- (5,5) -- (5,-.5);
  \draw[thick,ForestGreen] (4,-.5) -- (4,3) -- (7,3) -- (7,-.5);
  \draw[thick,ForestGreen] (6,-.5) -- (6,1) -- (9,1) -- (9,-.5);
  \draw[thick,ForestGreen] (8,-.5) -- (8,0) -- (10,0) -- (10,-.5);
  \node at (3,5) {$\bullet$};
  \node at (2,4) {$\bullet$};
  \node at (4,3) {$\bullet$};
  \node at (1,2) {$\bullet$};
  \node at (6,2) {$\whitecirc$};
  \node at (6,1) {$\bullet$};
  \node at (8,0) {$\bullet$};
\end{scope}
\draw[move,sdroop] (5.713,.588) -- (.787,-1.59); 
\draw[move,droop] (7.112,.588) -- (5.975,-1.59); 
\draw[move,sdroop] (8.512,.588) -- (11.162,-1.59); 
\draw[move,kdroop] (1.925,-5.912) -- (.787,-6.89); 
\draw[move,skdroop] (10.9,-5.912) -- (9.088,-6.89); 
\draw[move,kdroop] (13.7,-5.912) -- (13.238,-6.89); 
\draw[move,kdroop] (10.225,-11.213) -- (9.088,-12.19); 
\begin{scope}[shift={(0,-13.5)}]
  \draw[droop] (0,-0) -- ++(32pt,0) node[right] {droop};
  \draw[sdroop] (0,-.75) -- ++(32pt,0) node[right] {symplectic droop};
  \draw[kdroop] (0,-1.5) -- ++(32pt,0) node[right] {$K$-theoretic droop};
  \draw[skdroop] (0,-2.25) -- ++(32pt,0) node[right] {$K$-theoretic symplectic droop};
\end{scope}
\end{tikzpicture}
\caption{For $\tau = (1,7)(2,5)(3,4)(6,9)(8,10)$, the eight elements of $\bpdsp(\tau)$, connected by the moves shown in the legend; all four types of move occur.
\label{f:bpd-all-moves}}
\end{figure}

Next, we define operations on triangular bumpless pipe dreams.
Our ultimate goal will be to show that $\bpdsp(\tau)$ is connected by moves called \emph{droop moves}, \emph{symplectic droop moves}, and their $K$-theoretic analogues.
We prove this in \cref{c:droop-graph}; see \cref{f:bpd,f:bpd-all-moves} for examples.

We start by defining droop moves.
Suppose $i<k<j<\ell$, $B(j,i)=\seelbow$, and $B(\ell,k)=\nowire$ with no other elbows in the Type I region $[(j,i),(\ell,k)]_\preceq$.
The \newword{droop move} on this interval replaces $B(j,i)$ by $\nowire$ and $B(\ell,k)$ by $\nwelbow$, rerouting the pipe through $(j,i)$ by setting $B(j,k) = B(\ell,i) = \seelbow$.
See the image below:
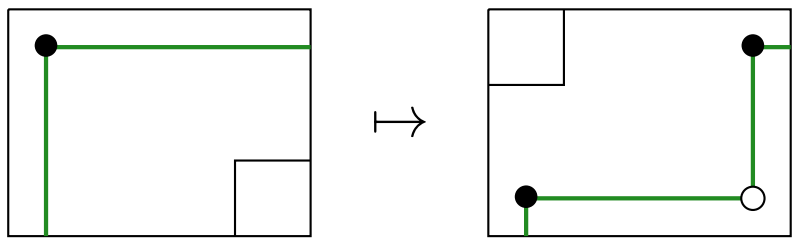
\[
\begin{tikzpicture}[scale=.5]
	\draw (0,0) -- (4,0) -- (4,-3) -- (0,-3) -- (0,0);
	\draw (3,-3) -- (3,-2) -- (4,-2);
	\draw[ForestGreen,thick] (.5,-3) -- (.5,-.5) -- (4,-.5);
	\node at (.5,-.5) {$\bullet$};
\end{tikzpicture}
\quad
\raisebox{1.75em}{$\mapsto$}
\quad
\begin{tikzpicture}[scale=.5]
	\draw (0,0) -- (4,0) -- (4,-3) -- (0,-3) -- (0,0);
	\draw (0,-1) -- (1,-1) -- (1,0);
	\draw[ForestGreen,thick] (.5,-3) -- (.5,-2.5) -- (3.5,-2.5)  -- (3.5,-.5) -- (4,-.5);
	\node at (.5,-2.5) {$\bullet$};
	\node at (3.5,-.5) {$\bullet$};
	\node at (3.5,-2.5) {$\whitecirc$};
\end{tikzpicture}
\]
A \newword{$K$-theoretic droop} starts instead from a $\nwelbow$ that is already present.
Suppose $i<i'<k<j<\ell$, $B(j,i)=B(\ell,i')=\seelbow$, and $B(\ell,k)=\nwelbow$, with no other elbows in the Type I region.
The $K$-theoretic droop replaces $B(j,i)$ by $\nowire$, $B(\ell,i')$ by $\cross$, and places $\seelbow$ tiles in positions $(j,i')$ and $(\ell,i)$, leaving $B(\ell,k) = \nwelbow$ unchanged.
See the image below:
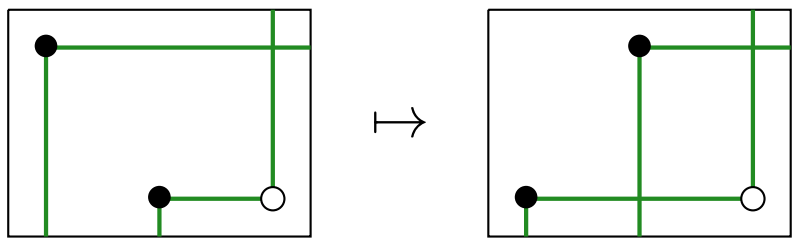
\[\begin{tikzpicture}[scale=.5]
	\draw (0,0) -- (4,0) -- (4,-3) -- (0,-3) -- (0,0);
	\draw[ForestGreen,thick] (.5,-3) -- (.5,-.5) -- (4,-.5);
    \draw[ForestGreen,thick] (2,-3) -- (2,-2.5) -- (3.5,-2.5)--(3.5,0);
	\node at (.5,-.5) {$\bullet$};
    \node at (2,-2.5) {$\bullet$};
    \node at (3.5,-2.5) {$\whitecirc$};
\end{tikzpicture}
\quad
\raisebox{1.75em}{$\mapsto$}
\quad
\begin{tikzpicture}[scale=.5]
	\draw (0,0) -- (4,0) -- (4,-3) -- (0,-3) -- (0,0);
	\draw[ForestGreen,thick] (2,-3) -- (2,-.5) -- (4,-.5);
    \draw[ForestGreen,thick] (.5,-3) -- (.5,-2.5) -- (3.5,-2.5)--(3.5,0);
	\node at (2,-.5) {$\bullet$};
    \node at (.5,-2.5) {$\bullet$};
    \node at (3.5,-2.5) {$\whitecirc$};
\end{tikzpicture}
\]
Again, there may be other (unpictured) pipes within this Type I region, but they do not bend there.
One can define an additional $K$-theoretic droop as the transpose of ours, but this move is not needed for our arguments.

\begin{lemma}
\label{lemma:closeddroopsandkdroops}
Let $B$ be a BPD or triangular BPD and $B\mapsto B'$ be a droop move or a $K$-theoretic droop move.
    \begin{enumerate}
        \item If $B\in \bpd(w)$ then $B'\in \bpd(w)$.
        \item If $B\in \bpdsp(\tau)$ then $B'\in \bpdsp(\tau)$.
    \end{enumerate}
\end{lemma}
\begin{proof}
 Part (1) is known for ordinary BPDs (see \cite{weigandt2021bumpless}*{Prop.~4.5}).
Part (2) follows from \cref{lemma:basiclocalreplacements}: the region of a droop or $K$-theoretic droop lies strictly below the diagonal, and the proof of Part (1) decomposes these moves into the replacements of that lemma. \end{proof}

A symplectic droop is the Type II analogue of the droop move.
Suppose $i<j<k<\ell$, $B(j,i)=\seelbow$, and $B(\ell,k)=\nowire$, with no other elbows in the Type II region, which is the interval $[(j,i),(k,\ell)]_\preceq$ folded into the lower triangle.
The \newword{symplectic droop} replaces $B(j,i)$ by $\nowire$, $B(\ell,k)$ by $\nwelbow$, $B(j,j)$ by $\notri$ and $B(k,k)$ by $\yestri$, and places $\seelbow$ tiles in positions $(k,i)$ and $(\ell,j)$.
It is a replacement of the form pictured below:
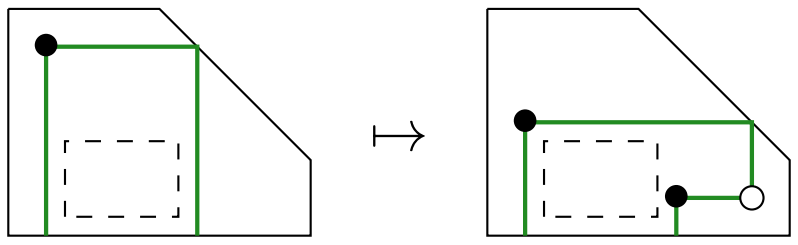
\[\begin{tikzpicture}[scale=.5]
	\draw (0,0) -- (2,0) -- (4,-2)--(4,-3)-- (0,-3) -- (0,0);
	\draw[ForestGreen,thick] (.5,-3) -- (.5,-.5) -- (2.5,-.5) -- (2.5,-3);
    \draw[dashed] (.75,-2.75) -- (.75,-1.75) -- (2.25,-1.75) -- (2.25,-2.75) -- (.75,-2.75);
	\node at (.5,-.5) {$\bullet$};
\end{tikzpicture}
\quad
\raisebox{1.5em}{$\mapsto$}
\quad
\begin{tikzpicture}[scale=.5]
\draw (0,0) -- (2,0) -- (4,-2)--(4,-3)-- (0,-3) -- (0,0);
	\draw[ForestGreen,thick] (.5,-3) -- (.5,-1.5) -- (3.5,-1.5) -- (3.5,-2.5) -- (2.5,-2.5) -- (2.5,-3);
    \draw[dashed] (.75,-2.75) -- (.75,-1.75) -- (2.25,-1.75) -- (2.25,-2.75) -- (.75,-2.75);
	\node at (.5,-1.5) {$\bullet$};
	\node at (2.5,-2.5) {$\bullet$};
	\node at (3.5,-2.5) {$\whitecirc$};
\end{tikzpicture}\]
The condition that there are no other elbows in the Type II region does not apply to the dashed inset.

Now suppose $i<j<j'<k<\ell$, $B(j,i)=B(\ell,j')=\seelbow$ and $B(\ell,k)=\nwelbow$ with no other elbows in the Type II region.
A \newword{$K$-theoretic symplectic droop} replaces $B(j,i)$ by $\nowire$, $B(\ell,j')$ by $\cross$, $B(j,j)$ by $\notri$, $B(j',j')$ by $\yestri$ and places $\seelbow$ tiles at $(j',i)$ and $(\ell,j)$.
The $\nwelbow$ in position $(\ell,k)$ is unchanged.
It is a replacement of the following form:
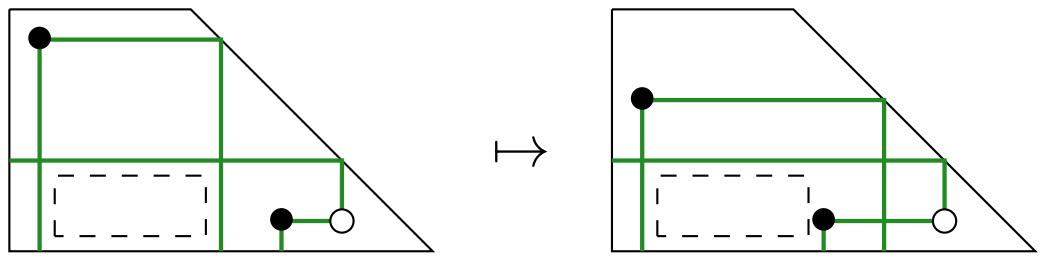
\[\begin{tikzpicture}[scale=.4]
	\draw (0,0) -- (3,0) -- (7,-4) -- (0,-4) -- (0,0);
	\draw[ForestGreen,thick] (.5,-4) -- (.5,-.5) -- (3.5,-.5) -- (3.5,-4);
    \draw[ForestGreen,thick] (4.5,-4) -- (4.5,-3.5) -- (5.5,-3.5) -- (5.5,-2.5) -- (0,-2.5);
    \draw[dashed] (.75,-3.75)--(.75,-2.75)--(3.25,-2.75)--(3.25,-3.75)--(.75,-3.75);
	\node at (.5,-.5) {$\bullet$};
    \node at (4.5,-3.5) {$\bullet$};
    \node at (5.5,-3.5) {$\whitecirc$};
\end{tikzpicture}
\quad
\raisebox{1.5em}{$\mapsto$}
\quad
\begin{tikzpicture}[scale=.4]
	\draw (0,0) -- (3,0) -- (7,-4) -- (0,-4) -- (0,0);
	\draw[ForestGreen,thick] (.5,-4) -- (.5,-1.5) -- (4.5,-1.5) -- (4.5,-4);
    \draw[ForestGreen,thick]   (3.5,-4) -- (3.5,-3.5) --  (5,-3.5) -- (5.5,-3.5) -- (5.5,-2.5) -- (0,-2.5);
    \draw[dashed] (.75,-3.75)--(.75,-2.75)--(3.25,-2.75)--(3.25,-3.75)--(.75,-3.75);
	\node at (.5,-1.5) {$\bullet$};
    \node at (5.5,-3.5) {$\whitecirc$};
    \node at (3.5,-3.5) {$\bullet$};
\end{tikzpicture}
\]
Again, the dashed inset can have additional elbows.
If $B \mapsto B'$ is any of the four moves above, we call the reverse replacement $B' \mapsto B$ the corresponding \newword{inverse} move.

There are similar moves that could reasonably be referred to as ($K$-theoretic) droops or symplectic droops.
However, as we show in \cref{l:deflation-droop}, these four moves suffice to realize symplectic deflation and connect $\bpdsp(\tau)$, so we restrict our focus to them alone.

\begin{lemma}
\label{lemma:closedsymplecticdroopsandkdroops}
Suppose $B$ is a triangular BPD.
Let $B\mapsto B'$ be a symplectic droop move or a $K$-theoretic symplectic droop move.
 If $B\in \bpdsp(\tau)$, then $B'\in \bpdsp(\tau)$.
\end{lemma}
\begin{proof}
By interpreting $B$ as a triangular planar history, we illustrate how to realize a symplectic droop.
First do a skew droop so that the newly created $\nwelbow$ sits above the dashed inset, and then use a sequence of the moves from \cref{lemma:basiclocalreplacements} to bring the green pipe into the desired location:
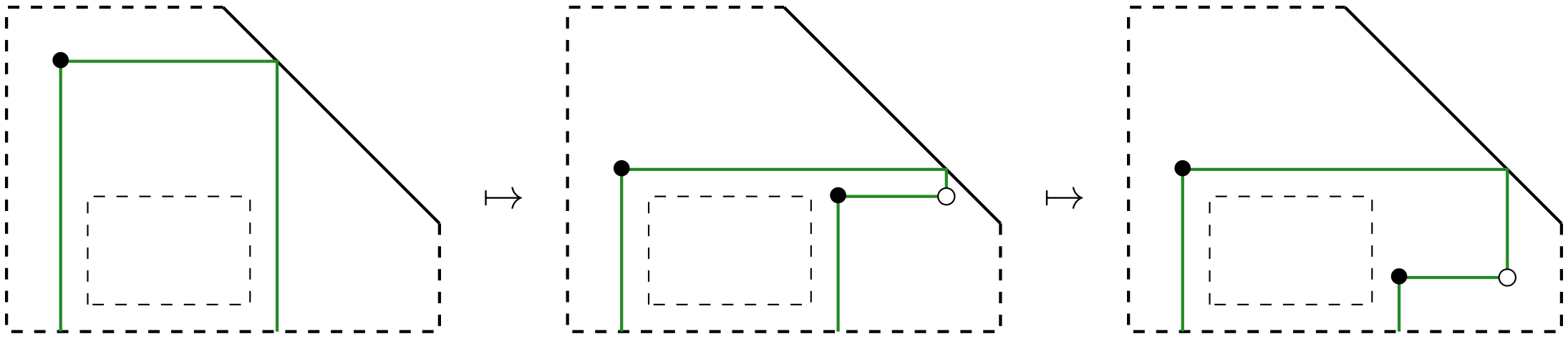
\[\begin{tikzpicture}[scale=1]
  \draw[thick] (2,0) -- (4,-2);
  \draw[thick,dashed] (4,-2) -- (4,-3) -- (0,-3) -- (0,0) -- (2,0);
  \draw[ForestGreen,thick] (.5,-3) -- (.5,-.5) -- (2.5,-.5) -- (2.5,-3);
  \draw[dashed] (.75,-2.75) -- (.75,-1.75) -- (2.25,-1.75) -- (2.25,-2.75) -- (.75,-2.75);
  \node at (.5,-.5) {$\bullet$};
\end{tikzpicture}
\quad
\raisebox{3em}{$\mapsto$}
\quad
\begin{tikzpicture}[scale=1]
  \draw[thick] (2,0) -- (4,-2);
  \draw[thick,dashed] (4,-2) -- (4,-3) -- (0,-3) -- (0,0) -- (2,0);
  \draw[ForestGreen,thick] (.5,-3) -- (.5,-1.5) -- (3.5,-1.5) -- (3.5,-1.75) -- (2.5,-1.75) -- (2.5,-3);
  \draw[dashed] (.75,-2.75) -- (.75,-1.75) -- (2.25,-1.75) -- (2.25,-2.75) -- (.75,-2.75);
  \node at (.5,-1.5) {$\bullet$};
  \node at (2.5,-1.75) {$\bullet$};
  \node at (3.5,-1.75) {$\whitecirc$};
\end{tikzpicture}
\quad
\raisebox{3em}{$\mapsto$}
\quad
\begin{tikzpicture}[scale=1]
  \draw[thick] (2,0) -- (4,-2);
  \draw[thick,dashed] (4,-2) -- (4,-3) -- (0,-3) -- (0,0) -- (2,0);
  \draw[ForestGreen,thick] (.5,-3) -- (.5,-1.5) -- (3.5,-1.5) -- (3.5,-2.5) -- (2.5,-2.5) -- (2.5,-3);
  \draw[dashed] (.75,-2.75) -- (.75,-1.75) -- (2.25,-1.75) -- (2.25,-2.75) -- (.75,-2.75);
  \node at (.5,-1.5) {$\bullet$};
  \node at (2.5,-2.5) {$\bullet$};
  \node at (3.5,-2.5) {$\whitecirc$};
\end{tikzpicture}\]
All unpictured elbows in this region lie in the dashed inset, which these moves avoid.

We now analyze the inverse $K$-theoretic symplectic droop.
In the sequence below, first perform a skew droop on the blue pipe, then the last replacement from \cref{lemma:basiclocalreplacements}, then apply an inverse droop to the blue pipe, and finally perform an inverse symplectic droop:
       
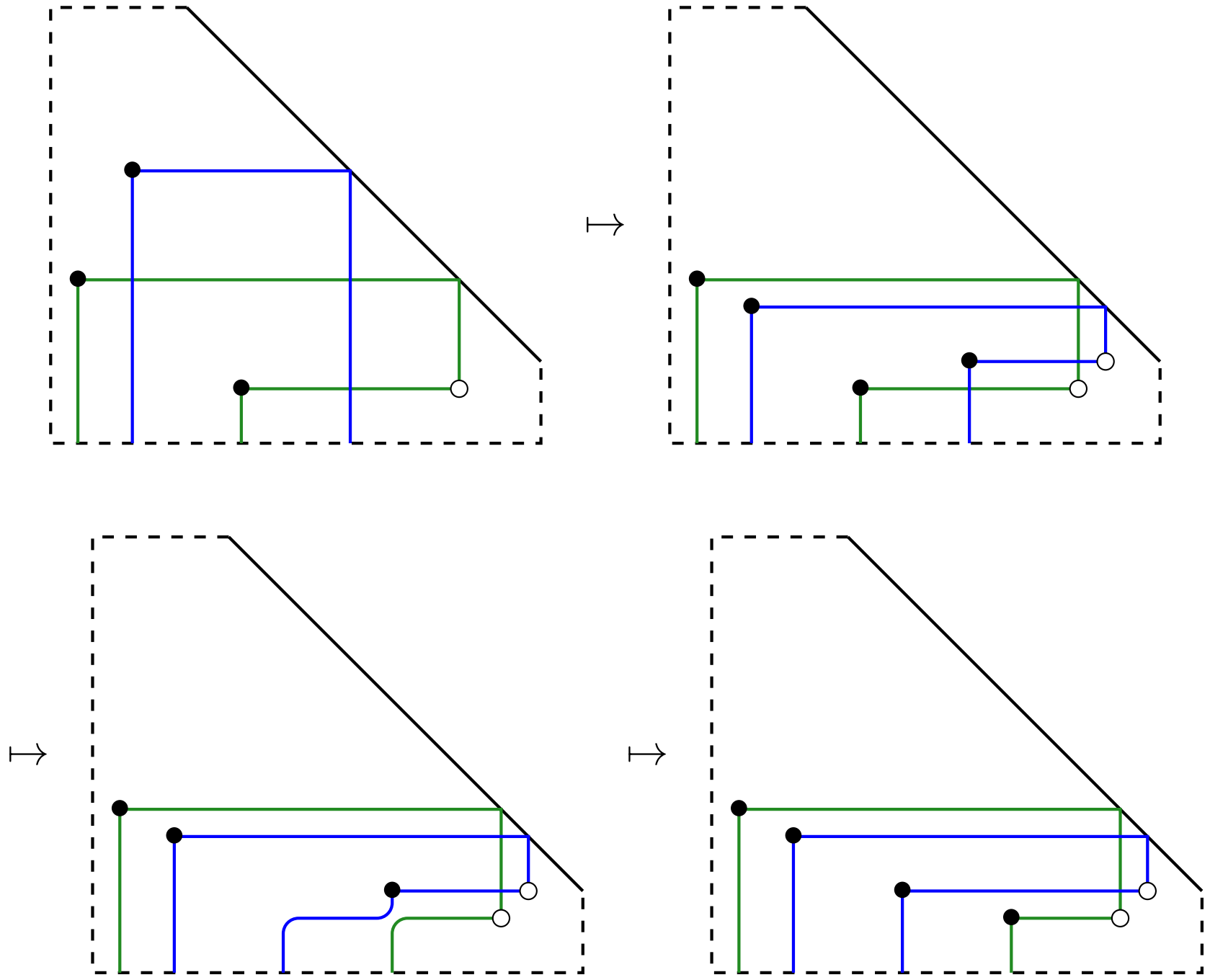
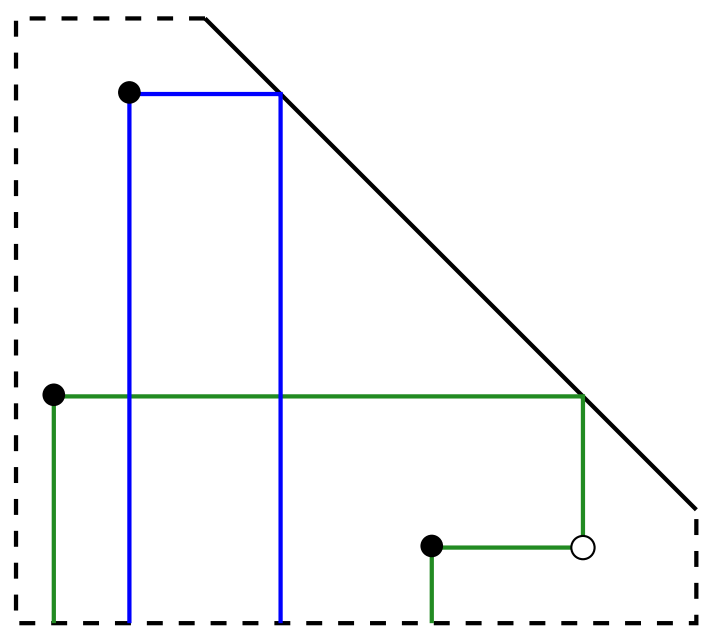
\[\begin{tikzpicture}[scale=1]
 \draw[thick] (-1.5,4.5) -- (1.75,1.25);
  \draw[thick,dashed] (-1.5,4.5) --(-2.75,4.5)--(-2.75,.5)--(1.75,.5)--(1.75,1.25);
  \draw[thick, ForestGreen]
    (-2.5,.5)--(-2.5,2) -- (1,2) -- (1,1) -- (-1,1)--(-1,.5);
  \draw[thick, blue]
    (-2,0.5) -- (-2,3) -- (0,3) -- (0,0.5);
  \node at (1,1) {$\whitecirc$};
  \node at (-2,3) {$\bullet$};
  \node at (-2.5,2) {$\bullet$};
  \node at (-1,1) {$\bullet$};
\end{tikzpicture}
\quad
\raisebox{5em}{$\mapsto$}
\quad
\begin{tikzpicture}[scale=1]
\draw[thick] (-1.5,4.5) -- (1.75,1.25);
  \draw[thick,dashed] (-1.5,4.5) --(-2.75,4.5)--(-2.75,.5)--(1.75,.5)--(1.75,1.25);
  \draw[thick, ForestGreen]
    (-2.5,.5)--(-2.5,2) -- (1,2) -- (1,1) -- (-1,1)--(-1,.5);
  \draw[thick, blue]
    (-2,0.5) -- (-2,1.75) -- (0,1.75) --(1.25,1.75)--(1.25,1.25)--(0,1.25)[rounded corners]-- (0,0.5);
  \node at (1,1) {$\whitecirc$};
  \node at (-2,1.75) {$\bullet$};
  \node at (1.25,1.25) {$\whitecirc$};
  \node at (0,1.25){$\bullet$};
  \node at (-2.5,2) {$\bullet$};
  \node at (-1,1) {$\bullet$};
  \end{tikzpicture}
\]

   \[\raisebox{5em}{$\mapsto$}
\quad
\begin{tikzpicture}[scale=1]
  \draw[thick] (-1.5,4.5) -- (1.75,1.25);
  \draw[thick,dashed] (-1.5,4.5) --(-2.75,4.5)--(-2.75,.5)--(1.75,.5)--(1.75,1.25);
  \draw[thick, ForestGreen]
    (-2.5,.5)--(-2.5,2) -- (1,2) -- (1,1)[rounded corners] --(0,1)-- (0,0.5) ;
  \draw[thick, blue]
    (-2,0.5) -- (-2,1.75) -- (0,1.75) --(1.25,1.75)--(1.25,1.25)--(0,1.25)[rounded corners]--(0,1)--(-1,1)[sharp corners]--(-1,.5);
  \node at (1,1) {$\whitecirc$};
  \node at (-2,1.75) {$\bullet$};
  \node at (1.25,1.25) {$\whitecirc$};
  \node at (0,1.25){$\bullet$};
  \node at (-2.5,2) {$\bullet$};
\end{tikzpicture}
\quad
\raisebox{5em}{$\mapsto$}
\quad
\begin{tikzpicture}[scale=1]
  \draw[thick] (-1.5,4.5) -- (1.75,1.25);
  \draw[thick,dashed] (-1.5,4.5) --(-2.75,4.5)--(-2.75,.5)--(1.75,.5)--(1.75,1.25);
  \draw[thick, ForestGreen]
    (-2.5,.5)--(-2.5,2) -- (1,2) -- (1,1)--(0,1)-- (0,0.5) ;
  \draw[thick, blue]
    (-2,0.5) -- (-2,1.75) -- (0,1.75) --(1.25,1.75)--(1.25,1.25)--(0,1.25)--(-1,1.25)--(-1,.5);
  \node at (1,1) {$\whitecirc$};
  \node at (-2,1.75) {$\bullet$};
  \node at (1.25,1.25) {$\whitecirc$};
  \node at (0,1){$\bullet$};
  \node at (-1,1.25){$\bullet$};
  \node at (-2.5,2) {$\bullet$};
\end{tikzpicture}
\]

\[
\raisebox{5em}{$\mapsto$}
\quad
\begin{tikzpicture}[scale=1]
  \draw[thick] (-1.5,4.5) -- (1.75,1.25);
  \draw[thick,dashed] (-1.5,4.5) --(-2.75,4.5)--(-2.75,.5)--(1.75,.5)--(1.75,1.25);
  \draw[thick, ForestGreen]
    (-2.5,.5)--(-2.5,2) -- (1,2) -- (1,1) --(0,1)-- (0,0.5) ;
  \draw[thick, blue]
    (-2,0.5) -- (-2,4)--(-1,4)--(-1,.5);
  \node at (1,1) {$\whitecirc$};
  \node at (-2.5,2) {$\bullet$};
  \node at (0,1){$\bullet$};
  \node at (-2,4) {$\bullet$};
\end{tikzpicture}\]
Unpictured straight segments are passed using \cref{lemma:basiclocalreplacements}, and every step preserves the symplectic Demazure product. 
\end{proof}

We now show that symplectic deflation is realized by the moves of the previous subsection, which allows us to deduce that $\bpdsp(\tau)$ is connected by ($K$-theoretic) droops and symplectic droops.

\begin{lemma}
    \label{l:deflation-droop}

Let $A' \in \asmsk(2n)$ such that $(b,a),U$ is $\Sp$-addable in $A'$, and let $A$ be the result of applying symplectic deflation.
Then $\Psi(A)$ can be constructed from $\Psi(A')$ by applying a sequence of ($K$-theoretic) droops and symplectic droops.
\end{lemma}

\begin{proof}
Let $U = \{(i_1,j_1),\dots,(i_h,j_h)\} \subseteq N^\Sp_{A'}(b,a)$ with $i_1 < \dots < i_h$.
For $1 \leq m \leq h$, let $A^{(m)}$ be the result of symplectic deflation at $(b,a)$ for $U_m := \{(i_1,j_1),\dots,(i_m,j_m)\}$, so that $A^{(h)} = A$.
Since each $U_m$ is non-empty, \cref{l:symplectic-deflation} implies that $A^{(m)}$ is a skew ASM in which $(b,a)$ is $\Sp$-removable.
We show that $\Psi(A^{(1)})$ can be obtained from $\Psi(A')$ by a single droop or symplectic droop.
Then for $2 \leq m \leq h$, we show the triangular bumpless pipe dream $\Psi(A^{(m)})$ can be obtained from $\Psi(A^{(m-1)})$ by a single $K$-theoretic droop move or $K$-theoretic symplectic droop move.

By maximality, $(i_1,j_1)$ is the only non-zero entry of $A'$ in $[(i_1,j_1),(b,a)]_\preceq$.
Direct inspection shows symplectic deflation at $(b,a), \{(i_1,j_1)\}$ corresponds to performing a droop on the same interval in $\Psi(A')$ when $(i_1,j_1)$ is a lower neighbor and a symplectic droop when it is an upper neighbor.

For $m\geq2$, maximality and the ordering of $U$ show that the only nonzero entries of $A^{(m-1)}$ in $[(i_m,j_m),(b,a)]_\preceq$ are the $1$'s at $(i_m,j_m)$ and $(b,j_{m-1})$ and the $-1$ at $(b,a)$.
Passing to $A^{(m)}$ replaces these two $1$'s by $1$'s at $(i_m,j_{m-1})$ and $(b,j_m)$, together with their transposes.
This is a $K$-theoretic droop for a lower neighbor and, after folding, a $K$-theoretic symplectic droop for an upper neighbor.
Since $A=A^{(h)}$, the result follows.
\end{proof}

\begin{corollary}
\label{c:droop-graph}	
For $\tau \in \mf_{2n}$, the graph on $\bpdsp(\tau)$ induced by droop moves, symplectic droop moves, $K$-theoretic droop moves and $K$-theoretic symplectic droop moves is connected.
Moreover, every $B \in \bpdsp(\tau)$ is obtained from $\Dsp_\tau$ by a sequence of such moves.
\end{corollary}

\begin{proof}
The result is immediate from \cref{t:symplectic-equality,l:deflation-droop} and the definition of $\keysp$.
\end{proof}

\section{Flagged fixed-point-free involutions and reverse plane partitions}
\label{section:flagged}

In this section, we show that symplectic bumpless pipe dreams for certain $\mf$-involutions we call \emph{flagged} are in bijection with reverse plane partitions on shifted partitions.

\subsection{Flagged fixed-point-free involutions}
We say a triangular bumpless pipe dream is \newword{flagged} if there are no $\cross$ tiles strictly northwest of a $\nowire$ tile, and that  $\tau \in \mf_{2n}$ is \newword{flagged} if $\Dsp_\tau$ is.

The matching $\tau$ \newword{$\mf$-contains} the matching $\rho$ if there are disjoint cycles in $\tau$ whose relative order is the same as that of the cycles of $\rho$.
Otherwise, say $\tau$ \newword{$\mf$-avoids} $\rho$.
There are several characterizations for when $\tau \in \mf_{2n}$ is \newword{$\mf$-vexillary}; one, from~\cite{hamaker2020fixed}*{Cor.~7.9}, is that $\tau$ $\mf$-avoids 16 specific matchings.
See~\cite{fang2022on} for further properties of $\mf$-avoidance.

\begin{proposition}
    \label{p:flagged-vexillary}
    A matching $\tau\in\mf_{2n}$ is flagged if and only if it is $\mf$-vexillary and $\mf$-avoids the matching pattern $(1,3)(2,4)(5,7)(6,8)$.
\end{proposition}

\begin{proof}
Whether $(r,c)$ is a $\nowire$ or $\cross$ in $\Dsp_\tau$ depends only on the relative order of $r$, $c$, $\tau(r)$ and $\tau(c)$, so if $\tau$ $\mf$-contains a non-flagged pattern then $\tau$ is not flagged.
Conversely, if $\tau$ is not flagged, choose a cross at $(r,c)$ and a blank at $(s,d)$ with $r<s$ and $c<d$.
The inequalities $\tau(r)<c$, $\tau(c)<r$, $\tau(s)>d$ and $\tau(d)>s$ imply that $r,c,s,d$ lie in four distinct cycles.
Restricting to these cycles gives a non-flagged matching on eight letters.
Direct inspection shows that these are exactly the displayed pattern and the thirteen eight-letter forbidden patterns in~\cite{hamaker2020fixed}*{Cor.~7.9}.
The remaining three forbidden patterns there contain the displayed pattern, so the result follows.
\end{proof}

To each flagged $\mf$-involution $\tau$ we associate the \newword{shifted shape} $\mu^{(\tau)}$ so that $D_{\mu^{(\tau)}}$ and $\Dsp_\tau$ have the same number of cells on each diagonal. 
Here, the $d$th diagonal is $\{(i+d,i): i \in [2n-d]\}$.
The following lemma guarantees such a shifted partition exists:

\begin{lemma}
    \label{l:flagged-shape}
    For $\tau$ a flagged $\mf$-involution, the $d$th diagonal of $\Dsp_\tau$ contains either the same number of cells as the $(d{+}1)$st diagonal or one more.
\end{lemma}

\begin{proof}
Let $h_d$ be the number of $\nowire$'s on the $d$th diagonal, so the result claims $h_d - h_{d+1} \in \{0,1\}$.
To track this difference, we use an auxiliary function
\[
f_d(c) = |\{j \leq c: \tau(j) > c + d\}|.
\]
Fix $c,d$ and let $r =  c+d$.
Note that $f_d(c)$ and $f_d(c{-}1)$ differ in two possible ways;
if $\tau(c) > r$ then $f_d(c)$ grows by one, while if $\tau(r) < c$ then $f_d(c)$ shrinks by one.
Equivalently, as an equation
\[
f_d(c) - f_d(c{-1}) = 1_{\{\tau(c) > r\}} - 1_{\{\tau(r)< c\}}. 
\]
This quantity is $1$ when $\tau(c) > r, \tau(r) > c$, which occurs if and only if $\Dsp_\tau(r,c) = \nowire$, and is $-1$ when $\tau(c) < r, \tau(r) < c$, which occurs if and only if $\Dsp_\tau(r,c) = \cross$.
Since $\tau$ is flagged, all $\nowire$'s on a diagonal occur before any $\cross$, so we have $h_d = \max_c f_d(c)$.
To finish the proof, observe by definition that 
$f_d(c) - f_{d+1}(c) = 1_{\{\tau(r+1) \leq c\}} \in \{0,1\}$ for every $c$.
Taking maxima gives $h_{d+1}\le h_d\le h_{d+1}+1$.
\end{proof}

For flagged $\mf$-involutions, every symplectic bumpless pipe dream is reduced.
\begin{lemma}
    \label{l:flagged-reduced}
    For $\tau \in \mf_{2n}$ flagged, each $B \in \bpdsp(\tau)$ is flagged and reduced.
\end{lemma}

\begin{proof}
    For $A \in \asmsk(2n)$ such that $(b,a),U$ is $\Sp$-addable, note that if $|U| > 1$, then two of its entries create a cross in $\Psi(A)$ strictly northwest of $(b,a)$, which contains a $\nowire$.
    Therefore, $\Psi(A)$ cannot be flagged.

    For $B \in \bpdsp(\tau)$, by \cref{rem:deflation-sequence} there exists a symplectic deflation sequence from the alternating sign matrix $\Psi^{-1}(\Dsp_\tau)$ to $\Psi^{-1}(B)$.
    We argue by induction with base case $\Dsp_\tau$.
    By the above observation and the inductive hypothesis, every deflation from $A$ to $A'$ is at a cell $(b,a)$ below the diagonal, for some $U$ with $|U|=1$.
	This deflation corresponds to a droop or symplectic droop from $B = \Psi(A)$ to $B' = \Psi(A')$, so $B'$ is reduced as desired.
	
    By assumption $B \in \bpdsp_\mr(\tau)$ is flagged.
	Let the droop or symplectic droop from $B$ to $B'$ start at $(j,i)$ and end at $(\ell,k)$.
    Since the $(\ell,k)$-entry of $B$ is $\nowire$, the flagged condition guarantees the only matched pair of pipes in the droop region passes through $(j,i)$.
    After performing the droop this continues to be the only matched pair in the droop region so $B'$ is also flagged.
    This completes our proof.
   	\end{proof}

Next, we introduce local moves on symplectic bumpless pipe dreams for flagged $\mf$-involutions.
They come in two types: (i) \newword{(classical) local moves}
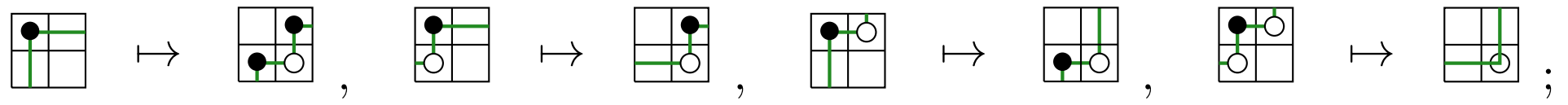
\begin{equation}
    \label{eq:classical-local}
  \begin{tikzpicture}[scale = .3]
    \draw (0,0) -- (0,2) -- (2,2) -- (2,0) -- (0,0);
    \draw (1,0) -- (1,2);
    \draw (0,1) -- (2,1);
    \draw[ForestGreen, thick] (.5,0) -- (.5,1.5) -- (2,1.5);
        \node at (.5,1.5) {$\bullet$};
\end{tikzpicture}
\quad
\raisebox{.5em}{$\mapsto$}
\quad
\begin{tikzpicture}[scale = .3]
    \draw (0,0) -- (0,2) -- (2,2) -- (2,0) -- (0,0);
    \draw (1,0) -- (1,2);
    \draw (0,1) -- (2,1);
    \draw[ForestGreen, thick] (.5,0) -- (.5,.5) -- (1.5,.5) -- (1.5,1.5) -- (2,1.5);
       \node at (.5,.5) {$\bullet$};
    \node at (1.5,1.5) {$\bullet$};
    \node at (1.5,.5) {$\whitecirc$};
\end{tikzpicture}\ ,
\quad
\begin{tikzpicture}[scale = .3]
    \draw (0,0) -- (0,2) -- (2,2) -- (2,0) -- (0,0);
    \draw (1,0) -- (1,2);
    \draw (0,1) -- (2,1);
    \draw[ForestGreen, thick] (0,.5) -- (.5,.5) -- (.5,1.5) -- (2,1.5);
        \node at (.5,1.5) {$\bullet$};
    \node at (.5,.5) {$\whitecirc$};
\end{tikzpicture}
\quad
\raisebox{.5em}{$\mapsto$}
\quad
\begin{tikzpicture}[scale = .3]
    \draw (0,0) -- (0,2) -- (2,2) -- (2,0) -- (0,0);
    \draw (1,0) -- (1,2);
    \draw (0,1) -- (2,1);
    \draw[ForestGreen, thick] (0,.5) -- (.5,.5) -- (1.5,.5) -- (1.5,1.5) -- (2,1.5);
        \node at (1.5,1.5) {$\bullet$};
    \node at (1.5,.5) {$\whitecirc$};
\end{tikzpicture}\ , \quad
\begin{tikzpicture}[scale = .3]
    \draw (0,0) -- (0,2) -- (2,2) -- (2,0) -- (0,0);
    \draw (1,0) -- (1,2);
    \draw (0,1) -- (2,1);
    \draw[ForestGreen, thick] (.5,0) -- (.5,1.5) -- (1.5,1.5) -- (1.5,2);
        \node at (.5,1.5) {$\bullet$};
        \node at (1.5,1.5) {$\whitecirc$};
\end{tikzpicture}
\quad
\raisebox{.5em}{$\mapsto$}
\quad
\begin{tikzpicture}[scale = .3]
    \draw[ForestGreen, thick] (.5,0) -- (.5,.5) -- (1.5,.5) -- (1.5,2);
        \node at (.5,.5) {$\bullet$};
    \node at (1.5,.5) {$\whitecirc$};
    \draw (0,0) -- (0,2) -- (2,2) -- (2,0) -- (0,0);
    \draw (1,0) -- (1,2);
    \draw (0,1) -- (2,1);
\end{tikzpicture}\ ,
\quad
\begin{tikzpicture}[scale = .3]
    \draw (0,0) -- (0,2) -- (2,2) -- (2,0) -- (0,0);
    \draw (1,0) -- (1,2);
    \draw (0,1) -- (2,1);
    \draw[ForestGreen, thick] (0,.5) -- (.5,.5) -- (.5,1.5) -- (1.5,1.5) -- (1.5,2);
        \node at (.5,1.5) {$\bullet$};
    \node at (.5,.5) {$\whitecirc$};
    \node at (1.5,1.5) {$\whitecirc$};
\end{tikzpicture}
\quad
\raisebox{.5em}{$\mapsto$}
\quad
\begin{tikzpicture}[scale = .3]
    \node at (1.5,.5) {$\whitecirc$};
    \draw (0,0) -- (0,2) -- (2,2) -- (2,0) -- (0,0);
    \draw (1,0) -- (1,2);
    \draw (0,1) -- (2,1);
    \draw[ForestGreen, thick] (0,.5) -- (.5,.5) -- (1.5,.5) -- (1.5,2);
\end{tikzpicture}\ ;
\end{equation}
and (ii) \newword{symplectic local moves}:
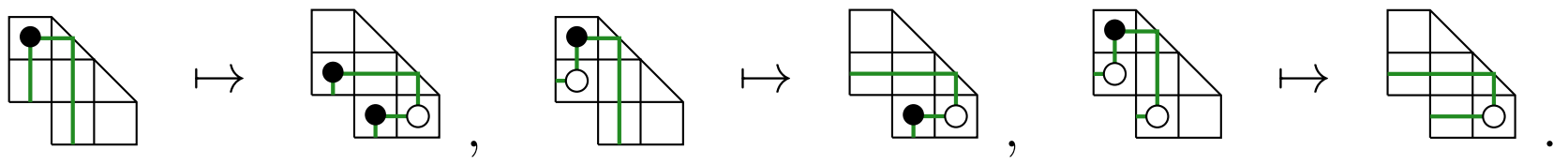
\begin{equation}
    \label{eq:symplectic-local}
\begin{tikzpicture}[scale=.3]
\draw (0,1) -- (0,3) -- (1,3) -- (3,1) -- (3,0) -- (1,0);
\draw (0,2) -- (2,2) -- (2,0);
\draw (0,1) -- (3,1);
\draw (1,0) -- (1,3);
\draw[ForestGreen,thick] (.5,1) -- (.5,2.5) -- (1.5,2.5) -- (1.5,0);
    \node at (.5,2.5) {$\bullet$};
\end{tikzpicture}
\quad
\raisebox{1em}{$\mapsto$}
\quad
    \begin{tikzpicture}[scale=.3]
\draw (0,1) -- (0,3) -- (1,3) -- (3,1) -- (3,0) -- (1,0);
\draw (0,2) -- (2,2) -- (2,0);
\draw (0,1) -- (3,1);
\draw (1,0) -- (1,3);
\draw[ForestGreen,thick] (.5,1) -- (.5,1.5) -- (2.5,1.5) -- (2.5,.5) -- (1.5,.5) -- (1.5,0);
\node at (.5,1.5) {$\bullet$};
\node at (1.5,.5) {$\bullet$};
\node at (2.5,.5) {$\whitecirc$};
\end{tikzpicture}\ , \quad 
\begin{tikzpicture}[scale=.3]
\draw (0,1) -- (0,3) -- (1,3) -- (3,1) -- (3,0) -- (1,0);
\draw (0,2) -- (2,2) -- (2,0);
\draw (0,1) -- (3,1);
\draw (1,0) -- (1,3);
\draw[ForestGreen,thick] (0,1.5) -- (.5,1.5) -- (.5,2.5) -- (1.5,2.5) -- (1.5,0);
\node at (.5,2.5) {$\bullet$};
\node at (.5,1.5) {$\whitecirc$};
\end{tikzpicture}
\quad
\raisebox{1em}{$\mapsto$}
\quad
    \begin{tikzpicture}[scale=.3]
\draw (0,1) -- (0,3) -- (1,3) -- (3,1) -- (3,0) -- (1,0);
\draw (0,2) -- (2,2) -- (2,0);
\draw (0,1) -- (3,1);
\draw (1,0) -- (1,3);
\draw[ForestGreen,thick] (0,1.5) -- (2.5,1.5) -- (2.5,.5) -- (1.5,.5) -- (1.5,0);
\node at (1.5,.5) {$\bullet$};
\node at (2.5,.5) {$\whitecirc$};
\end{tikzpicture}\ , \quad 
\begin{tikzpicture}[scale=.3]
\draw (0,1) -- (0,3) -- (1,3) -- (3,1) -- (3,0) -- (1,0);
\draw (0,2) -- (2,2) -- (2,0);
\draw (0,1) -- (3,1);
\draw (1,0) -- (1,3);
\draw[ForestGreen,thick] (0,1.5) -- (.5,1.5) -- (.5,2.5) -- (1.5,2.5) -- (1.5,.5) -- (1,.5);
\node at (.5,2.5) {$\bullet$};
\node at (.5,1.5) {$\whitecirc$};
\node at (1.5,.5) {$\whitecirc$};
\end{tikzpicture}
\quad
\raisebox{1em}{$\mapsto$}
\quad
    \begin{tikzpicture}[scale=.3]
\draw (0,1) -- (0,3) -- (1,3) -- (3,1) -- (3,0) -- (1,0);
\draw (0,2) -- (2,2) -- (2,0);
\draw (0,1) -- (3,1);
\draw (1,0) -- (1,3);
\draw[ForestGreen,thick] (0,1.5) -- (2.5,1.5) -- (2.5,.5) -- (1,.5);
\node at (2.5,.5) {$\whitecirc$};
\end{tikzpicture}\ 
.
\end{equation}
Note the configuration where the pipe enters the first column from below and exits the second to the left cannot occur, as these diagrams would have a $\cross$ between two pipes that meet on the diagonal.
For $B,B'$ differing by a local move it is clear by construction that $\demprodfpf(B) = \demprodfpf(B')$.

\begin{proposition}
    \label{p:local-flagged}
    For $\tau$ a flagged $\mf$-involution, $\bpdsp(\tau)$ is connected by local moves and symplectic local moves.
\end{proposition}

\begin{proof}
    By \cref{l:flagged-reduced,c:droop-graph}, it suffices to decompose each droop and symplectic droop into local moves.
    If the move is of Type I, it decomposes row by row from top to bottom into a composition of local moves.
    If it is of Type II, one can apply local moves in the highest row until one can apply a symplectic local move (this may require zero local moves), then apply local moves in the rightmost column until the corresponding $\nwelbow$ is in the bottom row of the droop region.
\end{proof}

\subsection{Reverse plane partitions} The \newword{subdiagonal} of $\lrtri_{2n}$ is $\{(i+1,i):i\in [2n-1]\}$, the first diagonal.
A \newword{reverse plane partition} on $D \subseteq \lrtri_{2n}$ is a function $f:D \to \mathbb{N}$ so that $(b,a) \preceq (d,c)$ implies $f((b,a)) \leq f((d,c))$.
Let $\RPP(D)$ be the set of reverse plane partitions on $D$ and $\RPP^E(D)$ be the subset whose subdiagonal entries are all even.
For $P \in \RPP(D)$ let 
\[
\RPP_P(D) = \{Q \in \RPP(D): Q_{ji} \leq P_{ji} \text{ for all } (j,i) \in D\},\] 
and \[ \RPP^E_P(D) = \RPP_P(D) \cap \RPP^E(D).
\]
As a shorthand, define $\RPP_k(D) = \RPP_{P^{(k)}}(D)$ where $P^{(k)}_{ji} = k$ for all $(j,i) \in D$, and likewise let $\RPP_k^E(D)=\RPP_{P^{(k)}}^E(D)$.
We are particularly interested in the case where $D = D_\lambda$ is the shifted Young diagram of a strict partition $\lambda$ (in our transposed convention).
For example,
\begin{equation}
\label{eq:rpp_ex}
P = \begin{ytableau}
0\\
0 & 1\\
1 & 1 & 4\\
2 & 2
\end{ytableau}\ ,
\quad 
Q = \begin{ytableau}
0\\
0 & 2\\
1 & 2 & 4\\
2 & 3
\end{ytableau},
\quad \mbox{and} \quad
R = \begin{ytableau}
0\\
0 & 2\\
1 & 2 & 2\\
1 & 2
\end{ytableau}
\end{equation}
are all in $\RPP_4(D_{431})$ and $P,R \in \RPP_Q(D_{431})$, but $P \notin \RPP^E_4(D_{431})$ since $P_{32}$ is an odd subdiagonal entry.
Additionally, $R \in \RPP^E_2(D_{431})$ while $P$ and $Q$ are not.

There is a standard construction describing reverse plane partitions as a collection of non-intersecting lattice paths.
Given a reverse plane partition of shifted shape, the $i$th path starts from the bottom left, separating the cells marked $i-1$ and $i$, terminating on the subdiagonal.
To make the lattice paths non-overlapping, one can shift the $i$th path and all cells labeled $i$ down and to the right $i$ units for all $i$.
For example, with $R$ as in \cref{eq:rpp_ex}, construct $R'$ from $R$ by setting $R'_{32} = 1$ (this is the second entry along the subdiagonal).
After performing the shifts to separate lattice paths we see $R$ and $R'$ correspond to
	\[
\ytableausetup{smalltableaux, centertableaux}
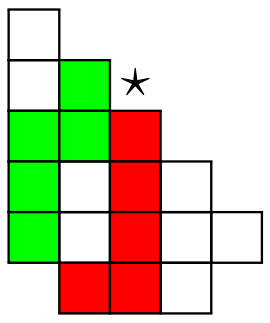
	\begin{ytableau}
		\\
		& *(green) & \none[\star]\\
		*(green) & *(green) &*(red)\\
		*(green) & & *(red) & \\
		*(green) & & *(red) & &\\
		\none & *(red) & *(red)&
	\end{ytableau} \quad \mbox{and} \quad 
	
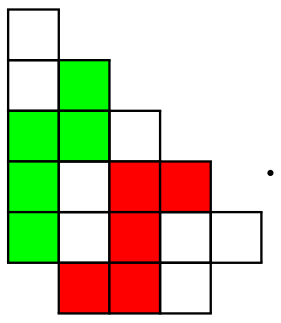
\begin{ytableau}
		\\
		& *(green)\\
		*(green) & *(green) &\\
		*(green) & & *(red) & *(red)\\
		*(green) & & *(red) & &\\
		\none & *(red) & *(red)&
	\end{ytableau}.
	\]
Here the green path separates values $0$ and $1$ while the red path separates values $1$ and $2$.
For $R$, we can connect the paths along the diagonal at $\star$, but this is not possible for $R'$.
More generally, observe that paths surrounding contiguous odd entries on the subdiagonal will not meet on the diagonal, but they will meet between even entries with different values.

For $P \in \RPP(D)$, the \newword{weight} of $P$ is
\begin{equation}
    \label{eq:rpp-weight}
    \wt_\RPP(P) = \prod_{(j,i) \in D} (x_{P_{ji}+j} + x_{P_{ji}+i}).
\end{equation}

\begin{figure}
	\begin{center}
	\begin{tikzpicture}[scale=.5]
\draw (.5,11.5) -- (.5,-.5) -- (12.5,-.5) -- (.5,11.5);
\draw (.5,10.5) -- (1.5,10.5) -- (1.5,-.5);
\draw (.5,9.5) -- (2.5,9.5) -- (2.5,-.5);
\draw (.5,8.5) -- (3.5,8.5) -- (3.5,-.5);
\draw (.5,7.5) -- (4.5,7.5) -- (4.5,-.5);
\draw (.5,6.5) -- (5.5,6.5) -- (5.5,-.5);
\draw (.5,5.5) -- (6.5,5.5) -- (6.5,-.5);
\draw (.5,4.5) -- (7.5,4.5) -- (7.5,-.5);
\draw (.5,3.5) -- (8.5,3.5) -- (8.5,-.5);
\draw (.5,2.5) -- (9.5,2.5) -- (9.5,-.5);
\draw (.5,1.5) -- (10.5,1.5) -- (10.5,-.5);
\draw (.5,0.5) -- (11.5,0.5) -- (11.5,-.5);
\draw[ForestGreen,thick] (1,-0.500000000000000) -- (1,8);
\draw[ForestGreen,thick] (1,8) -- (2,8);
\draw[ForestGreen,thick] (2,8) -- (2,9);
\draw[ForestGreen,thick] (2,9) -- (3,9);
\draw[ForestGreen,thick] (2,-0.500000000000000) -- (2,6);
\draw[ForestGreen,thick] (2,6) -- (3,6);
\draw[ForestGreen,thick] (3,6) -- (3,9);
\draw[ForestGreen,thick] (3,-0.500000000000000) -- (3,4);
\draw[ForestGreen,thick] (3,4) -- (4,4);
\draw[ForestGreen,thick] (4,4) -- (4,5);
\draw[ForestGreen,thick] (4,5) -- (5,5);
\draw[ForestGreen,thick] (5,5) -- (5,6);
\draw[ForestGreen,thick] (5,6) -- (6,6);
\draw[ForestGreen,thick] (4,-0.500000000000000) -- (4,3);
\draw[ForestGreen,thick] (4,3) -- (6,3);
\draw[ForestGreen,thick] (6,3) -- (6,6);
\draw[ForestGreen,thick] (5,-0.500000000000000) -- (5,2);
\draw[ForestGreen,thick] (5,2) -- (10,2);
\draw[ForestGreen,thick] (6,-0.500000000000000) -- (6,1);
\draw[ForestGreen,thick] (6,1) -- (11,1);
\draw[ForestGreen,thick] (7,-0.500000000000000) -- (7,3);
\draw[ForestGreen,thick] (7,3) -- (9,3);
\draw[ForestGreen,thick] (8,-0.500000000000000) -- (8,0);
\draw[ForestGreen,thick] (8,0) -- (12,0);
\draw[ForestGreen,thick] (9,-0.500000000000000) -- (9,3);
\draw[ForestGreen,thick] (10,-0.500000000000000) -- (10,2);
\draw[ForestGreen,thick] (11,-0.500000000000000) -- (11,1);
\draw[ForestGreen,thick] (12,-0.500000000000000) -- (12,0);
\node at (1,8) {$\bullet$};
\node at (2,9) {$\bullet$};
\node at (2,6) {$\bullet$};
\node at (3,4) {$\bullet$};
\node at (4,5) {$\bullet$};
\node at (4,3) {$\bullet$};
\node at (5,6) {$\bullet$};
\node at (5,2) {$\bullet$};
\node at (6,1) {$\bullet$};
\node at (7,3) {$\bullet$};
\node at (8,0) {$\bullet$};
\node at (2,8) {$\whitecirc$};
\node at (3,6) {$\whitecirc$};
\node at (4,4) {$\whitecirc$};
\node at (5,5) {$\whitecirc$};
\node at (6,3) {$\whitecirc$};
\node at (1,-1) {1};
\node at (2,-1) {2};
\node at (3,-1) {3};
\node at (4,-1) {4};
\node at (5,-1) {5};
\node at (6,-1) {6};
\node at (7,-1) {7};
\node at (8,-1) {8};
\node at (9,-1) {9};
\node at (10,-1) {10};
\node at (11,-1) {11};
\node at (12,-1) {12};
\node at (1,10) {$0$};
\node at (1,9) {$0$};
\node at (2,7) {$1$};
\node at (3,5) {$2$};
\node at (4,6) {$2$};
\node at (4,7) {$2$};
\node at (5,4) {$3$};
\node at (7,4) {$4$};
\draw[dashed,->] (1.8,7.2) -- (1.2,7.8);
\draw[dashed,->] (2.8,5.2) -- (1.2,6.8);
\draw[dashed,->] (3.8,6.2) -- (2.2,7.8);
\draw[dashed,->] (3.8,7.2) -- (2.2,8.8);
\draw[dashed,->] (4.8,4.2) -- (2.2,6.8);
\draw[dashed,->] (6.8,4.2) -- (3.2,7.8);
\end{tikzpicture}
\end{center}
\caption{For $\tau =(1,2)(3,4)(5,10)(6,11)(7,9)(8,12) = \Theta_4 \oplus \tau(431)$, we depict a symplectic bumpless pipe dream $B \in \bpdsp(\tau)$.
Each blank cell in $B$ is labeled with the number of pipes to its northwest.
Shifting these cells up the dashed lines, we see that $\Phi_\tau(B) = Q$ from~\cref{eq:rpp_ex} 
and can confirm $\wt_\fkS(B) =\wt_\RPP(Q)$.
}
\label{fig:bpd=rpp}
\end{figure}
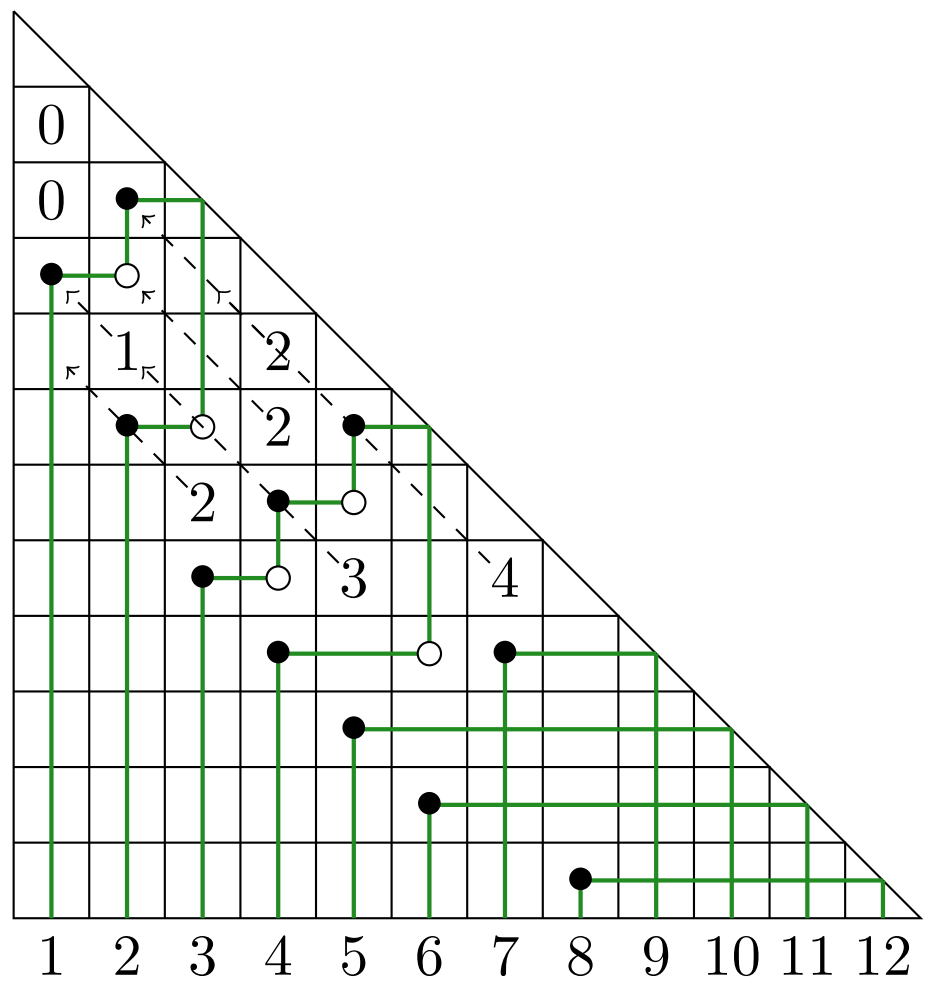

For $\tau$ flagged and $B \in \bpdsp(\tau)$, by \cref{p:local-flagged} the number of blank cells on the $d$th diagonal is the same as for $D_{\mu^{(\tau)}}$.
We define a map $\Phi_\tau:\bpdsp(\tau) \to \RPP(D_{\mu^{(\tau)}})$ as follows.
For $B \in \bpdsp(\tau)$ and $(j,i) \in D_{\mu^{(\tau)}}$ with $d = j-i$, let $(k+d,k)$ contain the $i$th $\nowire$ on diagonal $d$ of $B$.
We construct $P = \Phi_\tau(B)$ by setting $P_{ji} = k-i$, which is the number of pipes northwest of $(k+d,k)$.
See \cref{fig:bpd=rpp} for an example.
Additionally, define $P^\tau = \Phi_\tau(\Dsp_\tau)$ and $\RPP^E(\tau) = \RPP^E_{P^\tau}(D_{\mu^{(\tau)}})$.

\begin{theorem}
    \label{t:flagged-rpp}
    For $\tau \in \mf_{2n}$ flagged, the map $\Phi_\tau:\bpdsp(\tau) \to \RPP^E(\tau)$ is a bijection.  Furthermore, $\wt_\fkS(B) = \wt_\RPP(\Phi_\tau(B))$.
\end{theorem}

\begin{proof}
    We first show that $\Phi_\tau$ is well defined by showing that local moves correspond to simple transformations of shifted reverse plane partitions.
    Since the blank positions on consecutive diagonals interlace in $\Dsp_\tau$ by the proof of \cref{l:flagged-shape},  their diagonal displacements weakly increase along rows and columns.
    Therefore \(P^\tau\) is a reverse plane partition.
    To see $P^\tau \in \RPP^E_{P^\tau}(D_{\mu^{(\tau)}})$, observe pipes occupy consecutive entries of its subdiagonal.
    
    For $B \in \bpdsp(\tau)$, assume $\Phi_\tau(B)$ is defined.
    For $B'$ differing from $B$ by a classical local move, it is easy to see $\Phi_\tau(B')$ differs from $\Phi_\tau(B)$ in one entry that is not on the subdiagonal, which is decremented by one.
    Similarly, a symplectic local move decrements one subdiagonal entry by two.
    Moreover, any cell immediately above or to the left of the decremented tile must be to the northwest of one (classical) or two (symplectic) additional pipes, hence $\Phi_\tau(B')$ is also a reverse plane partition.
    Since there exists a sequence of (symplectic) local moves transforming $\Dsp_\tau$ to $B$ by \cref{p:local-flagged} with $P^\tau$ as a base case, $\Phi_\tau$ is well defined.

    To see $\Phi_\tau$ is a bijection, let $P \in \RPP^E(\tau)$.
    One can construct $P$ from $P^\tau$ by decrementing the entry in each cell a sufficient number of times, scanning row by row from top to bottom, left to right.
    We claim each such decrement corresponds to a local move.
    To see this, let $P$ be a reverse plane partition that will be decremented in cell $(j,i)$, and assume $P$ corresponds to $B \in \bpdsp(\tau)$ with the $\nowire$ in cell $(\ell,k)$ of $B$ corresponding to $(j,i) \in D_{\mu^{(\tau)}}$.
    Note the entries $(\ell{-1},k)$ and $(\ell,k{-1})$ in $B$ have the same number of pipes weakly northwest as $(\ell,k)$ does.
    Therefore, they cannot be $\nowire$, else they would correspond to cells $(j{-1},i)$ or $(j,i{-1})$ in $D_{\mu^{(\tau)}}$ with the same labels, contradicting the fact that $(j,i)$ can be decremented.
    Since $B$ is flagged, this implies a local move can be performed corresponding to decrementing $P_{ji}$ \textbf{unless} $(j,i)$ is on the subdiagonal.
    
    In the case where $(j,i)$ is on the subdiagonal, $(\ell,k)$ is also on the subdiagonal.
    Since $B(\ell,k{-1}) \neq \nowire$, it must be $\vwire$ or $\nwelbow$.
    The flagged condition and the fact $B(\ell,k) = \nowire$ combine to imply $B(\ell{-1},k{-1}) = \vwire$.
    If $B(\ell{-1},k{-2}) = \nowire$, we have $P_{j\,i{-1}} = P_{ji} -1$, so $P_{ji}$ cannot be decremented while preserving the even subdiagonal property.
    Since $B(\ell{-1},k{-2})$ cannot exit to its right, we see $B(\ell{-1},k{-2}) = \vwire$ or $\nwelbow$.
    For the pipes to connect, we must have $B(\ell{-2},k{-2}) = \seelbow$, hence the decrement can be performed via a symplectic local move.
    
    Note that $\cross$ tiles are preserved by (symplectic) local moves, so they are the same for every $B \in \bpdsp(\tau)$.
    Injectivity follows since $B$ is determined by its $\cross$ tiles and its $\nowire$ tiles, which are of the form $(j+P_{ji},i+P_{ji})$ for $(j,i) \in D_{\mu^{(\tau)}}$.
    This also shows the bijection is weight preserving.
\end{proof}

Recall $\Theta_{2n} = 2143\cdots 2n\, 2n{-1} = (1,2)(3,4) \dots (2n{-}1,2n)$.
For $\pi \in S_m$ and $\sigma \in S_n$, let $\pi \oplus \sigma = \pi_1 \cdots \pi_m\ (\sigma_1 {+} m) \cdots (\sigma_n {+}m)$.
For $\lambda$ a shifted shape, let $\tau(\lambda)$ be the $\mf$-dominant matching so that $\Dsp_{\tau(\lambda)} = D_{\lambda}$ as a set.

\begin{corollary}
    \label{c:rpp}
    For $\lambda$ a shifted shape, and $\tau = \Theta_{2k} \oplus \tau(\lambda)$,
    \[
    \fkSsp_\tau = \sum_{P \in \RPP^E_{2k}(D_\lambda)} \wt_\RPP(P).
    \]
\end{corollary}

\begin{proof}
    Since the first $2k$ pipes in $\Dsp_\tau$ never cross and these are the only pipes appearing weakly northwest of the $\nowire$'s in $\Dsp_\tau$, we see $\tau$ is flagged.
    Moreover, all of these pipes appear weakly northwest of each $\nowire$ in $\Dsp_\tau$, so $\Phi_\tau(\Dsp_\tau) = P^{(2k)}$.
    The result now follows from \cref{t:flagged-rpp}.
\end{proof}

Combining \cref{c:rpp} with~\cite{okada2021intermediate}*{Thm.~5.1}, we obtain the following enumerative corollary, which answers the $\mf$ case of~\cite{hamaker2022involution}*{Prob.~6.4}.

\begin{corollary}
	\label{c:rpp_count}
For $\lambda$ a shifted shape,
$
\fkSsp_{\Theta_{2k} \oplus \tau(\lambda)}(1/2,1/2,\dots,1/2)  = |\RPP^E_{2k}(D_\lambda)|.
$
In particular,
\[
\fkSsp_{\Theta_{2k} \oplus 2n\dots21}(1/2,1/2,\dots,1/2) = \prod_{i,j \in [n], i \neq j} \frac{i+j +2k-1}{i+j-1}.
\]
\end{corollary}

\begin{proof}
	For $P \in \RPP(D)$, note that $\wt_\RPP(P)(1,\dots,1) = 2^{|D|}$.
    Then by \cref{c:rpp},
\[
\fkSsp_\tau(1,1,\dots,1) = 2^{|\mu^{(\tau)}|} |\RPP^E_{2k}(D_\lambda)|.
\]
Since $\fkSsp_\tau$ is homogeneous of degree $\ellfpf(\tau) = |\mu^{(\tau)}|$, the first result follows.
The second is an immediate consequence of the fact $\tau(\lambda) = 2n\dots21 \in \mf_{2n}$ for $\lambda = (2n-2,2n-4,\dots,2)$ and~\cite{okada2021intermediate}*{Thm.~5.1}.
\end{proof}

The product formula for $\RPP^E_{2k}(D_{(2n{-2},\dots,4,2)})$ extends to shapes of the form $\lambda = (n{-}1,\dots,1) + (m{-1},\dots,1)$ with addition applied pointwise.
Therefore, \cref{c:rpp_count} specializes to a product formula for these cases as well.

\subsection{Set-valued reverse plane partitions}

Using marked symplectic bumpless pipe dreams, we can extend \cref{t:flagged-rpp} to a formula for $\fkGsp_\tau$ when $\tau$ is flagged.
To do so, we use \newword{set-valued reverse plane partitions} on $D \subseteq \lrtri_{2n}$, which are functions $f$ from $D$ to the non-empty finite subsets of $\mathbb{N}$ so that $(b,a) \prec (d,c)$ implies $\max f((b,a)) \leq \min f((d,c))$.
Let $\SetRPP(D)$ be the set of set-valued reverse plane partitions on $D$, and for $P \in \RPP(D)$ extend the definitions $\SetRPP_P(D)$, $\SetRPP^E(D)$, $\SetRPP^E_P(D)$ and $\SetRPP_k(D)$ analogously with restrictions now applied to the entire set of each entry.
For $P \in \SetRPP(D)$, let $\min(P) \in \RPP(D)$ be the reverse plane partition obtained by taking the minimum of each set.
Similarly, define $\SetRPP^E(\tau) := \SetRPP^E_{P^\tau}(D_{\mu^{(\tau)}})$.

For $P \in \SetRPP(D)$, define
\[
\wt_\SetRPP(P) = \prod_{(j,i) \in D} \prod_{k \in P_{ji}} \beta(x_{k +i} \oplus x_{k+j}).
\]
We now show for $\tau$ flagged that \cref{t:flagged-rpp} can be extended to the set-valued setting by defining $\Phi_\tau:\bpdsp_\mm(\tau)\rightarrow \SetRPP(D_{\mu^{(\tau)}})$.
To construct $P'=\Phi_\tau(B,M)$, first construct $P = \Phi_\tau(B)$.
Next, for each $(\ell,k) \in M$, find the nearest $\nowire$ northwest of $(\ell,k)$ on the same diagonal, let $(j,i) \in D_{\mu^{(\tau)}}$ be the cell corresponding to it, and add $k-i$ to $P_{ji}$.
For example, if each $\nwelbow$ in \cref{fig:bpd=rpp} is marked, we obtain
\[\ytableausetup{nosmalltableaux,centertableaux}
\begin{ytableau}
    0\\
    \scriptstyle {0,1} & 2\\
    \scriptstyle{1,2} & \scriptstyle{2,3} & 4  \\
    \scriptstyle{2,3} & \scriptstyle{{3,4}}
\end{ytableau}.
\]

\begin{theorem}
    \label{t:flagged-setrpp}
    For $\tau \in \mf_{2n}$ flagged, 
    \[
    \fkGsp_\tau = \beta^{-\ellfpf(\tau)} \sum_{P \in \SetRPP^E(\tau)} \wt_{\SetRPP}(P).
    \]
\end{theorem}

\begin{proof}
By \cref{t:flagged-rpp,c:marked}, it suffices to show that the
$\nwelbow$ tiles assigned to each cell of $P=\Phi_\tau(B)$
correspond to the allowable extra entries there, with weights
preserved.
We induct along a sequence of local moves from $\Dsp_\tau$ to $B$,
as in \cref{p:local-flagged}.
The base case has neither $\nwelbow$ tiles nor allowable extra entries.

A classical move decreasing $P_{ji}$ from $t$ to $t-1$ allows
$t$ as an extra entry at $(j,i)$ and removes it, if previously
allowed, from $(j-1,i)$ and $(j,i-1)$.
These are exactly the changes to the $\nwelbow$ tiles
in~\cref{eq:classical-local}.
For a symplectic move, $j=i+1$ and $t$ is even.
Decreasing $P_{ji}$ from $t$ to $t-2$ allows $t$ as an extra
entry at $(j,i)$ and removes $t$ and $t-1$, if previously
allowed, from $(j,i-1)$.
Again these are exactly the changes to the $\nwelbow$ tiles
in~\cref{eq:symplectic-local}.
All other assignments are unchanged, completing the induction.

Since allowable extra entries can be chosen independently,
markings of $B$ correspond bijectively to fillings with minimum $P$.
Finally, a minimum or extra entry $t$ at $(j,i)$ corresponds
to a blank or marked elbow at $(j+t,i+t)$, respectively,
with weight $\beta(x_{i+t}\oplus x_{j+t})$.
Thus the bijection preserves weights.
\end{proof}

\section{Final remarks}
\label{section:future}

We view the present paper as a launching point for many problems of interest.
Here are some we find particularly promising.

\subsection{Relation to shifted Lascoux polynomials}
In~\cite{marberg2023key}, Marberg and Scrimshaw introduce shifted key and Lascoux polynomials indexed by certain weak compositions.
They conjecture that for each $\tau \in \mf_{2n}$, $\fkGsp_\tau$ is a positive linear combination of shifted Lascoux polynomials.
In particular, when $\tau$ is $\mf$-vexillary (see~\cite{hamaker2020fixed}) this conjecture would imply $\fkGsp_\tau$ is a single shifted Lascoux polynomial.
In light of this conjecture, \cref{t:flagged-setrpp} should give a tableau formula for the shifted Lascoux polynomials associated to flagged $\mf$-involutions.
It would be interesting to characterize the weak compositions for which such a formula is possible and to determine how to write down the weak composition associated to flagged $\tau$.

We remark that identifying a tableau formula even for $\mf$-vexillary $\tau$ is more complicated than \cref{c:rpp} might suggest.
To see why, consider $\tau = (1,3)(2,4)(5,7)(6,8)$, which is $\mf$-vexillary but not flagged.
By \cref{c:schubert-bpd},
\[
\fkSsp_\tau = (x_1+x_2)\left((x_5+x_6)+(x_2+x_4)+(x_1+x_3) \right).
\]
The terms with $(x_2+x_4)$ and $(x_1+x_3)$ can be interpreted as coming from reverse plane partitions, but it is far from clear how the term with $(x_5+x_6)$ can be put on similar footing.

\subsection{Bijections for symplectic bumpless pipe dreams}
In~\cite{huang2023bijective}, Huang gives a bijective proof of Monk's rule for double Schubert polynomials using bumpless pipe dreams.
In our setting, this corresponds to the cohomology specialization of \cref{t:g_monk}.
It is natural to ask:
\begin{problem}
\label{prob:monk}
Find a bijective proof of \cref{t:g_monk} using symplectic bumpless pipe dreams.
\end{problem}
As discussed in \cref{rem:asm-tree}, symplectic bumpless pipe dreams are a bijective realization of the symplectic LS tree equations.
This can be viewed as a special case of \cref{prob:monk}. A complete solution would be significantly more complicated, and the $\beta=0$ case for reduced symplectic bumpless pipe dreams would still be of great interest.

Another natural bijective question is:
\begin{problem}
    \label{prob:pipe-dream}
    Construct a natural bijection between symplectic Hecke pipe dreams and marked symplectic bumpless pipe dreams for $\tau \in \mf_{2n}$.
\end{problem}
Building on a preliminary draft of our work, Pang proposed a conjectural solution to the reduced case of this problem in his MPhil thesis~\cite{pang2022bijections}*{Conj.~5.2.3}.

For pipe dreams, there is a bijective proof of Monk's rule for single Schubert polynomials based on the Little bijection~\cite{little2003combinatorial}, with details fully worked out in~\cite{billey2019bijective}.
By applying this bijection for the LS tree recurrence side by side with inflation, one can construct a recursive bijection, but this implicitly depends on a choice of recurrences.
It is a striking result that there is an intrinsic description of this map independent of the choice of recurrences~\cite{gao2021canonical}.
A similar Little bijection for $\mf$-involutions appears in~\cite{hamaker2018transition}, but the extension to symplectic pipe dreams has not been worked out.
This would give a recursive solution to \cref{prob:pipe-dream} for $\beta=0$, but an optimal solution would be independent of the choice of recurrences used.

\subsection{Back stable Schubert calculus}
The original motivation for introducing bumpless pipe dreams was to give a combinatorial formula for \emph{back stable Schubert polynomials}~\cite{lam2021back}, which represent classes in back stable Schubert calculus.
We now introduce this setting.
For $f:\Z\to\Z$, let $\operatorname{supp}(f)=\{i:f(i)\neq i\}$, and let $S_\Z$ be the permutations of $\Z$ with finite support, generated by $s_i=(i,i{+}1)$.
Define $\iota(s_i)=s_{i+1}$ and $\gamma(x_i)=x_{i+1}$, with $\gamma(\beta)=\beta$.
The (single) \newword{back stable Grothendieck polynomial} of $w\in S_\Z$ is
\[
\overleftarrow{\fkG}_w
=\lim_{n\to\infty}\gamma^{-n}\fkG_{\iota^n(w)}
=\lim_{n\to\infty}\gamma^{-2n}\fkG_{\iota^{2n}(w)}.
\]
These limits exist coefficientwise~\cite{lam2023back}: each fixed monomial in the $x_i$ has an eventually constant coefficient in $\Z[\beta]$.
The resulting formal series use variable indices bounded above, but may have infinitely many terms in a fixed degree.

Pawlowski introduced back stable symplectic Schubert polynomials $\overleftarrow{\fkS}^\Sp_\tau$ to represent cohomology classes of certain $\mathrm{GL}_n$-orbit closures in the product of an orthogonal Grassmannian and a flag variety~\cite{pawlowski2019universal}*{Thm.~5.11}.
We sketch how to extend our formula to this setting.
Let $\Theta_\Z=\cdots(-3,-2)(-1,0)(1,2)(3,4)\cdots$, and define
\[
\mf_\Z=\{\Theta_\Z\iprod w:w\in S_\Z\}\setminus\{\varnothing\},
\qquad
\mb(\tau)=\{w\in S_\Z:\Theta_\Z\iprod w=\tau\}.
\]
The set $\mb(\tau)$ is finite.
Indeed, Bruhat rank comparison with $\tau$ forces every intermediate matching in a word for $\tau$ to retain the adjacent pairs outside a fixed finite interval with odd left endpoint and even right endpoint containing every position where $\tau$ differs from $\Theta_\Z$; a generator outside that interval would either break one of these pairs across a preserved cut or give the null matching.
Thus all such words use generators in a fixed finite symmetric group.

We set $\ellfpf(\tau)=\min\{\ell(w):w\in\mb(\tau)\}$, which agrees with the rank function when $\tau\in\mf_{2n}$.
Define the \newword{back stable symplectic Grothendieck polynomial} by
\[
\overleftarrow{\fkG}^\Sp_\tau
=\beta^{-\ellfpf(\tau)}\sum_{w\in\mb(\tau)}
       \beta^{\ell(w)}\overleftarrow{\fkG}_w.
\]
Setting $\beta=0$ recovers Pawlowski's definition; the finite version is the usual expansion of $\fkGsp_\tau$ into ordinary Grothendieck polynomials~\cite{marberg2020k}*{Thm.~3.12}.

Extend the diagram of $\Dsp_\tau$ to $\lrtri_\Z=\{(j,i)\in\Z^2:i<j\}$  with adjacent trivial pipes connecting $-i,-i{+1}$ for $i$ odd.
In particular, the northwest tail consists of pipes of the form
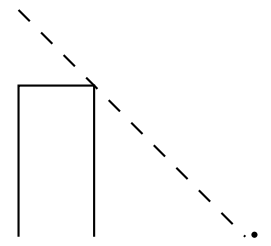
\[
\begin{tikzpicture}
\draw[dashed] (.5,-.5) -- (2,-2);
\draw (.5,-2) -- (.5,-1) -- (1,-1) -- (1,-2);
\end{tikzpicture}.
\]
Let $\bpdspback(\tau)$ be the set obtained from this extension by a finite number of ($K$-theoretic) droops and symplectic droops.
Its diagrams have finitely many $\nowire$ and $\nwelbow$ tiles, so $\wt$ is defined by the same finite products as before.

\begin{proposition}
\label{p:backstable}
For $\tau\in\mf_\Z$, the following sum is coefficientwise well defined and satisfies
\[
\overleftarrow{\fkG}^\Sp_\tau
=\beta^{-\ellfpf(\tau)}\sum_{B\in\bpdspback(\tau)}\wt(B).
\]
\end{proposition}

\begin{proof}
Fix an even $N$ beyond every position where $\tau$ differs from $\Theta_\Z$.
For sufficiently negative even $m$, let $\bpdspback_m(\tau)$ be the diagrams obtained by extending a finite diagram on the labels $[m+1,N]$ with adjacent trivial pipes.
Shifting labels by $-m$ identifies this set with $\bpdsp(\tau^{(m)})$, where $\tau^{(m)}$ is the restriction of $\tau$ to $[m+1,N]$, shifted to $[1,N-m]$.
By \cref{c:sbpd-formula},
\[
\beta^{-\ellfpf(\tau)}
 \sum_{B\in\bpdspback_m(\tau)}\wt(B)
 =\gamma^m\fkGsp_{\tau^{(m)}}.
\]
The finite ordinary Grothendieck expansion, followed by ordinary back stability, shows that the right-hand side tends coefficientwise to $\overleftarrow{\fkG}^\Sp_\tau$ as $m\to-\infty$; the indexing set $\mb(\tau)$ is finite and independent of $m$.
On the left, adjoining an adjacent pair corresponds to $A\mapsto\Theta_2\oplus A$ on ASMs, introducing no weighted tile.
Hence these diagram sets are nested with unchanged weights.
Their expanded weights have nonnegative coefficients, so coefficientwise stabilization also proves that only finitely many diagrams contribute to any fixed monomial, allowing us to take their union.
Adjoining pairs at the southeast end likewise preserves weights; the finite formula and ordinary forward stability show that increasing $N$ adds no diagrams beyond these extensions.
The union is therefore $\bpdspback(\tau)$, giving the result.
\end{proof}

\bibliographystyle{amsrn}
\bibliography{deriv}

\end{document}